\documentclass[12pt]{article}
\usepackage{latexsym}
\usepackage{amsfonts}
\usepackage{amssymb}
\usepackage{amsmath}
\newtheorem{theorem}{Theorem}
\newtheorem{corollary}{Corollary}
\newtheorem{proposition}{Proposition}

\begin{document}
\author{Gheorghe Minea}
\title{On the proper meaning of the curvature tensor and its general framework}
\maketitle
\begin{titlepage}
\begin{abstract}
We make evident a curvature tensor for every vector subbundle of an arbitrary manifold tangent bundle which reduces to the curvature 
tensor of an Ehresmann connection in the case of the horizontal subbundle of the tangent bundle to the total space of the nonlinear 
fiber bundle on which the connection is defined. Then the classical theorem of Frobenius would characterize the complete integrability 
of a vector subbundle of the tangent bundle by a zero curvature tensor in the sense of our definition here. A basic tool is a result about the
curvature tensor of the natural lift of the vector subbundle to a manifold of maps with values in the base of that subbundle. Another is a 
localization property for a Lie algebra of vector fields over this manifold of maps. These allow to prove an additive formula for the curvature 
tensors of two supplementary subbundles. The main result consists in identifying a natural linear parallel transport on a supplementary vector 
subbundle along any tangent path to the vector subbundle under study, which is the right generalization of a linear connection parallel transport
on a vector bundle along the projection in the base of that path. Then we derive the differential equation of the quotient of respective parallel 
transport operators induced by two different supplementary subbundles to the subbundle in question in terms of its curvature. Using this we 
obtain the equation of the infinitesimal variation of tangent paths to a vector subbundle, defined by its curvature, that appears as the root 
for the Jacobi equation of the infinitesimal variation of geodesics.
\end{abstract} 
\end{titlepage}

\textbf{Contents}\\

\S 1. Introduction. The curvature tensor of a vector subbundle of the tangent bundle.\\

\S 2. The definition of the curvature tensor and its expression in local coordinates. The two nondegenerate  cases and some integrability
questions.\\

\S 3. The lift of the curvature tensor to a manifold of maps with values in the base of the subbundle.\\

\S 4. On a localization property.\\

\S 5. An additive formula for the curvature tensors of two supplementary subbundles.\\

\S 6. The parallel transport in a supplementary vector subbundle along a tangent path to the vector subbundle under study.\\

\S 7. The equation of infinitesimal variation through tangent paths to the given vector subbundle.\\

\section{Introduction. The curvature tensor of a vector subbundle of the tangent bundle}

This paper is intended to fill a gap, due to a remarkable blind hens' phenomenon, regarding the notion of curvature. After the important work of 
Ehresmann, this missing step forward needed sixty four years! I became aware of the right definition of the curvature tensor some seven years 
ago, but I could not believe it is not known: in the paper [3] I took this definition for granted and I sent for it to the book [7] - where it 
does not appear! Only after looking through the beautiful book of Sternberg [6] I decided first to write a half page paper devoted only to this 
definition, in fact the content of Theorem 1 here. This approach was suggested by the double formulation (and proof) of the Frobenius theorem
from the book [4] of Narasimhan. Indeed, there is used the formula (\ref{18}) from here to pass from one formulation to the other. And precisely 
this formula allows to introduce our \textit{definition of the curvature tensor for every smooth vector subbundle of any manifold tangent bundle}.
In fact, for $H_{p}\subseteq T_{p}M$ smooth vector subbundle, the curvature tensor at $p\in M$ is a linear operator  
$C^{H}_{p}\in\mathrm{Hom}\,(H_{p}\wedge H_{p},T_{p}M/H_{p})$. Then the theorem of Frobenius asserts the equivalence between the complete 
integrability of such a vector subbundle and the equality to zero, at every point of the manifold, of its curvature tensor, in the sense of our 
definition here. The local expression of the curvature tensor shows that the Christoffel symbols appear even in this general case and that their 
meaning comes from the local representation of the respective vector subbundle of $TM$ using charts on $M$ and on a Grassmann manifold.\\
Next, in section \S 2, it is shown that for any $p_{0}\in M,\, H_{0}\subseteq T_{p_{0}}M$ vector subspace and 
$C_{0}\in\mathrm{Hom}\,(H_{0}\wedge H_{0},T_{p_{0}}M/H_{0})$ there is a locally defined, around $p_{0}$, vector subbundle $H$ of $TM$ such that
$H_{p_{0}}  = H_{0}$ and $C^{H}_{p_{0}} = C_{0}$. Taking into account the integrability property of the subbundle $H$ in the most degenerate case, 
when $C^{H}_{p} = 0$ for all $p\in M$, we point out in Proposition 4, dual properties of non-integrability of $H$ corresponding to the two
cases when the linear operator $C^{H}_{p}$ is non-degenerate, hence injective or surjective, respectively.\\
Having in view to study tangent paths to $H$ in $M$ and variations of them through tangent paths to $H$ alike, we consider in section \S 3 the
manifolds $\textsl{C}^{\infty}(D,M)$, for $D$ a compact domain, meant to be of dimension 1 or 2 respectively. If $H \longrightarrow M$ is the 
vector subbundle of $TM$, there is a natural lift of it to a vector bundle $\textsl{C}^{\infty}(D,H)\longrightarrow\textsl{C}^{\infty}(D,M)$,
which appears to be a vector subbundle of $T\textsl{C}^{\infty}(D,M)$, if we consider the (formal) isomorphism
$T\textsl{C}^{\infty}(D,M)\widetilde{\longrightarrow}\textsl{C}^{\infty}(D,TM)$. The result of \S 3, Theorem 2, essential in all the paper, is a 
simple and natural formula for the curvature of $\textsl{C}^{\infty}(D,H)$, as a vector subbundle of $T\textsl{C}^{\infty}(D,M)$, in terms of the 
curvature of $H$.\\
The section \S 4 is devoted to a Lie subalgebra of vector fields over $\textsl{C}^{\infty}(D,M)$ which is isomorphic to the set of global sections 
of a sheaf of Lie algebras over $D\times M$. The localization property of this Lie algebra structure is used in \S 5 to prove a certain identity 
(Theorem 4) that entails the formula expressing the sum of curvatures of two supplementary subbundles by such a Lie bracket (Theorem 5). But the
formula for the Lie bracket of this sheaf sections (Theorem 3) will play also a role in section \S 6.\\
This section \S 6 is the core of the whole paper. We show that for every tangent to the subbundle $H$ path $\gamma$ and for every supplementary 
to $H$ vector subbundle $K$ there exists a natural linear parallel transport of the fibers of $K$ along the path $\gamma$. It is the right 
generalization of the case when $M$ is the total space of a vector bundle, over $B$ say, and $\gamma$ is the zero lift of an arbitrary path in 
$B$. Then $\gamma$ is tangent to the horizontal subbundle $H$ defining a linear connection on $M$. Next, the linear fibers of $M$, over $B$, can
be identified with the tangent to them in zero, giving thus the supplementary subbundle $K$ along $\gamma$. And, in this case, the parallel 
transport of the fibers of $M$ along the original path in the base $B$, given by the linear connection, coincides with the parallel transport of
the fibers of $K$ along $\gamma$. If $M$ is the total space of a nonlinear fiber bundle on which is defined an Ehresmann connection by $H$ and
$K$ is the tangent to the nonlinear fibers, then the linear parallel transport operators of $K$ along $\gamma$, found by us, appear to be the 
tangent operators, in the points of the tangent to $H$ path $\gamma$, to the diffeomorphisms of nonlinear fibers given by the Ehresmann 
connection parallel transport along the projection of $\gamma$ in the base (Theorem 8).\\
Coming back to our general framework, if we consider the linear parallel transport operators of the fibers of $TM/H$ defined by two supplementary
subbundles $K^{1}$ and $K^{2}$ to $H$, then there is an interesting equation of evolution of their quotient, along $\gamma$, in terms of the 
curvature of $H$ (Theorem 10). This equation is the consequence of Theorem 9, where we find the equation of the infinitesimal variation of 
tangent paths to $H$ using a supplementary subbundle $K$ to $H$. Then in section \S 7 we prove the reciprocal of this statement and show that the 
equation, defined by the curvature of $H$, is independent of $K$ (Theorems 11 and 12). Finally, we show in Theorem 13 that this equation is the 
root for the equation of Jacobi of infinitesimal variation of geodesics.

\section{The definition of the curvature tensor\\ and its expression in local coordinates.\\
The two nondegenerate cases and some\\ integrability questions}

\noindent Let $M$ be a $\textsl{C}^{\infty}$ (finite dimensional) manifold and $H_{p}\subseteq T_{p} M, p\in M,$ be a $\textsl{C}^{\infty}$ 
vector subbundle of its tangent bundle. For a diffeomorphism $\varphi :M\longrightarrow N$ we denote $\varphi_{\ast} H$ the vector subbundle 
of $TN$:
\begin{equation}\label{1}
 (\varphi_{\ast} H)_{\varphi (p)} = T_{p}\varphi\cdot H_{p},\;\; p\in M.
\end{equation}
For a Banach vector space $E$ we consider
\begin{equation}\label{2}
  \varPsi^{E} : E\times E\longrightarrow T E,\;\varPsi^{E}(v,w) = \dfrac{\mathrm{d}}{\mathrm{d}t} (v + t w)\rvert_{t = 0} \in T_{v} E,
\end{equation}
which is a bijection and $ \varPsi^{E} (x,\cdot) : E\longrightarrow T_{x} E$ is an isomorphism $\forall x\in E$. Then, if $M$ is a Banach 
manifold, $E$ a Banach vector space and $f : M\longrightarrow E$ is a $\textsl{C}^{1}$ - function, we denote
\begin{equation}\label{3}
 \mathrm{d}_{p} f : T_{p} M\longrightarrow E,\;\;\mathrm{d}_{p} f = \varPsi^{E} (f(p),\cdot )^{-1}\cdot T_{p} f
\end{equation}
its differential at $p\in M$. And if $\chi : U\longrightarrow E,\;\; U = \mathring{U} \subseteq M$, is a local chart on $M$, we consider the
subspaces  of $E$
\begin{equation}\label{4}
  \widetilde{(\chi_{\ast} H)}_{\chi (p)} =  \mathrm{d}_{p}\chi\cdot H_{p},\;\;p\in U.
\end{equation}
Using a chart in the neighbourhood of $p_{0}\in M$ we may construct a local trivialization of $H$ around $p_{0}$ in the following way: let
$V := \widetilde{(\chi_{\ast} H)}_{\chi (p_{0})}$ and $W$ such that $V\dot{+} W = E$, where $\dot{+}$ stands for the interior direct sum. 
Then a natural isomorphism allows to consider a modified chart
\begin{equation}\label{5}
 \chi : U\longrightarrow V\times W
\end{equation}
such that $ \mathrm{d}_{p_{0}}\chi\cdot H_{p_{0}} = V\times \{0_{W}\}$. The smoothness of $H$ is equivalent to the smoothness of the map 
$p\mapsto \mathrm{d}_{p}\chi\cdot H_{p}$ taking values in the Grassmann manifold of subspaces of $V\times W$ isomorphic with $V$. In this way we 
find $\varGamma (x,y)\in\mathrm{Hom} (V,W)$ depending smoothly on $(x,y)\in V\times W$ in a neighbourhood of $(x_{0},y_{0}) = \chi (p_{0})$
such that $\varGamma (x_{0},y_{0}) = 0_{\mathrm{Hom} (V,W)}$ and
\begin{equation}\label{6}
  \widetilde{(\chi_{\ast} H)}_{(x,y)} = \mathrm{graph}\,\varGamma (x,y) = \{(v,\varGamma (x,y) v)\;\rvert\;v\in V\}.
\end{equation}
In that neighbourhood we consider the isomorphism
\begin{equation}\label{7}
 \varLambda^{V}_{(x,y)} : V\widetilde{\longrightarrow}\widetilde{(\chi_{\ast} H)}_{(x,y)},\;\;\varLambda^{V}_{(x,y)} v = (v,\varGamma (x,y) v),
\end{equation}
and finally the trivialization of $H$
\begin{equation}\label{8}
 \vartheta_{p} : H_{p}\widetilde{\longrightarrow} V,\;\; \vartheta_{p} = (\varLambda^{V}_{\chi (p)})^{-1}\cdot \mathrm{d}_{p}\chi\rvert_{H_{p}}
 ,\;\; p\in U.
\end{equation}
We will call a chart (\ref{5}) where $H$ is represented by (\ref{6}) an \textit{adapted chart} for the vector subbundle $H$. If 
$H_{p}^{\bot}\subseteq T^{\ast}_{p} M$ denotes the orthogonal with respect to the duality $\{T_{p} M,T_{p}^{\ast} M\}$ then $H^{\bot}$ becomes a 
$\textsl{C}^{\infty}$ vector subbundle of $T^{\ast} M$. In the case of a $H$-adapted chart we denote
\begin{equation}\label{9}
 \varLambda^{W}_{(x,y)} : W\widetilde{\longrightarrow}E/\mathrm{graph}\,\varGamma(x,y),\;\;\varLambda^{W}_{(x,y)} w = (0_{V},w) + 
 \mathrm{graph}\,\varGamma(x,y),
\end{equation}
where $E = V\times W$, and then its transposed establishes an isomorphism
\begin{equation}\label{10}
 (\varLambda^{W}_{(x,y)})^{\ast} : (\mathrm{graph}\,\varGamma(x,y))^{\bot}\longrightarrow W^{\ast}.
\end{equation}
On the other hand $((\mathrm{d}_{p}\chi)^{-1})^{\ast} : T_{p}^{\ast} M\longrightarrow E^{\ast}$ is an isomorphism between $H_{p}^{\bot}$ and\\
$(\mathrm{d}_{p}\chi\cdot H_{p})^{\bot}$, so that
\begin{equation}\label{11}
 \vartheta_{p}^{\bot} : H_{p}^{\bot}\widetilde{\longrightarrow} W^{\ast},\;\; \vartheta_{p}^{\bot} = 
 (\varLambda^{W}_{\chi (p)})^{\ast}\cdot ((\mathrm{d}_{p}\chi)^{-1})^{\ast}\rvert_{ H_{p}^{\bot}}
\end{equation}
gives a trivialization of $H^{\bot}$ in the same neighbourhood. There, a section $X\in\textsl{C}^{\infty}\Gamma (H)$ is represented by 
$f = \vartheta\circ X\circ\chi^{-1}\in\textsl{C}^{\infty} (U,V)$ such that
\begin{equation}\label{12}
 \mathrm{d}_{\chi^{-1} (x,y)}\chi\cdot X_{\chi^{-1} (x,y)} = (f(x,y),\varGamma (x,y)\;f(x,y))
\end{equation}
and a section $\alpha\in\textsl{C}^{\infty}\Gamma (H^{\bot})$ by $\varphi = 
\vartheta_{p}^{\bot}\circ\alpha\circ\chi^{-1}\in\textsl{C}^{\infty} (U,W^{\ast})$ for which
\begin{equation}\label{13}
 <((\mathrm{d}_{\chi^{-1}(x,y)}\chi)^{\ast})^{-1}\alpha_{\chi^{-1}(x,y)}, (v,w)> = <\varphi (x,y),w - \varGamma (x,y) v>.
\end{equation}

Our starting point is the following elementary
\begin{theorem}. For $X, Y\in\textsl{C}^{\infty}\Gamma (H)$ and $\alpha\in\textsl{C}^{\infty}\Gamma (H^{\bot})$ the following equality holds\\
$\forall p\in M$:
\begin{equation}\label{14}
  \mathrm{d}_{p}\alpha(X_{p},Y_{p}) = - <\alpha_{p},\big[X,Y\big]_{p}>. 
\end{equation}
It follows that for $p_{0}\in M$ fixed it is well defined a trilinear map
\begin{equation}\label{15}
 \tau_{p_{0}} : H_{p_{0}}\times H_{p_{0}}\times H_{p_{0}}^{\bot}\longrightarrow \textbf{R}
\end{equation}
by
\begin{equation}\label{16}
 \tau_{p_{0}}(u,v,\phi) = <\phi,\big[X,Y\big]_{p_{0}}> = -\mathrm{d}_{p_{0}}\alpha (u,v),
\end{equation}
where $X, Y\in\textsl{C}^{\infty}\Gamma (H)$, $\alpha\in\textsl{C}^{\infty}\Gamma (H^{\bot})$ are arbitrary with
\begin{equation}\label{17}
 X_{p_{0}} = u, Y_{p_{0}} = v, \alpha_{p_{0}} = \phi.
\end{equation}
\end{theorem}
\textbf{Proof}. The equality (\ref{14}) follows from the general identity
\begin{equation}\label{18}
 \mathrm{d}_{p}\alpha(X_{p},Y_{p}) = <\mathrm{d}_{p}<\alpha,Y>,X_{p}> - <\mathrm{d}_{p}<\alpha,X>,Y_{p}> - <\alpha_{p},\big[X,Y\big]_{p}>,
\end{equation}
for $X, Y\in\textsl{C}^{\infty}\Gamma (TM), \alpha\in\textsl{C}^{\infty}\Gamma (T^{\ast}M)$. Then in (\ref{16}) the second term shows
that $\tau_{p_{0}}$ depends only on $\phi$ while the third term shows that $\tau_{p_{0}}$ depends only on $u$ and $v$. On the other hand, the 
representations (\ref{12}) and (\ref{13}) allow to construct sections for $H$ and $H^{\bot}$ with prescribed value in a given point. Using this 
procedure we can also verify that $\tau_{p_{0}}$ is trilinear $\blacksquare$\\
We consider then the following\\
\textbf{Definition}. If $H$ is a vector subbundle of $TM$ we define its \textit{curvature tensor} in $p_{0}\in M$ as the operator
\begin{equation}\label{19}
 C_{p_{0}}^{H}\in\mathrm{Hom}(H_{p_{0}}\wedge H_{p_{0}}, T_{p_{0}}M /H_{p_{0}}),
\end{equation} 
\begin{equation}\label{20}
  C_{p_{0}}^{H}(u\wedge v) := P_{p_{0}}^{H}\big[X,Y\big]_{p_{0}},
\end{equation}
if 
\begin{equation}\label{21}
 P_{p_{0}}^{H} :  T_{p_{0}}M\longrightarrow T_{p_{0}}M /H_{p_{0}}
\end{equation}
denotes the canonical projection and $X, Y\in\textsl{C}^{\infty}\Gamma (H)$ are chosen so that 
$X_{p_{0}} = u,\\
Y_{p_{0}} = v\;\; \blacksquare$

Then the classical \textit{Frobenius theorem} may be phrased as stating that $H$ \textit{is completely integrable if and only if} 
$C_{p}^{H} = 0, \forall p\in M$ (see Narasimhan [4]). \\
We look now for a local representation of the curvature tensor corresponding to the local form (\ref{6}) of the vector subbundle $H$. For any
diffeomorphism $\varphi : M\longrightarrow N$ we have also the canonical isomorphism 
\begin{equation}\label{22}
 T_{p}\varphi/H_{p} : T_{p} M/H_{p}\widetilde{\longrightarrow} T_{\varphi(p)} N/(\varphi_{\ast}H)_{\varphi (p)}
\end{equation}
(see (\ref{1})). If $C^{\varphi_{\ast} H}$ is the curvature tensor of the subbundle $\varphi_{\ast} H$ of $TN$, we have
\begin{equation}\label{23}
 T_{p}\varphi/H_{p}\cdot C^{H}_{p} (X_{p}\wedge Y_{p}) = C^{\varphi_{\ast} H}_{\varphi (p)} (T_{p}\varphi\cdot X_{p}\wedge T_{p}\varphi\cdot Y_{p}),
\end{equation}
for $X_{p},\;Y_{p}\in H_{p},\;p\in M$. If $\chi : U\longrightarrow E,\;\; U = \mathring{U} \subseteq M$, is a local chart on $M$, we denote
(see (\ref{4}))
\begin{equation}\label{24}
 \mathrm{d}_{p}\chi/H_{p} : T_{p} M/H_{p}\longrightarrow E/\widetilde{(\chi_{\ast} H)}_{\chi (p)}
\end{equation}
and define $\widetilde{C^{\chi_{\ast} H}_{e}}\in\mathrm{Hom} (\widetilde{(\chi_{\ast} H)}_{e}\wedge\widetilde{(\chi_{\ast} H)}_{e},
E/\widetilde{(\chi_{\ast} H)}_{e})$, for $e\in E$, by
\begin{equation}\label{25}
\mathrm{d}_{p}\chi / H_{p}\cdot C^{H}_{p} (X_{p}\wedge Y_{p}) = 
\widetilde{C^{\chi_{\ast} H}_{\chi (p)}} (\mathrm{d}_{p}\chi\cdot X_{p}\wedge\mathrm{d}_{p}\chi\cdot Y_{p}).
\end{equation}
For $X\in\textsl{C}^{\infty}\Gamma (TM)$ we have $\chi_{\ast} X = (\chi^{-1})^{\ast} X\in\textsl{C}^{\infty}\Gamma (TE)$ defined as usual
\begin{equation}\label{26}
 (\chi_{\ast} X)_{e} = T_{\chi^{-1} (e)}\chi\cdot X_{\chi^{-1} (e)}
\end{equation}
and for $X\in\textsl{C}^{\infty}\Gamma (TE)$ we denote
\begin{equation}\label{27}
 \tilde{X} = p_{2}\circ (\varPsi^{E})^{-1}\circ X,\;\;\tilde{X} : E\longrightarrow E,
\end{equation}
(where $p_{2} : E\times E\longrightarrow E$ is canonical), such that
\begin{equation}\label{28}
 X_{e} = \varPsi^{E}(e,\tilde{X}(e)),\;\;e\in E.
\end{equation}
From (\ref{26}) and (\ref{27}) we get
\begin{equation}\label{29}
 \widetilde{\chi_{\ast} X}(e) = \mathrm{d}_{\chi^{-1}(e)}\chi\cdot X_{\chi^{-1}(e)}.
\end{equation}
And if we define for $F,\;G\in\textsl{C}^{\infty} (E,E)$
\begin{equation}\label{30}
 \big[F,G\big](e) = G^{\prime}(e)\cdot F(e) - F^{\prime}(e)\cdot G(e),
\end{equation}
for $X,\;Y\in\textsl{C}^{\infty}\Gamma (TE)$ we have
\begin{equation}\label{31}
 \widetilde{\big[X,Y\big]} = \big[\tilde{X},\tilde{Y}\big].
\end{equation}
Since $(\chi^{-1})^{\ast}\big[X,Y\big] = \big[(\chi^{-1})^{\ast} X,(\chi^{-1})^{\ast} Y\big]$, we have
\begin{equation}\label{32}
 \widetilde{\chi_{\ast}\big[X,Y\big]} = \big[\widetilde{\chi_{\ast} X}, \widetilde{\chi_{\ast} Y}\big]
\end{equation}
for $X,\;Y\in\textsl{C}^{\infty}\Gamma (TM)$.\\
In order to distinguish the type of application to different vectors we use here the notations
\begin{equation}\label{33}
 <\dfrac{\partial F}{\partial x}(x,y); f> = \dfrac{\mathrm{d}}{\mathrm{d} t}F(x + tf, y)\rvert_{t=0},\;\;\;
<\dfrac{\partial F}{\partial y}(x,y); h> = \dfrac{\mathrm{d}}{\mathrm{d} t}F(x, y + th))\rvert_{t=0},
\end{equation}
for $F : O\longrightarrow Z,\;\;Z$ vector space and $O = \mathring{O}\subseteq V\times W$. In the case  $Z = \mathrm{Hom}(V,W)$  
$<\dfrac{\partial F}{\partial x}(x,y); f>g\in W$ has a clear meaning for
$f, g\in V$; also $<\dfrac{\partial F}{\partial y}(x,y); h>g\in W$, for $h\in W, g\in V$.\\
Coming back to (\ref{32}) we get for
\begin{equation}\label{34}
 \widetilde{\chi_{\ast} X}(x,y) = (f(x,y),\varGamma (x,y)\;f(x,y)),\;\; \widetilde{\chi_{\ast} Y}(x,y) = (g(x,y),\varGamma (x,y)\;g(x,y))
\end{equation}
$\big[\widetilde{\chi_{\ast} X},\widetilde{\chi_{\ast} Y}\big](x,y) = \big ( <\dfrac{\partial g}{\partial x};f> + 
<\dfrac{\partial g}{\partial y};\varGamma\;f> - <\dfrac{\partial f}{\partial x};g> - <\dfrac{\partial f}{\partial y};\varGamma\;g> ,\\
<\dfrac{\partial \varGamma}{\partial x};f>g + \varGamma <\dfrac{\partial g}{\partial x};f> + 
<\dfrac{\partial \varGamma}{\partial y};\varGamma f>g + \varGamma <\dfrac{\partial g}{\partial y};\varGamma f> - 
<\dfrac{\partial \varGamma}{\partial x};g>f -\\
- \varGamma <\dfrac{\partial f}{\partial x};g> - <\dfrac{\partial \varGamma}{\partial y};\varGamma g>f - 
\varGamma <\dfrac{\partial f}{\partial y};\varGamma g>\big )$.\\
Let us consider
\begin{equation}\label{35}
 P^{\widetilde{\chi_{\ast} H}}_{e} : E\longrightarrow E/\widetilde{(\chi_{\ast} H)}_{e},\;\;e\in E,
\end{equation}
canonical and remark that (see (\ref{24}))
\begin{equation}\label{36}
 \mathrm{d}_{p}\chi/H_{p}\cdot P^{H}_{p} =  P^{\widetilde{\chi_{\ast} H}}_{\chi (p)} \cdot \mathrm{d}_{p}\chi,
\end{equation}
$p\in U$. Then\\

$P^{\widetilde{\chi_{\ast} H}}_{(x,y)}\big[\widetilde{\chi_{\ast} X},\widetilde{\chi_{\ast} Y}\big](x,y) =$\\

$= P^{\widetilde{\chi_{\ast} H}}_{(x,y)}\big (0,
<\dfrac{\partial \varGamma}{\partial x};f>g\; - <\dfrac{\partial \varGamma}{\partial x};g>f + 
<\dfrac{\partial \varGamma}{\partial y};\varGamma f>g\; - <\dfrac{\partial \varGamma}{\partial y};\varGamma g>f\big )$,\\

or (see(\ref{9}))
\begin{multline}\label{37}
 P^{\widetilde{\chi_{\ast} H}}_{(x,y)}\big[\widetilde{\chi_{\ast} X},\widetilde{\chi_{\ast} Y}\big](x,y) = \\
 = \varLambda^{W}_{(x,y)}\big (<\dfrac{\partial \varGamma}{\partial x};f>g\; - <\dfrac{\partial \varGamma}{\partial x};g>f + 
<\dfrac{\partial \varGamma}{\partial y};\varGamma f>g\; - <\dfrac{\partial \varGamma}{\partial y};\varGamma g>f\big )
\end{multline}
Let us consider
\begin{equation}\label{38}
 \widehat{C^{\chi_{\ast}H}_{(x,y)}}\in\mathrm{Hom} (V\wedge V,W)
\end{equation}
defined by
\begin{multline}\label{39}
 \widehat{C^{\chi_{\ast}H}_{(x,y)}}(f\wedge g) = 
<\dfrac{\partial \varGamma}{\partial x}(x,y);f>g\; - <\dfrac{\partial \varGamma}{\partial x}(x,y);g>f +\\
+ <\dfrac{\partial \varGamma}{\partial y}(x,y);\varGamma(x,y) f>g\; - <\dfrac{\partial \varGamma}{\partial y}(x,y);\varGamma(x,y) g>f.
\end{multline}
As $\widetilde{(\chi_{\ast}[X,Y])}_{\chi (p)} = \mathrm{d}_{p}\chi\cdot \big[X,Y\big]_{p}$, from (\ref{32}), (\ref{37}) and (\ref{39}) we get
\begin{equation}\label{40}
\mathrm{d}_{p}\chi / H_{p}\cdot C^{H}_{p} ((\mathrm{d}_{p}\chi)^{-1}\varLambda^{V}_{\chi (p)} f\wedge
 (\mathrm{d}_{p}\chi)^{-1}\varLambda^{V}_{\chi (p)} g) = \varLambda^{W}_{\chi (p)}\cdot\widehat{C^{\chi_{\ast}H}_{\chi (p)}}(f\wedge g),
\end{equation}
for $f,\;g\in V$. We say that (\ref{39}) is the \textit{local representation of the curvature tensor} in an adapted chart.
We recall that here $H$ is defined by $\varGamma$ according to (\ref{6}).\\
The formula (\ref{39}) shows the proper and general meaning of the \textit{Christoffel symbols}, as defined by $\varGamma$ from (\ref{6}),
in virtue of the role of $\varGamma$ in the expression of the curvature tensor (see, for instance, Kobayashi+Nomizu, vol I [2], 
Sternberg [6]): they simply give locally the subbundle whose curvature is computed.\\
We examine now the structure of the curvature tensor at a fixed point.
\begin{proposition}. For $p_{0}\in M$ fixed, the operator $C_{p_{0}}^{H}$ may take any value from \\
 $\mathrm{Hom}(H_{p_{0}}\wedge H_{p_{0}}, T_{p_{0}}M /H_{p_{0}})$,
 when $H_{p}$ varies with $p$ around $p_{0}$ and $H_{p_{0}}$ is prescribed. 
\end{proposition}
\textbf{Proof}. We consider
\begin{equation}\label{41}
 M=V\times W
\end{equation}
and $H_{(x,y)}$ given by (\ref{6}) and $C_{(x,y)}^{H}$ by (\ref{39}). We take next $C\in\mathrm{Hom}(V\wedge V, W)$ arbitrary and define
$H$ by
\begin{equation}\label{42}
 \varGamma (x,y)v = \dfrac{1}{2} C(x\wedge v),\;\;x, v\in V, \;y\in W.
\end{equation}
Then $<\dfrac{\partial\varGamma}{\partial x}(x,y); u>v = \dfrac{1}{2} C(u\wedge v),\;\dfrac{\partial\varGamma}{\partial y}(x,y) = 0$,
wherefrom, using also (\ref{39}), for this $H$:
\begin{equation}\label{43}
  C_{(x,y)}^{H} = C,\;\;\forall (x,y)\in V\times W
\end{equation}
$\blacksquare$\\
The following fact comes straight from the definition of the curvature tensor.
\begin{proposition}. Let $H_{p} \subseteq K_{p}\subseteq T_{p}M$ be two smooth vector subbundles and let
\begin{equation}\label{44}
 j_{p} : H_{p}\longrightarrow K_{p},\;\;q_{p} :T_{p}M/H_{p}\longrightarrow T_{p}M/K_{p}
\end{equation}
be the canonical maps. Then $\forall p\in M$
\begin{equation}\label{45}
 C_{p}^{K} (j_{p}u\wedge j_{p} v) = q_{p} C_{p}^{H} (u\wedge v),\;\;\forall u, v\in H_{p}.
\end{equation}
\end{proposition}
$\blacksquare$\\
Next we try to identify in suitable properties of the curvature tensor the existence of any of two types of integral submanifolds.
\begin{proposition}.i) If $N\subseteq M$ is a submanifold with
\begin{equation}\label{46}
 T_{p}N\subseteq H_{p},\;\;\forall p\in N,
\end{equation}
then
\begin{equation}\label{47}
  C_{p}^{H}(u\wedge v) = 0,\;\;\forall p\in N,\;\;\forall u, v\in T_{p}N.
\end{equation}
 ii) If $N\subseteq M$ is a submanifold with
\begin{equation}\label{48}
 T_{p}N\supseteq H_{p},\;\;\forall p\in N,
\end{equation}
then
\begin{equation}\label{49}
 C_{p}^{H}(u\wedge v)\in T_{p} N/H_{p},\;\;\forall p\in N,\;\;\forall u,v\in H_{p}  
\end{equation}
 (as $T_{p} N/H_{p}\subseteq T_{p} M/H_{p})$).
\end{proposition}
\textbf{Proof}. It is easy to see that for each submanifold $N$ and each $p\in N$ there exists a neighbourhood around $p$ in $M$ and a 
completely integrable subbundle $K\subseteq TM$ defined on it such that
\begin{equation}\label{50}
 K_{p} = T_{p}N,\;\;\forall p\in N.
\end{equation}
i) In this case $K_{p}\subseteq H_{p}$ and $C_{p}^{K} = 0$; then from (\ref{45}) we infer (\ref{47}).\\
ii) In this case $K_{p}\supseteq H_{p}$ and $C_{p}^{K} = 0$; then (\ref{45}) gives $q_{p} C_{p}^{H} (u\wedge v) = 0,\;\;\forall p\in N,
\;\;u, v\in H_{p}$. But $\mathrm{ker}\, q_{p} = K_{p}/H_{p} = T_{p}N/H_{p}$ (see (\ref{44}) and (\ref{50}))$\blacksquare$\\

We will say that the curvature tensor is \textit{nondegenerate} at $p_{0}$ if the operator\\ 
$C_{p_{0}}^{H}\in\mathrm{Hom}\,(H_{p_{0}}\wedge H_{p_{0}}, T_{p_{0}}M /H_{p_{0}})$ has maximal rank. Hence there are two cases, depending on
$m = \mathrm{dim}\, H_{p_{0}}$ and $n = \mathrm{dim}\, T_{p_{0}} M$ : in the case that $\binom{m}{2} \leqslant n - m, \;\;C_{p_{0}}^{H}$ should
be injective and in the case that $\binom{m}{2} \geqslant n - m, \;\;C_{p_{0}}^{H}$ should be surjective. Of course, if the curvature tensor 
is nondegenerate at a certain point, then it is nondegenerate at each point from a neighbourhood of it.

\begin{proposition}.i) If $C_{p_{0}}^{H}$ is injective then there does not exist a submanifold $N\subseteq M,\\ N\ni p_{0}$, such that
\begin{equation}\label{51}
 T_{p} N\subseteq H_{p},\;\;\forall p\in N,\;\;\mathrm{dim}\, N\geqslant 2.
\end{equation}
ii) If  $C_{p_{0}}^{H}$ is surjective then there does not exist a submanifold $N\subseteq M, N\ni p_{0}$, such that
\begin{equation}\label{52}
 T_{p} N\supseteq H_{p},\;\;\forall p\in N,\;\;\mathrm{codim}\, N\geqslant1.
\end{equation}
\end{proposition}
\textbf{Proof}.i)  From (\ref{47}), Proposition 3 i) and the injectivity it follows that\\
$T_{p} N\wedge T_{p} N\subseteq \mathrm{ker}\, C_{p}^{H} =\{0\},\;\;\forall p\in N$ and in a neighbourhood of $p_{0}$, and then 
$\mathrm{dim}\, N\leqslant 1$.\\
ii) Using Proposition 3 ii) and the surjectivity we get $T_{p} N/H_{p} = T_{p} M/H_{p}$, hence \\
$T_{p} N = T_{p} M$ and $\mathrm{codim}\, N = 0\;\;\blacksquare$ 

But the previous Proposition can be supplemented with the counterexamples from the next
\begin{proposition}.i) In the case that $C_{p}^{H}$ is injective $\forall p\in M$ it may exist a submanifold $N\subseteq M$ with
\begin{equation}\label{53}
 T_{p} N\supseteq H_{p},\;\;\forall p\in N,\;\;\mathrm{codim}\, N\geqslant1. 
\end{equation}
More precisely, this is possible in the case of (\ref{41}), (\ref{42}) with $C$ injective and
\begin{equation}\label{54}
 \mathrm{dim}\, (V\times W) = 4,\;\; \mathrm{dim}\, V = \mathrm{dim}\, H = 2,
\end{equation}
when we can find $N$ of
\begin{equation}\label{55}
 \mathrm{dim}\, N = 3. 
\end{equation}
ii) In the case that $C_{p}^{H}$ is surjective $\forall p\in M$ it may exist a submanifold $N\subseteq M$ with 
 \begin{equation}\label{56}
 T_{p} N\subseteq H_{p},\;\;\forall p\in N,\;\;\mathrm{dim}\, N\geqslant 2.
 \end{equation}
 Again, this can be realized for (\ref{41}), (\ref{42}) with $C$ surjective,
 \begin{equation}\label{57}
  \mathrm{dim}\,(V\times W) = 4,\;\; \mathrm{dim}\,V = \mathrm{dim}\,H =3, 
 \end{equation}
with a certain $N$ of
\begin{equation}\label{58}
 \mathrm{dim}\,N = 2. 
\end{equation}
\end{proposition}
\textbf{Proof}.i) We take bases such that $V = \textbf{R} e_{1}\dot{+} \textbf{R} e_{2},\;\; W = \textbf{R} e_{3}\dot{+} \textbf{R} e_{4}$ and
define $C$ by
\begin{equation}\label{59}
 C(e_{1}\wedge e_{2}) = e_{3}.
\end{equation}
Then from (\ref{42}) we have $\mathrm{im}\,\varGamma (x^{\prime},x^{\prime\prime}) = \textbf{R} e_{3},\;\forall x^{\prime}\in V,\;
x^{\prime} \neq 0,\;\forall x^{\prime\prime}\in W$ and \\$H_{(x^{\prime},x^{\prime\prime})} = 
\{v\oplus\varGamma(x^{\prime},x^{\prime\prime})v \;\arrowvert\ v\in V\}\subset \textbf{R} e_{1}\dot{+} \textbf{R} e_{2}\dot{+}\textbf{R} e_{3}$.
Therefore the hypersurfaces $N_{\lambda}\subset V\times W$,
\begin{equation}\label{60}
 N_{\lambda} = \{x_{4} = \lambda\},
\end{equation}
$\lambda\in \textbf{R}$, with 
$T_{(x^{\prime},x^{\prime\prime})} N_{\lambda}= \textbf{R} e_{1}\dot{+} \textbf{R} e_{2}\dot{+}\textbf{R} e_{3}$, have the desired properties 
(\ref{53}), (\ref{55}).\\
ii) We take $V =  \textbf{R} e_{1}\dot{+} \textbf{R} e_{2}\dot{+}\textbf{R} e_{3},\; W = \textbf{R} e_{4}, \;C : V\wedge V\longrightarrow W$ by 
\begin{equation}\label{61}
 C(e_{1}\wedge e_{2}) =  C(e_{2}\wedge e_{3}) = 0,\; C(e_{1}\wedge e_{3}) = e_{4}.
\end{equation}
Let us consider the surfaces
\begin{equation}\label{62}
N_{\lambda} = \{x_{3} = 0,\;x_{4} = \lambda\}, 
\end{equation}
$\lambda\in \textbf{R}$. For $x=(x^{\prime},x^{\prime\prime})\in N_{\lambda},\;x^{\prime}\in V,\;x^{\prime\prime}\in W$, we have $x^{\prime}\in
 \textbf{R} e_{1}\dot{+} \textbf{R} e_{2}$. On the other hand $T_{(x^{\prime},x^{\prime\prime})} N_{\lambda} = 
 \textbf{R} e_{1}\dot{+} \textbf{R} e_{2}$ and if $z\in T_{(x^{\prime},x^{\prime\prime})} N_{\lambda}$ then $C(x^{\prime}\wedge z) = 0$ wherefrom
 $z= z\dot{+} \dfrac{1}{2} C(x^{\prime}\wedge z) = z\dot{+}\varGamma(x^{\prime},x^{\prime\prime}) z\in H_{(x^{\prime},x^{\prime\prime})}$; thus
(\ref{56}) holds$\;\blacksquare$

\section{The lift of the curvature tensor to a manifold\\
of maps with values in the base of the subbundle}

The formula that is proved here involves the Lie bracket of vector fields on the infinite dimensional manifold $\textsl{C}^{\infty}(D,M)$, where 
$D$ is an arbitrary \textit{compact} manifold (possibly with boundary) and $M$ is the manifold where the vector subbundles of $TM$ are 
considered. The special form of the vector fields that appear in it allows to work with a certain Lie algebra that does not use  
a differentiable structure on $\textsl{C}^{\infty}(D,M)$ , being enough to compute
on the manifolds $\textsl{C}^{k}(D,M)$ modelled on Banach spaces. We refer to R. Palais [5] for the facts concerning these Banach manifolds that
we use in the sequel. The suitable framework is that of $\textsl{C}^{\infty}$ nonlinear fiber bundles $\pi : \eta\longrightarrow D$ and of the 
respective Banach manifolds $\textsl{C}^{k}\Gamma (\eta)$ of $\textsl{C}^{k}$ sections.\\
The main tool in order to introduce a structure of
$\textsl{C}^{\infty}$ manifold on $\textsl{C}^{k}\Gamma (\eta)$ is the existence, for each $\sigma_{0}\in\textsl{C}^{0}\Gamma (\eta)$, of a 
\textit{vector bundle neighbourhood} of it, that is an open subset $E\subseteq \eta$ with $E\supseteq\sigma_{0}(D)$ such that the fiber bundle
$\pi\rvert_{E} : E\longrightarrow D$ be endowed with a compatible structure of $\textsl{C}^{\infty}$ vector bundle. The point is that the open 
subsets $\textsl{C}^{k}\Gamma (E)$ of $\textsl{C}^{k}\Gamma (\eta)$, with their own structure of Banach space, give, through their embedding,
the inverse of charts for that structure.\\
Let 
\begin{equation}\label{63}
 \tau_{\eta}\rvert_{TF\eta} :TF\eta\longrightarrow \eta,\;\;TF_{e}\; \eta = T_{e}\; \eta_{\pi (e)},\;e\in \eta,
\end{equation}
be the vector bundle of tangents to the fibers and let, for $\sigma\in\textsl{C}^{k}\Gamma (\eta)$,
\begin{equation}\label{64}
 \sigma^{\ast} (TF\eta)\longrightarrow D,\;\;[\sigma^{\ast} (TF\eta)]_{\zeta} = T_{\sigma_{\zeta}} \;\eta_{\zeta},
\end{equation}
be the vector bundle pulled back through $\sigma$. As a consequence of the specific differential structure of $\textsl{C}^{k}\Gamma (\eta)$
\begin{equation}\label{65}
\varPhi^{k}_{\sigma} : T_{\sigma}\textsl{C}^{k}\Gamma (\eta)\widetilde{\longrightarrow}\textsl{C}^{k}\Gamma (\sigma^{\ast}(TF\eta)),\;\;
[\varPhi^{k}_{\sigma (0)}(\dfrac{\mathrm{d}\sigma(t)}{\mathrm{d}t}\big |_{t=0})]_{\zeta} = 
\dfrac{\mathrm{d}[\sigma(t)_{\zeta}]}{\mathrm{d}t}\big |_{t=0},\;\zeta\in D,
\end{equation}
is an isomorphism. On the other hand
\begin{equation}\label{66}
 \pi\circ\tau_{\eta}\rvert_{TF\eta} : TF\eta\longrightarrow D,\;\;(TF\eta)_{\zeta} = T(\eta_{\zeta}),\;\zeta\in D,
\end{equation}
is a nonlinear fiber bundle over $D$ for which the manifold of sections  $\textsl{C}^{k}\Gamma (TF\eta)$ becomes a vector bundle
over $\textsl{C}^{k}\Gamma (\eta)$ through the projection
\begin{equation}\label{67}
 l(\tau_{\eta}\rvert_{TF\eta}) :  \textsl{C}^{k}\Gamma (TF\eta) \longrightarrow\textsl{C}^{k}\Gamma (\eta),\;\;
 l(\tau_{\eta}\rvert_{TF\eta})(\varphi) = \tau_{\eta}\circ\varphi,
\end{equation}
\begin{equation}\label{68}
 \textsl{C}^{k}\Gamma (TF\eta)_{\sigma} = [l(\tau_{\eta}\rvert_{TF\eta})]^{-1}(\{\sigma\}) = \textsl{C}^{k}\Gamma (\sigma^{\ast}(TF\eta)).
\end{equation}
In this way $\varPhi^{k}$ becomes an isomorphism of vector bundles, of Banach space fibers, over $\textsl{C}^{k}\Gamma (\eta)$:
\begin{equation}\label{69}
\varPhi^{k} :  T\textsl{C}^{k}\Gamma (\eta)\widetilde{\longrightarrow} \textsl{C}^{k}\Gamma (TF\eta),\;\;l(\tau_{\eta}\rvert_{TF\eta})\circ
\varPhi^{k} = \tau_{\textsl{C}^{k}\Gamma (\eta)}.
\end{equation}\\
Remark that, if $E\subseteq\eta$ is a vector bundle neighbourhood in $\eta$, then $TFE$ becomes a vector bundle neighbourhood for $TF\eta$
(again over $D$). Indeed,
\begin{equation}\label{70}
 TFE = (\tau_{\eta}\rvert_{TF\eta})^{-1}(E)
\end{equation}
is open in $TF\eta$ and its fiber $(TFE)_{\zeta} = T(E_{\zeta})$ is the total space of the tagent bundle to the open subset $E_{\zeta}$ of 
$\eta_{\zeta}$, a vector space itself. Recall the bijection $\varPsi^{E} : E\times E\longrightarrow TE$, for an arbitrary vector space $E$ 
(see (\ref{2})), that gives a structure of vector space on $TE$. Therefore 
$\varPsi^{E_{\zeta}} : E_{\zeta}\times E_{\zeta}\longrightarrow T(E_{\zeta}),\;\zeta\in D$, makes $E\times_{D} E$ a $\textsl{C}^{\infty}$ vector 
bundle over $D$ and
\begin{equation}\label{71}
 \varPsi^{E} : E\times_{D} E\widetilde{\longrightarrow} TFE
\end{equation}
a vector bundle isomorphism. In this way $TFE$ is a vector bundle neighbourhood for $TF\eta$ if $E$ is a vector bundle neighbourhood for $\eta$.
And the open subset $\textsl{C}^{k}\Gamma (TFE)$ of $\textsl{C}^{k}\Gamma (TF\eta)$ is also a Banach vector space with the same differentiable 
structure as that induced by the Banach manifold  $\textsl{C}^{k}\Gamma (TF\eta)$. For the vector bundle $E$ we have also the isomorphisms of 
Banach spaces
\begin{equation}\label{72}
  l(\varPsi^{E}) :\textsl{C}^{k}\Gamma(E\times_{D} E)\widetilde{\longrightarrow}\textsl{C}^{k}\Gamma (TFE),\;\;l(\varPsi^{E})(\sigma) = 
  \varPsi^{E}\circ\sigma,
\end{equation}
\begin{equation}\label{73}
 \varSigma^{k} : \textsl{C}^{k}\Gamma(E)\times\textsl{C}^{k}\Gamma(E)\widetilde{\longrightarrow}\textsl{C}^{k}\Gamma(E\times_{D} E),\;\;
 \varSigma^{k}(\sigma,\tau)_{\zeta} = (\sigma_{\zeta},\tau_{\zeta}),
\end{equation}
\begin{equation}\label{74}
 \varPsi^{\textsl{C}^{k}\Gamma(E)} : \textsl{C}^{k}\Gamma(E)\times\textsl{C}^{k}\Gamma(E)\widetilde{\longrightarrow}T\textsl{C}^{k}\Gamma(E),
\end{equation}
\begin{equation}\label{75}
 \varPhi^{k} : T\textsl{C}^{k}\Gamma (E)\widetilde{\longrightarrow}\textsl{C}^{k}\Gamma (TFE),
\end{equation}
that make a commutative diagram:
\begin{equation}\label{76}
 \varPhi^{k}\cdot \varPsi^{\textsl{C}^{k}\Gamma(E)} =  l(\varPsi^{E})\cdot \varSigma^{k}. 
\end{equation}
Now, for every positive integer $r$ we denote
\begin{equation}\label{77}
 \textsl{C}^{\infty}\Gamma_{r} (T\textsl{C}^{k}\Gamma (\eta)) =\{ X\in\textsl{C}^{\infty}(\textsl{C}^{k+r}\Gamma (\eta),
 T\textsl{C}^{k}\Gamma (\eta))\;\rvert\;
 \tau_{\textsl{C}^{k}\Gamma (\eta)}\circ X = \mathrm{id}_{\textsl{C}^{k+r}\Gamma (\eta)}\}.
\end{equation}
We see that $\textsl{C}^{\infty}\Gamma_{0} (T\textsl{C}^{k}\Gamma (\eta)) = \textsl{C}^{\infty}\Gamma (T\textsl{C}^{k}\Gamma (\eta))$ in 
the usual sense of smooth vector fields on the manifold $\textsl{C}^{k}\Gamma (\eta)$. Next we consider
\begin{equation}\label{78}
 \textsl{C}^{\infty}\Gamma_{r} (T\textsl{C}^{\infty}\Gamma (\eta)) = 
 \bigcap_{k=0}^{\infty}\textsl{C}^{\infty}\Gamma_{r} (T\textsl{C}^{k}\Gamma (\eta)).
\end{equation}
Here we keep in mind that $\textsl{C}^{\infty}\Gamma (\eta)$ is dense in $\textsl{C}^{k}\Gamma (\eta)$, that  $\varPhi^{k}$ from (\ref{69})
is independent of $k$ in the sense that
\begin{equation}\label{79}
 \varPhi^{k}\rvert_{T\textsl{C}^{k+1}\Gamma (\eta)} = \varPhi^{k+1},
\end{equation}
for the inclusion $T\textsl{C}^{k+1}\Gamma (\eta)\subseteq T\textsl{C}^{k}\Gamma (\eta)$ is tangent to the inclusion 
$\textsl{C}^{k+1}\Gamma (\eta)\subseteq\textsl{C}^{k}\Gamma (\eta)$. In this way we obtain
\begin{equation}\label{80}
 \varPhi : \bigcap_{k=0}^{\infty} T\textsl{C}^{k}\Gamma (\eta)\longrightarrow\textsl{C}^{\infty}\Gamma (TF\eta)
\end{equation}
and for $X\in\textsl{C}^{\infty}\Gamma_{r} (T\textsl{C}^{\infty}\Gamma (\eta))$ we get
\begin{equation}\label{81}
 \varPhi\circ X : \textsl{C}^{\infty}\Gamma (\eta)\longrightarrow\textsl{C}^{\infty}\Gamma (TF\eta)
\end{equation}
with the property that (see (\ref{69}))
\begin{equation}\label{82}
 l(\tau_{\eta}\rvert_{TF\eta})\circ\varPhi\circ X = \mathrm{id}_{\textsl{C}^{\infty}\Gamma (\eta)}.
\end{equation}
Finally we denote (see(\ref{78}))
\begin{equation}\label{83}
 \textsl{C}^{\infty}\Gamma_{\bullet} (T\textsl{C}^{\infty}\Gamma (\eta)) = 
 \bigcup_{r=0}^{\infty}\textsl{C}^{\infty}\Gamma_{r} (T\textsl{C}^{\infty} \Gamma (\eta)). 
\end{equation}
The elements of $\textsl{C}^{\infty}\Gamma_{r} (T\textsl{C}^{\infty}\Gamma (\eta))$ will be called smooth vector fields on 
$\textsl{C}^{\infty}\Gamma (\eta)$ \textit{of order} $r$ and those from 
$\textsl{C}^{\infty}\Gamma_{\bullet} (T\textsl{C}^{\infty} \Gamma (\eta))$ smooth vector fields \textit{of finite order}. We show now that
$\textsl{C}^{\infty}\Gamma_{\bullet} (T\textsl{C}^{\infty} \Gamma (\eta))$  becomes a graduated Lie algebra, in the sense that for a suitable
bracket
\begin{equation}\label{84}
 [\textsl{C}^{\infty}\Gamma_{r} (T\textsl{C}^{\infty}\Gamma (\eta)),\textsl{C}^{\infty}\Gamma_{s} (T\textsl{C}^{\infty}\Gamma (\eta))]\subseteq
 \textsl{C}^{\infty}\Gamma_{r+s} (T\textsl{C}^{\infty}\Gamma (\eta)),\;\;\forall r, s.
\end{equation}
Let $X\in\textsl{C}^{\infty}\Gamma_{r} (T\textsl{C}^{\infty}\Gamma (\eta))$ and $E$ be a vector bundle neighbourhood in $\eta$. As 
$\textsl{C}^{k}\Gamma (E)$ is open in $\textsl{C}^{k}\Gamma (\eta)$, for $\sigma\in\textsl{C}^{k}\Gamma (E),\;\;
T_{\sigma}\textsl{C}^{k}\Gamma (\eta) = T_{\sigma}\textsl{C}^{k}\Gamma (E)$ so that $X_{\sigma}\in T_{\sigma}\textsl{C}^{k}\Gamma (E)$ for 
$\sigma\in\textsl{C}^{k+r}\Gamma (E)\subseteq\textsl{C}^{k}\Gamma (E)$, we find that $X$ maps $\textsl{C}^{k+r}\Gamma (E)$ into  
$T\textsl{C}^{k}\Gamma (E)$ for all $k$. If we denote
\begin{equation}\label{85}
 p_{2}^{k} : \textsl{C}^{k}\Gamma (E)\times\textsl{C}^{k}\Gamma (E)\longrightarrow\textsl{C}^{k}\Gamma (E)
\end{equation}
the canonical projection, then the representation of $X$ in the chart defined by $E$ would be
\begin{equation}\label{86}
 \tilde{X}^{k} =  p_{2}^{k}\circ (\varPsi^{\textsl{C}^{k}\Gamma (E)})^{-1} \circ X.
\end{equation}
Since $\varPsi^{\textsl{C}^{k}\Gamma (E)}$, as well as $ p_{2}^{k}$, are independent of $k$ in the sense already defined (see(\ref{79})), it 
follows that
$\tilde{X}^{k}\rvert_{\textsl{C}^{k+r+1}\Gamma (E)} = \tilde{X}^{k+1},\;\;\forall k\geqslant 0$.
We get $\tilde{X} : \textsl{C}^{\infty}\Gamma (E)\longrightarrow\textsl{C}^{\infty}\Gamma (E)$ such that 
$\tilde{X}\in\textsl{C}^{\infty}(\textsl{C}^{k+r}\Gamma (E),\textsl{C}^{k}\Gamma (E)),\;\;\forall k\geqslant 0$. In fact, from (\ref{76}) and
(\ref{86}) we get
\begin{equation}\label{87}
 \varPhi (X_{\sigma})_{\zeta} = \varPsi^{E_{\zeta}}(\sigma_{\zeta},\tilde{X}(\sigma)_{\zeta}).
\end{equation}
Analogously, for
$Y\in\textsl{C}^{\infty}\Gamma_{s} (T\textsl{C}^{\infty}\Gamma (\eta))$ we get 
$\tilde{Y}\in\bigcap_{k=0}^{\infty}\textsl{C}^{\infty}(\textsl{C}^{k+s}\Gamma (E),\textsl{C}^{k}\Gamma (E))$. Then for 
$\sigma\in\textsl{C}^{k+r+s}\Gamma (E)$ 
\begin{equation}\label{88}
 \big[\tilde{X},\tilde{Y}\big](\sigma) = <\tilde{Y}^{\prime}(\sigma);\tilde{X}(\sigma)> - <\tilde{X}^{\prime}(\sigma);\tilde{Y}(\sigma)>
\end{equation}
is well defined in $\textsl{C}^{k}\Gamma (E)$ and an easy inspection shows that in fact\\ 
\begin{equation}\label{89}
 \big[\tilde{X},\tilde{Y}\big]\in\bigcap_{k=0}^{\infty}\textsl{C}^{\infty}(\textsl{C}^{k+r+s}\Gamma (E),\textsl{C}^{k}\Gamma (E)).
\end{equation}
It remains to verify that, if we define
\begin{equation}\label{90}
\big[X,Y\big]_{\sigma} =  \varPsi^{\textsl{C}^{k}\Gamma (E)}(\sigma,\big[\tilde{X},\tilde{Y}\big](\sigma))
\end{equation}
for $\sigma\in\textsl{C}^{k+r+s}\Gamma (E)$, the definition does not depend on the vector bundle neighbourhood $E\supseteq \sigma (D)$. Let then
$E_{1},\;E_{2}$ be two such neighbourhoods of $\sigma$ in $\eta$. We need to distingush $\vartheta_{j} : E_{j}\longrightarrow\eta,\;j = 1, 2$,
the respective fiber preserving embeddings and consider $\psi = \vartheta_{2}^{-1}\circ\vartheta_{1}$ the nonlinear fiber preserving mapping
from an open subset $O_{1}\subseteq E_{1}$ onto an open subset  $O_{2}\subseteq E_{2}$. Then $\textsl{C}^{k}\Gamma (O_{j})$ are open in
 $\textsl{C}^{k}\Gamma (E_{j})$ and
 $l(\psi) : \textsl{C}^{k}\Gamma (O_{1})\longrightarrow\textsl{C}^{k}\Gamma (O_{2}),\\
 l(\psi)(\sigma) = \psi\circ\sigma$, is a diffeomorphism for each $k\geqslant 0$. $ l(\psi)$ is precisely the map that changes 
coordinates, in $\textsl{C}^{k}\Gamma (\eta)$, from $\textsl{C}^{k}\Gamma (E_{1})$ into coordinates from $\textsl{C}^{k}\Gamma (E_{2})$. But
because $l(\psi)$ does not depend on $k$, all the computations go on as if all happens in the Banach spaces 
$\textsl{C}^{k}\Gamma (E_{j}),\;\;j =1, 2$, for one fixed $k$. These considerations dealt with the graduated Lie algebra (\ref{83}) of 
smooth vector fields on $\textsl{C}^{\infty}\Gamma (\eta)$ of finite order, for $\eta\longrightarrow D$ deneral fiber bundle.\\
Now let us consider the case when $\eta = D\times M$ is the trivial nonlinear fiber bundle $p_{D} : D\times M\longrightarrow D$. Then the 
identification of $\sigma\in \textsl{C}^{k}\Gamma ( D\times M)$ with $\beta\in \textsl{C}^{k}(D,M)$ given by
\begin{equation}\label{91}
 \textsl{C}^{k}\Gamma ( D\times M)\widetilde{\longrightarrow}\textsl{C}^{k}(D,M),\;\; \sigma_{\zeta} = (\zeta, \beta (\zeta)),\;\;\zeta\in D,
\end{equation}
and that of 
\begin{equation}\label{92}
 TF(D\times M)\widetilde{\longrightarrow} D\times TM,\;\;T_{(\zeta,p)}(\{\zeta\}\times M)\widetilde{\longrightarrow} \{\zeta\}\times T_{p} M,
\end{equation}
lead to the isomorphism of vector bundles over $\textsl{C}^{k}(D,M)$ (see (\ref{65}), (\ref{69})):
\begin{equation}\label{93}
 \varPhi^{k} :  T\textsl{C}^{k} (D,M)\widetilde{\longrightarrow} \textsl{C}^{k}(D,TM),
\end{equation}
\begin{equation}\label{94}
  \varPhi^{k}_{\beta} :  T_{\beta}\textsl{C}^{k} (D,M)\widetilde{\longrightarrow} \textsl{C}^{k}\Gamma (\beta^{\ast}(TM)),
\end{equation}
\begin{equation}\label{95}
[\varPhi^{k}_{\beta (0)}(\dfrac{\mathrm{d}\beta(t)}{\mathrm{d}t}\big |_{t=0})]_{\zeta} = 
\dfrac{\mathrm{d}[\beta(t)(\zeta)]}{\mathrm{d}t}\big |_{t=0},\;\zeta\in D.
\end{equation}
Still, in order to introduce charts on $\textsl{C}^{k} (D,M)$ in the neighbourhood of $\beta_{0}$, we have to consider vector bundle 
neighbourhoods $E$ of $\mathrm{graph}(\beta_{0})$ in $D\times M$. If $\chi$ is the chart defined by $E$, then as in (\ref{91})
\begin{equation}\label{96}
 \chi (\beta) = \sigma,\;\;\sigma_{\zeta} = (\zeta,\beta (\zeta)),\;\;\zeta\in D.
\end{equation}
We will then denote
\begin{equation}\label{97}
 \textsl{C}^{\infty}\Gamma_{\bullet} (T\textsl{C}^{\infty}(D,M))
\end{equation}
the graduated Lie algebra of smooth vector fields on $\textsl{C}^{\infty}(D,M)$ of finite order.\\
There are two important instances of smooth vector fields of finite order on $\textsl{C}^{\infty}(D,M)$. First, for 
$X\in \textsl{C}^{\infty}\Gamma (TM)$
\begin{equation}\label{98}
 l(X)_{\beta}(\zeta) =: X_{\beta(\zeta)},\;\;\beta\in\textsl{C}^{\infty}(D,M), \;\;\zeta\in D,
\end{equation}
and second, for $X\in \textsl{C}^{\infty}\Gamma (TD)$
\begin{equation}\label{99}
 r(X)_{\beta}(\zeta) =: T_{\zeta}\beta\cdot X_{\zeta},\;\;\beta\in\textsl{C}^{\infty}(D,M), \;\;\zeta\in D.
\end{equation}
$l(X)$ is of order $0$, while $r(X)$ is of order $1$. And we have
\begin{equation}\label{100}
 l([X,Y]) = [l(X),l(Y)],\;\;\forall X,\;Y\in\textsl{C}^{\infty}\Gamma (TM),
\end{equation}
\begin{equation}\label{101}
 r([X,Y]) = - [r(X),r(Y)],\;\;\forall X,\;Y\in\textsl{C}^{\infty}\Gamma (TD).
\end{equation}
The first equality comes from the Lie group morphism

 $L : \textsl{C}^{\infty}\mathrm{Diff} M\longrightarrow \textsl{C}^{\infty}\mathrm{Diff}\textsl{C}^{k} (D,M),\;
 L(\varphi)(\beta) = \varphi\circ\beta,\\
 \varphi\in \textsl{C}^{\infty}\mathrm{Diff} M,\;\;\beta\in\textsl{C}^{k} (D,M)$,

while from the Lie group antimorphism

 $R :  \textsl{C}^{\infty}\mathrm{Diff} D\longrightarrow  \textsl{C}^{\infty}\mathrm{Diff}\textsl{C}^{k} (D,M),\;
 R(\varphi)(\beta) = \beta\circ\varphi,\\
 \varphi\in \textsl{C}^{\infty}\mathrm{Diff} D,\;\;\beta\in\textsl{C}^{k} (D,M)$,

comes the second equality. Indeed, we have

$l = T_{\mathrm{id}_{M}} L,\;\;r =  T_{\mathrm{id}_{D}} R$

as the Lie algebra morphism, or antimorphism, corresponding to them.\\
According to (\ref{93}) and (\ref{67}) we will understand
\begin{equation}\label{102}
 l(\tau_{M}) : \textsl{C}^{k} (D,TM)\longrightarrow\textsl{C}^{k} (D,M)
\end{equation}
as a vector bundle, $\forall k\geqslant 0$. A more general case is the following: for any $\pi : G\longrightarrow M$, smooth vector bundle,
$l(\pi) : \textsl{C}^{k} (D,G)\longrightarrow\textsl{C}^{k} (D,M)$ is the smooth vector bundle of Banach space fiber
\begin{equation}\label{103}
 \textsl{C}^{k} (D,G)_{\beta} = \textsl{C}^{k}\Gamma (\beta^{\ast}(G)).
\end{equation}
When $\sigma : M\longrightarrow G$ is a section of $\pi$, $l(\sigma)$ is a section of $l(\pi)$, if $l(\sigma)(\beta) = 
\sigma\circ\beta$. 
We denote, as in (\ref{77})

 $\textsl{C}^{\infty}\Gamma_{r} (\textsl{C}^{k}(D,G)) =\{ X\in\textsl{C}^{\infty}(\textsl{C}^{k+r}(D,M),
 \textsl{C}^{k}(D,G))\;\rvert\;
l(\pi)\circ X = \mathrm{id}_{\textsl{C}^{k+r}(D,M)}\}$. 

Then 
\begin{equation}\label{104}
  \textsl{C}^{\infty}\Gamma_{\bullet} (\textsl{C}^{\infty}(D,G)) = 
  \bigcup_{r=0}^{\infty}\bigcap_{k=0}^{\infty}\textsl{C}^{\infty}\Gamma_{r} (\textsl{C}^{k}(D, G))
\end{equation}
will be the space of smooth sections, of finite order, for the vector bundle
\begin{equation}\label{105}
 l(\pi) : \textsl{C}^{\infty} (D,G)\longrightarrow\textsl{C}^{\infty} (D,M).
\end{equation}
We call (\ref{105}) \textit{the lift of the vector bundle} $\pi : G\longrightarrow M$ \textit{to} $\textsl{C}^{\infty} (D,M)$. 
Of special interest is the case when $G = H\subseteq TM$ is the vector subbundle under study and also the case when $G = TM/H$. 
We remark that, if $H\subseteq TM$ is a vector subbundle, then $\textsl{C}^{k}(D,H)$ is a vector subbundle, over $\textsl{C}^{k}(D,M)$, of 
$\textsl{C}^{k}(D,TM),\;\;\forall k\geqslant 0$. In the case that $P$ is a smooth section of $T^{\ast} M\otimes TM$ such that $P_{m}$ is a 
projection on $H_{m},\;\;\forall m\in M$ (as, for instance, when $P_{m}$ is the orthogonal projection on $H_{m}$ defined by a smooth metric 
on $M$) we may consider for $X\in \textsl{C}^{\infty}\Gamma_{\bullet} (\textsl{C}^{\infty}(D,TM))$
\begin{equation}\label{106}
 l(P)\cdot X\in\textsl{C}^{\infty}\Gamma_{\bullet} (\textsl{C}^{\infty}(D,H))
\end{equation}
defined $\forall \beta\in\textsl{C}^{\infty}(D,M)$ and $\forall\zeta\in D$ by
\begin{equation}\label{107}
  (l(P)\cdot X)_{\beta}(\zeta) = l(P)_{\beta}(\zeta)\cdot X_{\beta}(\zeta) = P_{\beta (\zeta)}\cdot X_{\beta}(\zeta).
\end{equation}
Using such a projection (\ref{106}), we see that the isomorphism
\begin{equation}\label{108}
TM \widetilde{\longrightarrow} H\times_{M} TM/H
\end{equation}
entails the isomorphism
\begin{equation}\label{109}
\textsl{C}^{\infty} (D,TM) \widetilde{\longrightarrow} \textsl{C}^{\infty}(D,H)\times_{\textsl{C}^{\infty}(D,M)}\textsl{C}^{\infty}(D,TM/H)
\end{equation}
and then also that
\begin{equation}\label{110}
Q_{\beta} := \textsl{C}^{\infty} (D,TM)_{\beta}/\textsl{C}^{\infty}(D,H)_{\beta}\widetilde{\longrightarrow}\textsl{C}^{\infty}(D,TM/H)_{\beta},
\end{equation}
$\forall\beta\in\textsl{C}^{\infty}(D,M)$.
Let us recall the vector bundle morphism (\ref{21}) $P^{H} :TM\longrightarrow TM/H$. Analogously, we can consider
\begin{equation}\label{111}
P^{\textsl{C}^{\infty}(D,H)} : \textsl{C}^{\infty} (D,TM)\longrightarrow \textsl{C}^{\infty} (D,TM)/ \textsl{C}^{\infty}(D,H).
\end{equation}
If we define $l(P^{H}) : \textsl{C}^{\infty} (D,TM)\longrightarrow\textsl{C}^{\infty}(D,TM/H)$ by $l(P^{H})(X) = P^{H}\circ X$ we see that
\begin{equation}\label{112}
Q_{\beta}\cdot P^{\textsl{C}^{\infty}(D,H)}_{\beta} = l(P^{H})_{\beta}.
\end{equation}
Taking into account the equality (\ref{112}) and the formula (\ref{20}), the following result gives the curvature of the lift to 
$\textsl{C}^{\infty}(D,M)$ in terms of the curvature of the subbundle. It is also in agreement with the relation (\ref{100}).

\begin{theorem}. Let $H\subseteq TM$ be a smooth vector subbundle, over $M$, of curvature tensor $C^{H}$ and $D$ 
be a smooth compact manifold. Then 
\begin{equation}\label{113}
 P^{H}_{\beta (\zeta)} \big[X,Y\big]_{\beta} (\zeta) = C^{H}_{\beta(\zeta)} (X_{\beta} (\zeta)\wedge Y_{\beta} (\zeta)),
\end{equation}
$\forall\;X,\;Y\in \textsl{C}^{\infty}\Gamma_{\bullet} (\textsl{C}^{\infty}(D,H)),\;\;\forall\; \beta\in\textsl{C}^{\infty}(D,M),\;\;\forall\; 
\zeta\in D$, the Lie bracket from the left being taken as for smooth vector fields of finite order on $\textsl{C}^{\infty}(D,M)$.
\end{theorem}
\textbf{Proof}. We intend to prove (\ref{113}) for a fixed $\beta_{0}$ and a fixed $\zeta_{0}$. Let then $E$ be a vector bundle neighbourhood of
$\mathrm{graph} (\beta_{0})$ in $D\times M$ and $\iota : E\longrightarrow D\times M$ be the respective fiber preserving open embedding of it. If
$\pi : E\longrightarrow D$ is the respective structural projection and $p_{2} : D\times M\longrightarrow M$ is canonical, then for
\begin{equation}\label{114}
 \kappa := p_{2} \circ \iota,\;\;\kappa : E\longrightarrow M,
\end{equation}
we have $\iota (e) = (\pi (e),\kappa (e)),\;\;\forall\;e\in E$. And for  every $\zeta\in D$
\begin{equation}\label{115}
 \kappa_{\zeta}:= \kappa\rvert_{E_{\zeta}},\;\; \kappa_{\zeta} : E_{\zeta}\longrightarrow M,
\end{equation}
is a diffeomorphism on an open subset of $M$. Let $k\geqslant 0$ be fixed and let
\begin{equation}\label{116}
 O^{k} = \{\beta\in\textsl{C}^{k}(D,M)\;\rvert\;\mathrm{graph} (\beta )\subset E\}
\end{equation}
be the open domain of the chart on $\textsl{C}^{k}(D,M)$:
\begin{equation}\label{117}
 \chi : O^{k}\longrightarrow \textsl{C}^{k}\Gamma (E),\;\;\chi (\beta)_{\zeta} = (\zeta , \beta (\zeta)),\;\forall\;\zeta\in D,
\end{equation}
(see also (\ref{96})). Therefore
\begin{equation}\label{118}
 \chi^{-1} (\sigma) = \kappa \circ \sigma,\;\;\forall\; \sigma\in\textsl{C}^{k}\Gamma (E),
\end{equation}
or, for a fixed $\zeta_{0}$
\begin{equation}\label{119}
 \chi (\beta)_{\zeta_{0}} = \kappa_{\zeta_{0}}^{-1} (\beta (\zeta_{0})),\;\;\forall\;\beta\in O^{k}.
\end{equation}
In what follows $\kappa_{\zeta_{0}}^{-1}$ \textit{will be considered as an} $E_{\zeta_{0}}$-\textit{valued chart on} $M$ and will be denoted
\begin{equation}\label{120}
 \mathring{\chi} := \kappa_{\zeta_{0}}^{-1}.
\end{equation}
Then (\ref{119}) reads
\begin{equation}\label{121}
 \chi (\beta)_{\zeta_{0}} = \mathring{\chi} (\beta (\zeta_{0})),\;\;\forall\;\beta\in O^{k}.
\end{equation}
Let us consider the map 
\begin{equation}\label{122}
 \delta_{\zeta_{0}} : \textsl{C}^{k}(D,M)\longrightarrow M,\;\; \delta_{\zeta_{0}} (\beta) = \beta (\zeta_{0}).
\end{equation}
It is easy to see that (see (\ref{94}))
\begin{equation}\label{123}
 T_{\beta} \delta_{\zeta_{0}}\cdot X = X_{\zeta_{0}},\;\;\forall\;X\in \textsl{C}^{k}\Gamma (\beta^{\ast}(TM)).
\end{equation}
Taking now the differential with respect to $\beta$ in (\ref{121}) we get
\begin{equation}\label{124}
 (\mathrm{d}_{\beta}\chi\cdot X)_{\zeta_{0}} = \mathrm{d}_{\beta (\zeta_{0})} \mathring{\chi}\cdot X_{\zeta_{0}},
 \;\;\forall X\in \textsl{C}^{k}\Gamma (\beta^{\ast}(TM)),
\end{equation}
since the mapping $:\textsl{C}^{k}\Gamma (E)\longrightarrow E_{\zeta_{0}},\;\;\sigma\mapsto \sigma_{\zeta_{0}}$, is linear and continuous. On
the other hand, for a smooth vector field $X$ of order $r$ on $\textsl{C}^{\infty}(D,M)$, if $\widetilde{\chi_{\ast} X}$ is defined by
\begin{equation}\label{125}
 ((\chi^{-1})^{\ast} X)_{\sigma} = (\sigma, \widetilde{\chi_{\ast} X} (\sigma)),\;\;\forall\;\sigma\in\textsl{C}^{k+r}\Gamma (E), 
\end{equation}
then $\widetilde{\chi_{\ast} X} : \textsl{C}^{k+r}\Gamma (E)\longrightarrow\textsl{C}^{k}\Gamma (E)$ as a $\textsl{C}^{\infty}$ map (see(\ref{86})
and (\ref{87})). Also, as in (\ref{29})
\begin{equation}\label{126}
 \widetilde{\chi_{\ast} X} (\sigma)= \mathrm{d}_{\chi^{-1} (\sigma)}\chi\cdot X_{\chi^{-1} (\sigma)},\;\;
 \forall\;\sigma\in\textsl{C}^{k+r}\Gamma (E).
\end{equation}
Combining this with (\ref{124}) and (\ref{121}) we come to
\begin{equation}\label{127}
 \widetilde{\chi_{\ast} X} (\sigma)_{\zeta_{0}} = 
 \mathrm{d}_{\mathring{\chi}^{-1}(\sigma_{\zeta_{0}})}\mathring{\chi} \cdot (X_{\chi^{-1} (\sigma)})(\zeta_{0}),
 \;\;\forall\;\sigma\in\textsl{C}^{k+r}\Gamma (E).
\end{equation}
Next, for $X, Y\in\textsl{C}^{\infty}\Gamma_{\bullet} (T\textsl{C}^{\infty}(D,M))$ (see (\ref{97}))
\begin{equation}\label{128}
 \widetilde{\chi_{\ast}[X,Y]} = [\widetilde{\chi_{\ast} X},\widetilde{\chi_{\ast} Y}]
\end{equation}
(compare with (\ref{32}) as a consequence of the definition (\ref{88}) for $[\tilde{X},\tilde{Y}]$ and (\ref{90}) for $[X,Y]$. Then from 
(\ref{127}) we get
\begin{equation}\label{129}
  [\widetilde{\chi_{\ast} X},\widetilde{\chi_{\ast} Y}](\sigma)_{\zeta_{0}} = 
  \mathrm{d}_{\mathring{\chi}^{-1}(\sigma_{\zeta_{0}})}\mathring{\chi}\cdot ([X,Y]_{\chi^{-1} (\sigma)})(\zeta_{0}),
 \;\;\forall\;\sigma\in\textsl{C}^{k+r+s}\Gamma (E),
\end{equation}
(if $Y$ is of order $s$) and finally
\begin{multline}\label{130}
 \mathrm{d}_{\beta_{0}(\zeta_{0})}\mathring{\chi}\cdot ([X,Y]_{\beta_{0}})(\zeta_{0}) =\\
 = \dfrac{\mathrm{d}}{\mathrm{d} t}[\widetilde{\chi_{\ast} Y}(\chi (\beta_{0}) + t\;\widetilde{\chi_{\ast} X} (\chi (\beta_{0})))_{\zeta_{0}} -
 \widetilde{\chi_{\ast} X}(\chi (\beta_{0}) + t\;\widetilde{\chi_{\ast} Y} (\chi (\beta_{0})))_{\zeta_{0}}]\rvert_{t=0}.
\end{multline}
Let us consider the subspaces of $E_{\zeta_{0}}$ (see(\ref{4}))
\begin{equation}\label{131}
 \widetilde{\mathring{\chi}_{\ast} H}_{\mathring{\chi} (m)} := \mathrm{d}_{m}\mathring{\chi}\cdot H_{m},\;\;
 m\in \kappa_{\zeta_{0}} (E_{\zeta_{0}}).
\end{equation}
In virtue of the hypothesis that $X,\;Y\in\textsl{C}^{\infty}\Gamma_{\bullet} (\textsl{C}^{\infty}(D,H))$ we have\\
$(X_{\beta})(\zeta_{0}),\;(Y_{\beta})(\zeta_{0})\in H_{\beta(\zeta_{0})},\;\;\forall\;\beta\in\textsl{C}^{\infty}(D,M).$ It results from 
(\ref{127}) that 
\begin{equation}\label{132}
 \widetilde{\chi_{\ast} X} (\sigma)_{\zeta_{0}},\;\widetilde{\chi_{\ast} Y} (\sigma)_{\zeta_{0}}\in
 \widetilde{\mathring{\chi}_{\ast} H}_{\sigma_{\zeta_{0}}},\;\;\forall\;\sigma\in\textsl{C}^{\infty}\Gamma (E).
\end{equation}
Then, for $\sigma$ in a neighbourhood of $\chi (\beta_{0})$ in $\textsl{C}^{\infty}\Gamma (E)$, $\sigma_{\zeta_{0}}$ remains in a given
neighbourhood $U$ of
\begin{equation}\label{133}
  \chi (\beta_{0})_{\zeta_{0}} = \mathring{\chi} (\beta_{0} (\zeta_{0}))
\end{equation}
(see(\ref{121})). Let us denote
\begin{equation}\label{134}
 V = \widetilde{\mathring{\chi}_{\ast} H}_{\mathring{\chi}(\beta_{0}(\zeta_{0}))}
\end{equation}
and choose $W\subset E_{\zeta_{0}}$ a supplementary for $V$ subspace, i.e. such that $E_{\zeta_{0}} = V\oplus W$. Ignoring in notation a natural 
isomorphism, we will write simply
\begin{equation}\label{135}
E_{\zeta_{0}} = V\times W.
\end{equation}
Next, we choose the neighbourhood $U$ of $\chi (\beta_{0})_{\zeta_{0}}$ such that
\begin{equation}\label{136}
 \widetilde{\mathring{\chi}_{\ast} H}_{(x,y)} = \mathrm{graph} (\varGamma (x,y))
\end{equation}
for $(x,y)\in U$ and certain $\varGamma (x,y)\in\mathrm{Hom} (V,W)$ (compare with (\ref{6})). In this way $\mathring{\chi}$ becomes an 
adapted  to $H$ chart in the neighbourhood of $\beta_{0}(\zeta_{0})$. In order to simplify notation we denote
\begin{equation}\label{137}
 (x_{0},y_{0}) = \mathring{\chi} (\beta_{0} (\zeta_{0}))
\end{equation}
and then (see (\ref{134})) $\varGamma (x_{0},y_{0}) = 0$.  We will put also
\begin{equation}\label{138}
 \sigma_{0}:= \chi (\beta_{0})
\end{equation}
and then
\begin{equation}\label{139}
 (\sigma_{0})_{\zeta_{0}} =  (x_{0},y_{0}).
\end{equation}
We recall the notation $P^{\widetilde{\mathring{\chi}_{\ast} H}}_{e}$ from (\ref{35}) and the relation (\ref{36}).
Then, if we apply $P^{\widetilde{\mathring{\chi}_{\ast} H}}_{(x_{0},y_{0})}$ on both sides of (\ref{130}), we get 
\begin{multline}\label{140}
  \mathrm{d}_{\beta_{0} (\zeta_{0})}\mathring{\chi}/H_{\beta_{0} (\zeta_{0})}\cdot P^{H}_{\beta_{0} (\zeta_{0})}
([X,Y]_{\beta_{0}})(\zeta_{0}) =\\
=P^{\widetilde{\mathring{\chi}_{\ast} H}}_{(x_{0},y_{0})} \dfrac{\mathrm{d}}{\mathrm{d} t}[\widetilde{\chi_{\ast} Y}(\sigma_{0} +
t\;\widetilde{\chi_{\ast} X} (\sigma_{0}))_{\zeta_{0}} -
 \widetilde{\chi_{\ast} X}(\sigma_{0} + t\;\widetilde{\chi_{\ast} Y} (\sigma_{0}))_{\zeta_{0}}]\rvert_{t=0}.
\end{multline}
From (\ref{132}) and (\ref{136}) we infer the existence of smooth functions $v_{X},\;v_{Y}$ on a neighbourhood of $\sigma_{0}$ in 
$\textsl{C}^{\infty}\Gamma (E)$ taking values in $V$, such that
\begin{equation}\label{141}
 \widetilde{\chi_{\ast} X} (\sigma)_{\zeta_{0}} = (v_{X} (\sigma),\varGamma (\sigma_{\zeta_{0}})\cdot v_{X} (\sigma))\in V\times W,
\end{equation}
and analogously for $\widetilde{\chi_{\ast} Y} (\sigma)_{\zeta_{0}}$. Then the right hand side of (\ref{140}) is\\

$P^{\widetilde{\mathring{\chi}_{\ast} H}}_{(x_{0},y_{0})} 
[<\widetilde{\chi_{\ast} Y}^{\prime}(\sigma_{0})_{\zeta_{0}};\widetilde{\chi_{\ast} X} (\sigma_{0})> -
<\widetilde{\chi_{\ast} X}^{\prime}(\sigma_{0})_{\zeta_{0}};\widetilde{\chi_{\ast} Y} (\sigma_{0})>] = $\\

$= P^{\widetilde{\mathring{\chi}_{\ast} H}}_{(x_{0},y_{0})}\big ([<v^{\prime}_{Y} (\sigma_{0});\widetilde{\chi_{\ast} X} (\sigma_{0})> -
 <v^{\prime}_{X} (\sigma_{0});\widetilde{\chi_{\ast} Y} (\sigma_{0})>],$\\
 
$\varGamma ((\sigma_{0})_{\zeta_{0}})[<v^{\prime}_{Y} (\sigma_{0});\widetilde{\chi_{\ast} X} (\sigma_{0})> -
 <v^{\prime}_{X} (\sigma_{0});\widetilde{\chi_{\ast} Y} (\sigma_{0})>]\big) +$\\

$+ P^{\widetilde{\mathring{\chi}_{\ast} H}}_{(x_{0},y_{0})}\big (0,<\varGamma^{\prime}(x_{0},y_{0});
\widetilde{\chi_{\ast} X} (\sigma_{0})_{\zeta_{0}}> v_{Y} (\sigma_{0})\; - 
<\varGamma^{\prime}(x_{0},y_{0});\widetilde{\chi_{\ast} Y} (\sigma_{0})_{\zeta_{0}}> v_{X} (\sigma_{0})\big ) =$\\

$=  P^{\widetilde{\mathring{\chi}_{\ast} H}}_{(x_{0},y_{0})}\big (0,<\dfrac{\partial\varGamma}{\partial x} (x_{0},y_{0}); v_{X} (\sigma_{0})> 
v_{Y}(\sigma_{0}) -
<\dfrac{\partial\varGamma}{\partial x} (x_{0},y_{0}); v_{Y} (\sigma_{0})> v_{X}(\sigma_{0}) +$\\

$+ <\dfrac{\partial\varGamma}{\partial y} (x_{0},y_{0}); \varGamma (x_{0},y_{0}) v_{X} (\sigma_{0})> v_{Y}(\sigma_{0}) -
<\dfrac{\partial\varGamma}{\partial y} (x_{0},y_{0}); \varGamma (x_{0},y_{0}) v_{Y} (\sigma_{0})> v_{X}(\sigma_{0})\big ) =$\\

$= \varLambda^{W}_{(x_{0},y_{0})}\big (<\dfrac{\partial\varGamma}{\partial x} (x_{0},y_{0}); v_{X} (\sigma_{0})> 
v_{Y}(\sigma_{0}) -
<\dfrac{\partial\varGamma}{\partial x} (x_{0},y_{0}); v_{Y} (\sigma_{0})> v_{X}(\sigma_{0}) +$\\

$+ <\dfrac{\partial\varGamma}{\partial y} (x_{0},y_{0}); \varGamma (x_{0},y_{0}) v_{X} (\sigma_{0})> v_{Y}(\sigma_{0}) -
<\dfrac{\partial\varGamma}{\partial y} (x_{0},y_{0}); \varGamma (x_{0},y_{0}) v_{Y} (\sigma_{0})> v_{X}(\sigma_{0})\big ) =$\\

$=  \mathrm{d}_{\beta_{0} (\zeta_{0})}\mathring{\chi}/H_{\beta_{0} (\zeta_{0})}\cdot C^{H}_{\beta_{0} (\zeta_{0})}
((\mathrm{d}_{\beta_{0} (\zeta_{0})}\mathring{\chi})^{-1}\varLambda^{V}_ {\mathring{\chi} (\beta_{0} (\zeta_{0}))}\cdot v_{X} (\sigma_{0})
\wedge(\mathrm{d}_{\beta_{0} (\zeta_{0})}\mathring{\chi})^{-1}\varLambda^{V}_ {\mathring{\chi} (\beta_{0} (\zeta_{0}))}\cdot v_{Y} (\sigma_{0}))
 =$\\

$=  \mathrm{d}_{\beta_{0} (\zeta_{0})}\mathring{\chi}/H_{\beta_{0} (\zeta_{0})}\cdot C^{H}_{\beta_{0} (\zeta_{0})}
(X_{\beta_{0}} (\zeta_{0})\wedge Y_{\beta_{0}} (\zeta_{0}))$\\

in virtue of (\ref{141}), (\ref{127}) and (\ref{40}). Comparing with the left hand side of (\ref{140}) we obtain (\ref{113}) in $\beta_{0}$ and 
$\zeta_{0}$. The proof is complete $\blacksquare$\\

\section{On a localization property}

In the Theorem 4, \S 5, we state an equality of the form $\big[ A,B\big] = C$ in the Lie algebra of vector fields of finite order over 
$\textsl{C}^{\infty} (D,M)$, for $A, B, C$ certain vector fields of a simpler form. The equality is easy to prove for $D$ and $M$ diffeomorphic 
to open subsets of vector spaces, but the argument for the possibility to localize the equality, even in this particular case, is quite involved. 
We preferred to allocate the entire section \S 4 to this question and to base our proof on the fact that $A, B, C$ above belong to a Lie 
subalgebra of $\textsl{C}^{\infty}\Gamma_{\bullet} (T\textsl{C}^{\infty}(D,M))$, isomorphic to the set of global sections of a sheaf of Lie 
algebras over $D\times M$. On the other hand, the formula of the Lie bracket in this subalgebra, from Theorem 3 below, will be used further in 
section \S 6. But we presume that this Lie subalgebra and sheaf are interesting in themselves.\\
Coming back to the general case of a nonlinear fiber bundle $\pi : \eta\longrightarrow D$, 
$\zeta_{0}\in D$ and $\sigma\in\textsl{C}^{r}\Gamma (\eta),\;\;r\geqslant 0$, we denote $j^{r}_{\zeta_{0}}\sigma$ the $r$-jet of $\sigma$ at 
$\zeta_{0}$ and $J^{r}_{\zeta_{0}}\eta = \{j^{r}_{\zeta_{0}}\sigma\;\rvert\;\sigma$ local section of $\eta$ around $\zeta_{0}\}$ (see Palais [5] 
for definitions). The disjoint union $J^{r}\eta$ becomes a fiber bundle over $D$ of fiber $J^{r}_{\zeta}\eta$ at $\zeta\in D$ and $J^{0}\eta =
\eta$. On the other hand the projection
\begin{equation}\label{142}
 \pi^{r}_{0} : J^{r}\eta\longrightarrow\eta,\;\;\pi^{r}_{0}(j^{r}_{\zeta}\sigma) = \sigma_{\zeta},
\end{equation}
makes $J^{r}\eta$ a fiber bundle over $\eta$ of fiber
\begin{equation}\label{143}
 (J^{r}\eta)_{e} = \{j^{r}_{\zeta}\sigma\;\rvert\;\zeta = \pi (e),\;\sigma_{\zeta} = e\},\;\;e\in \eta.
\end{equation}
If 
\begin{equation}\label{144}
 \pi^{r} : J^{r}\eta\longrightarrow D
\end{equation}
denotes the canonical projection, then $\pi^{0} = \pi$ and $\pi^{r} = \pi\circ\pi^{r}_{0}$. For $\sigma\in\textsl{C}^{r}\Gamma (\eta)$, 
$j^{r}\sigma$ becomes a section of $J^{r}\eta$ by $(j^{r}\sigma)_{\zeta} = j^{r}_{\zeta}\sigma,\;\zeta\in D$, such that $\forall\; k\geqslant 0$
\begin{equation}\label{145}
 j^{r} : \textsl{C}^{r+k}\Gamma (\eta)\longrightarrow\textsl{C}^{k}\Gamma (J^{r}\eta)
\end{equation}
is a $\textsl{C}^{\infty}$ map. Then if $\varXi : J^{r}\eta\longrightarrow TF\eta$ is a  $\textsl{C}^{\infty}$ fiber preserving map of fiber 
bundles over $\eta$, i.e.
\begin{equation}\label{146}
 \varXi_{e} : (J^{r}\eta)_{e}\longrightarrow T_{e}\eta_{\pi (e)},\;\;e\in\eta,
\end{equation}
then
\begin{equation}\label{147}
 X_{\sigma}(\zeta) := \varXi_{\sigma_{\zeta}}(j^{r}_{\zeta}\sigma ),\;\;\zeta\in D,\;\sigma\in\textsl{C}^{r+k}\Gamma (\eta),
\end{equation}
is a vector field of order $r$ on  $\textsl{C}^{\infty}\Gamma (\eta )$, i.e.
$X\in\textsl{C}^{\infty}\Gamma_{r} (T\textsl{C}^{\infty}\Gamma (\eta )$. In this case we say that $X$ is a vector field \textit{of differential 
order $r$}. We denote $\textsl{C}^{\infty}FB_{\eta}(J^{r}\eta,TF\eta)$ the set of smooth maps 
$\varXi$ of the form (\ref{146}) (Since $\varXi : J^{r}\eta\longrightarrow TF\eta$ is also fiber preserving map over $D$, $X$ defined by 
(\ref{147}) is a nonlinear differential operator of order $r$ from $\textsl{C}^{\infty}\Gamma (\eta )$ to $\textsl{C}^{\infty}\Gamma(TF\eta )$
in the sense of Palais [5]). Let next, for $0\leqslant s\leqslant r$ (compare with (\ref{142}))
\begin{equation}\label{148}
 \pi^{r}_{s} : J^{r}\eta\longrightarrow J^{s}\eta,\;\; \pi^{r}_{s}(j^{r}_{\zeta}\sigma) = j^{s}_{\zeta}\sigma.
\end{equation}
Being surjective, this fiber preserving map defines an embedding
\begin{equation}\label{149}
 \textsl{C}^{\infty}FB_{\eta}(J^{s}\eta,TF\eta)\hookrightarrow\textsl{C}^{\infty}FB_{\eta}(J^{r}\eta,TF\eta),\;\;
 \varXi\mapsto\varXi\circ \pi^{r}_{s}\;,
\end{equation}
that corresponds, through (\ref{147}), to the embedding (see(\ref{78}))
\begin{equation}\label{150}
 \textsl{C}^{\infty}\Gamma_{s} (T\textsl{C}^{\infty}\Gamma (\eta))\subseteq \textsl{C}^{\infty}\Gamma_{r} (T\textsl{C}^{\infty}\Gamma (\eta)).
\end{equation}
We may define accordingly the inductive limit
\begin{equation}\label{151}
\textsl{C}^{\infty}FB_{\eta}(J^{\bullet}\eta,TF\eta) = \varinjlim_{r}\textsl{C}^{\infty}FB_{\eta}(J^{r}\eta,TF\eta)
\end{equation}
and say that a vector field of the form (\ref{147}) with $\varXi\in\textsl{C}^{\infty}FB_{\eta}(J^{\bullet}\eta,TF\eta)$ is a vector field 
\textit{of finite differential order} on $\textsl{C}^{\infty}\Gamma(\eta)$ and that $\varXi$ is its \textit{total symbol}. Our aim now is to show 
that these vector fields form a Lie subalgebra of the Lie algebra of vector fields of finite order, which means to identify a structure of 
graduated Lie algebra on $\textsl{C}^{\infty}FB_{\eta}(J^{\bullet}\eta,TF\eta)$ for which the map $\varXi\longmapsto X$ given by (\ref{147})
becomes a morphism of graduated Lie algebras.\\
This structure enjoys of a certain localization property that makes the total symbols from
$\textsl{C}^{\infty}FB_{\eta}(J^{\bullet}\eta,TF\eta)$ global sections of a fine sheaf (see Warner [7]) of graduated Lie algebras over $\eta$.
The use of this property needs a more tractable notion of fiber bundle and jets of sections. In fact, in 
what follows $\eta$ will be a smooth manifold without boundary and $\pi : \eta\longrightarrow D$ an arbitrary submersion. A local section will be a 
map $\sigma : U\longrightarrow \eta$ defined on an open subset $U$ of $D$ such that $\pi\circ\sigma = \mathrm{id}_{U}$. Then its jet 
$j^{k}_{\zeta}\sigma$ of order $k$ is correctly defined for $\zeta\in U$. We define $TF_{e}\eta =: \mathrm{ker}\;T_{e}\pi\,,\;\forall e\in\eta$,
and denote $\eta_{\zeta} =: \pi^{-1} (\{\zeta\})\,,\;\zeta\in D$. Then every open subset $\xi = \mathring{\xi}\subseteq \eta$ will be considered
as a fiber bundle with the projection $\pi\rvert_{\xi}$. In this case $(J^{r}\xi)_{e} = (J^{r}\eta)_{e},\;TF_{e}\xi = TF_{e}\eta,\;\forall e\in\xi$
(see (\ref{143})). We may consider the vector spaces
\begin{equation}\label{152}
  \textsl{V}^{r}(\eta)_{e} := \textsl{C}^{\infty}((J^{r}\eta)_{e},TF_{e}\eta)\,,\;r\geqslant 0,
\end{equation}
\begin{equation}\label{153}
 \textsl{V}(\eta)_{e} :=  \varinjlim_{r}\textsl{C}^{\infty}((J^{r}\eta)_{e},TF_{e}\eta)\,,\;e\in \eta,
\end{equation}
and the associated vector bundle over $\eta$, so that every $\varXi\in\textsl{C}^{\infty}FB_{\eta}(J^{\bullet}\eta,TF\eta)$ is a smooth
section of the vector bundle $\textsl{V}(\eta)$. As, for $\xi$ open in $\eta$ and $e\in\xi$, we have 
$\textsl{V}^{r}(\xi)_{e} =\textsl{V}^{r}(\eta)_{e}\,,\;\textsl{V}(\xi)_{e} = \textsl{V}(\eta)_{e}$, the restriction 
$\varXi\rvert_{\xi}\in\textsl{C}^{\infty}FB_{\xi}(J^{\bullet}\xi,TF\xi)$ has a clear meaning.\\
As it is shown in Palais [5], for every $F : \xi\longrightarrow\eta$ fiber preserving map over $D$, there exists a well defined fiber preserving 
map (over $D$) $J^{r}(F) : J^{r}\xi\longrightarrow J^{r}\eta$ by 
\begin{equation}\label{154}
 J^{r}_{\zeta}(F)(j^{r}_{\zeta}\sigma) = j^{r}_{\zeta}(F\circ\sigma),\;\;J^{r}_{\zeta}(F) = J^{r}(F)\rvert_{J^{r}_{\zeta}\xi}\;.
\end{equation}
In the case of the fiber preserving map over $D$ (see(\ref{63})) $\tau_{\eta}\rvert_{TF\eta} :TF\eta\longrightarrow \eta$, the map 
\begin{equation}\label{155}
J^{r}_{\zeta}(\tau_{\eta}\rvert_{TF\eta}) :  J^{r}_{\zeta}(TF\eta)\longrightarrow J^{r}_{\zeta}\eta
\end{equation}
makes $J^{r}_{\zeta}(TF\eta)$ a fiber bundle over $J^{r}_{\zeta}\eta$ of fiber $(J^{r}_{\zeta}(TF\eta))_{j^{r}_{\zeta}\sigma}:= 
[J^{r}_{\zeta}(\tau_{\eta}\rvert_{TF\eta})]^{-1}(\{j^{r}_{\zeta}\sigma\})$. It is not difficult to verify that
\begin{equation}\label{156}
 (J^{r}_{\zeta}(TF\eta))_{j^{r}_{\zeta}\sigma} = J^{r}_{\zeta}(\sigma^{\ast}(TF\eta)).
\end{equation}
Indeed, we remark first that in the definition of $j^{r}_{\zeta}\sigma$ enters only the values of $\sigma$ on an arbitrary neighbourhood of 
$\zeta$ in $D$ and similarly $J^{r}_{\zeta}\eta$ depends only on the restriction of $\eta$ to such a neighbourhood. On the other hand, the 
restriction of the fiber bundle $\sigma^{\ast}(TF\eta)$ to a neighbourhood of $\zeta\in D$ depends only on the values of $\sigma$ on that 
neighbourhood. Then, in order to prove (\ref{156}) we may first restrict to a neighbourhood of $\sigma(D)$ where $\eta$ is a vector bundle and 
then to a neighbourhood of $\zeta$ where $\eta$ is a trivial vector bundle.\\
We will need in the sequel to identify the tangent space $T_{j^{r}_{\zeta}\sigma}J^{r}_{\zeta}\eta$ and in order of that we rephrase a result 
from Palais [5] (theorem 17.1, p.82) in a form more suitable for our purpose. Namely, we state it as an isomorphism of nonlinear fiber bundles 
over $D$
\begin{equation}\label{157}
 \varOmega^{r} : TFJ^{r}\eta\widetilde{\longrightarrow} J^{r}TF\eta,
\end{equation}
that is, given by diffeomorphisms depending smoothly on $\zeta\in D$
\begin{equation}\label{158}
 \varOmega^{r}_{\zeta} : T(J^{r}_{\zeta}\eta)\widetilde{\longrightarrow}J^{r}_{\zeta}(TF\eta),\;\;\zeta\in D,
\end{equation}
and these, taking into account (\ref{156}), as isomorphisms of vector bundles over $J^{r}_{\zeta}\eta$, being given by isomorphisms of vector 
spaces
\begin{equation}\label{159}
 (\varOmega^{r}_{\zeta})_{j^{r}_{\zeta}\sigma} : T_{j^{r}_{\zeta}\sigma}J^{r}_{\zeta}\eta\widetilde{\longrightarrow}
 J^{r}_{\zeta}(\sigma^{\ast}(TF\eta))
\end{equation}
depending smoothly of $j^{r}_{\zeta}\sigma\in J^{r}_{\zeta}\eta$ and $\zeta\in D$. More precisely
\begin{equation}\label{160}
 (\varOmega^{r}_{\zeta})_{j^{r}_{\zeta}\sigma}\big(\dfrac{\mathrm{d}(j^{r}_{\zeta}\sigma(t))}{\mathrm{d}t}\big |_{t=0}\big): = 
 j^{r}_{\zeta}\big (\dfrac{\mathrm{d}\sigma (t)}{\mathrm{d} t}\big |_{t=0}\big)\;,
\end{equation}
where $\sigma (0) = \sigma$ so that in the right hand side 
$\dfrac{\mathrm{d}\sigma (t)}{\mathrm{d} t}\big |_{t=0}\in\textsl{C}^{r}\Gamma(\sigma^{\ast}(TF\eta))$. To prove that 
$(\varOmega^{r}_{\zeta})_{j^{r}_{\zeta}\sigma}$, so defined by (\ref{160}), is an isomorphism depending smoothly on 
$j^{r}_{\zeta}\sigma$ and $\zeta$ we may use again a local trivial and linear restriction of $\eta$. From the definition (\ref{160}) we get (see 
(\ref{65})) 
\begin{equation}\label{161}
 (\varOmega^{r}_{\zeta})_{j^{r}_{\zeta}\sigma}\cdot T_{\sigma} j^{r}_{\zeta} = j^{r}_{\zeta}\cdot\varPhi_{\sigma}^{k+r}
\end{equation}
and
\begin{equation}\label{162}
 l(\varOmega^{r})\circ\varPhi^{k}\circ Tj^{r} = j^{r}\circ\varPhi^{k+r}
\end{equation}
on $T\textsl{C}^{k+r}\Gamma (\eta)$. Here, in the left hand side $\varPhi^{k}$ and $j^{r}$ refer to $\eta$, while in the right hand side $j^{r}$
and $\varPhi^{k+r}$ refer to $TF\eta$.\\
In the formula of the Lie bracket in the space of total symbols $\textsl{C}^{\infty}FB_{\eta}(J^{\bullet}\eta,TF\eta)$ we will use 
$\omega_{M}$, the canonical involution on $T(TM)$, for an arbitrary manifold $M$ (see, for instance, Abraham and Robbin [1]). We will also need 
the following explanatory
\begin{proposition}. Let $M$ be a smooth manifold, $p\in M,\;\;u, v\in T_{p}M$ and $A, B\in T(TM)$ be such that
\begin{equation}\label{163}
 A\in T_{u}TM,\;B\in T_{v}TM,\;\omega_{M}(A)\in T_{v}TM,\;\omega_{M}(B)\in T_{u}TM.
\end{equation}
Then the difference in $T_{u}TM$ 
\begin{equation}\label{164}
 \omega_{M}(B) - A\in T_{u}(T_{p}M),
\end{equation}
as a subspace of $T_{u}TM$, and analogously the difference in $T_{v}TM$
\begin{equation}\label{165}
 \omega_{M}(A) - B\in T_{v}(T_{p}M).
\end{equation}
Moreover, in $T_{p}M$ we have
\begin{equation}\label{166}
 [\varPsi^{T_{p}M}(v,\cdot)]^{-1}(\omega_{M}(A) - B) = -\, [\varPsi^{T_{p}M}(u,\cdot)]^{-1}(\omega_{M}(B) - A).
\end{equation}
\end{proposition}
\textbf{Proof}. We know the important property (see [1])
\begin{equation}\label{167}
 T\tau_{M}\circ \omega_{M} = \tau_{TM}.
\end{equation}
And also that
\begin{equation}\label{168}
 T_{u}(T_{p}M) = \mathrm{ker}\, T_{u}\tau_{M}.
\end{equation}
From (\ref{163}) we get $u = \tau_{TM}(A),\;v = \tau_{TM}(B)$ and using also (\ref{167}) we obtain $v = T\tau_{M}(A)$,\\
$u = T\tau_{M}(B)$. Therefore $T_{u}\tau_{M}(\omega_{M}(B) - A) = \tau_{TM}(B) - T_{u}\tau_{M}(A) = v - v = 0$ and analogously 
$T_{v}\tau_{M}(\omega_{M}(A) - B) = \tau_{TM}(A) - T_{v}\tau_{M}(B) = u - u = 0$. To verify (\ref{166}) we will use the covariance property
\begin{equation}\label{169}
 T(T\phi)\circ\omega_{M} = \omega_{N}\circ T(T\phi)
\end{equation}
for any diffeomorphism $\phi : M\longrightarrow N$, the localization property
\begin{equation}\label{170}
 \omega_{U} = \omega_{M}\rvert_{T(TU)}
\end{equation}
for $U$ open in $M$, and the formula for $M = V$ vector space
\begin{equation}\label{171}
 \omega_{V}(\varPsi^{TV}(\varPsi^{V}(x,y),\varPsi^{V}(u,v))) = \varPsi^{TV}(\varPsi^{V}(x,u),\varPsi^{V}(y,v)),
\end{equation}
$\forall\, x, y, u, v\in V$. (Here $TV$ with the structure of vector space induced by $\varPsi^{V}$ from $V\times V$). Using a 
local chart in the neighbourhood of $p\in M$ we are left to prove (\ref{166}) in a vector space $V$. Let then $u = \varPsi^{V}(p,s),\;
v = \varPsi^{V}(p,t),\;A = \varPsi^{TV}(\varPsi^{V}(p,s),\varPsi^{V}(t,x))$, so that $B = \varPsi^{TV}(\varPsi^{V}(p,t),\varPsi^{V}(s,y))$ 
(in virtue of (\ref{163}) and (\ref{171})). Then\\
$\omega_{V}(A) - B = \varPsi^{TV}(\varPsi^{V}(p,t),\varPsi^{V}(s,x) - \varPsi^{V}(s,y))$,\\
where $\varPsi^{V}(s,x) - \varPsi^{V}(s,y) = 
\varPsi^{V}(0,x-y)$, since the differences are taken in $T_{v}(TV)$, isomorphic to $TV$ through $\varPsi^{TV}(v,\cdot)$, and $TV$ is isomorphic 
to $V\times V$ through $\varPsi^{V}$. So that
\begin{equation}\label{172}
\omega_{V}(A) - B = \varPsi^{TV}(\varPsi^{V}(p,t),\varPsi^{V}(0,x-y)).
\end{equation}
Remark now that
\begin{equation}\label{173}
 \varPsi^{TV}(\varPsi^{V}(p,t),\varPsi^{V}(0,w)) = \varPsi^{T_{p}V}(\varPsi^{V}(p,t),\varPsi^{V}(p,w)),
\end{equation}
$\,\forall p, t, w\in V$. Indeed, $\varPsi^{TV}(\varPsi^{V}(p,t),\varPsi^{V}(0,w)) = \dfrac{\mathrm{d}}{\mathrm{d}h}\big (\varPsi^{V}(p,t)
+ h\varPsi^{V}(0,w)\big )\big\rvert_{h=0} =\\
= \dfrac{\mathrm{d}}{\mathrm{d}h}\varPsi^{V}(p,t+hw)\big\rvert_{h=0} = \dfrac{\mathrm{d}}{\mathrm{d}h}\big (\varPsi^{V}(p,t)+
h\varPsi^{V}(p,w)\big )\big\rvert_{h=0} = \varPsi^{T_{p}V}(\varPsi^{V}(p,t),\varPsi^{V}(p,w))$, since the second time the sum is taken in 
$T_{p}V$, which is isomorphic through $\varPsi^{V}(p,\cdot)$ to $V$. Combining (\ref{172}) with (\ref{173}) we get
\begin{equation}\label{174}
 [\varPsi^{T_{p}V}(\varPsi^{V}(p,t),\cdot)]^{-1}(\omega_{V}(A) - B) = \varPsi^{V}(p,x-y).
\end{equation}
Analogously, we get first
\begin{equation}\label{175}
 \omega_{V}(B) - A = \varPsi^{TV}(\varPsi^{V}(p,s),\varPsi^{V}(0,y-x))
\end{equation}
and then
\begin{equation}\label{176}
 [\varPsi^{T_{p}V}(\varPsi^{V}(p,s),\cdot)]^{-1}(\omega_{V}(B) - A) = \varPsi^{V}(p,y-x). 
\end{equation}
But in $T_{p}V$ we have $\varPsi^{V}(p,y-x) = -\,\varPsi^{V}(p,x-y)\,\blacksquare$\\
Now for a nonlinear fiber bundle $\pi : \eta\longrightarrow D$, we consider the following\\

\textbf{Axiom(A)}. \textit{For every $\zeta_{0}\in D$ and every smooth, locally defined around $\zeta_{0}$, {section} $\sigma_{0}$ of $\pi$
there exists a smooth global section $\sigma$ of $\pi$ that coincides with $\sigma_{0}$ on an even smaller neighbourhood of $\zeta_{0}$.}\\
Remark that in both cases of a trivial nonlinear fiber bundle and that of a vector bundle the axiom \textbf{(A)} is satisfied. Of course, it is
not satisfied by a fiber bundle without smooth global sections.\\ 
Let $f_{sr}$ denote the canonical fiber preserving map (see Palais [5]) $f_{sr} : J^{s+r}\eta\longrightarrow J^{s}(J^{r}\eta)$
\begin{equation}\label{177}
 (f_{sr})_{\zeta} : J^{s+r}_{\zeta}\eta\longrightarrow J^{s}_{\zeta}(J^{r}\eta),\;\;(f_{sr})_{\zeta}(j^{s+r}_{\zeta}\sigma) = 
 j^{s}_{\zeta}(j^{r}\sigma);
\end{equation}
it is distinct from  $f_{rs} : J^{r+s}\eta\longrightarrow J^{r}(J^{s}\eta)$ defined by $(f_{rs})_{\zeta}(j^{r+s}_{\zeta}\sigma) = 
 j^{r}_{\zeta}(j^{s}\sigma)$.\\
 
\begin{theorem}. Let us define, for $\varXi\in\textsl{C}^{\infty}FB_{\eta}(J^{r}\eta,TF\eta)$ and 
$\varTheta\in\textsl{C}^{\infty}FB_{\eta}(J^{s}\eta,TF\eta)$, 
\begin{multline}\label{178}
\big[\varXi,\varTheta\big](j^{r+s}_{\zeta}\sigma) = [\varPsi^{T_{\sigma_{\zeta}}\eta_{\zeta}}(\varXi(j^{r}_{\zeta}\sigma),\cdot)]^{-1}\cdot
[\omega_{\eta_{\zeta}}((T\varTheta\circ(\varOmega^{s})^{-1}\circ J^{s}(\varXi)\circ f_{sr})(j^{r+s}_{\zeta}\sigma)) -\\
- (T\varXi\circ(\varOmega^{r})^{-1}\circ J^{r}(\varTheta)\circ f_{rs})(j^{r+s}_{\zeta}\sigma)],
\end{multline}
the difference being taken in $T_{\varXi(j^{r}_{\zeta}\sigma)}(T\eta_{\zeta})$ with result  in 
$T_{\varXi(j^{r}_{\zeta}\sigma)}(T_{\sigma_{\zeta}}\eta_{\zeta})$ according to the previous \textbf{Proposition 6}. Then we have 
$\forall\,k, l\geqslant 0$ (see (\ref{148}))
\begin{equation}\label{179}
 \big[\varXi\circ\pi^{r+k}_{r},\varTheta\circ\pi^{s+l}_{s}\big] = \big[\varXi,\varTheta\big]\circ\pi^{r+s+k+l}_{r+s}
\end{equation}
which shows that the bracket (\ref{178}) is well defined on $\textsl{C}^{\infty}FB_{\eta}(J^{\bullet}\eta,TF\eta)$ (see(\ref{151})) and that,
for every $\xi$ open in $\eta$,
\begin{equation}\label{180}
 \big[\varXi,\varTheta\big]\rvert_{\xi} = \big[\varXi\rvert_{\xi},\varTheta\rvert_{\xi}\big].
\end{equation}
For any locally trivial fiber bundle $\eta$ and for
\begin{equation}\label{181}
X_{\sigma}(\zeta) := \varXi(j^{r}_{\zeta}\sigma),\;Y_{\sigma}(\zeta) := \varTheta (j^{s}_{\zeta}\sigma) 
\end{equation}
we have
\begin{equation}\label{182}
 \big[\varXi,\varTheta\big](j^{r+s}_{\zeta}\sigma) = \big[X,Y\big]_{\sigma}(\zeta),
\end{equation}
wherefrom we will infer that $\big[\varXi,\varTheta\big]$ satisfies the Jacobi identity for all fiber bundle. Therefore $\varXi\longmapsto X$ is a 
morphism of graduated Lie algebras from $\textsl{C}^{\infty}FB_{\eta}(J^{\bullet}\eta,TF\eta)$ to
$\textsl{C}^{\infty}\Gamma_{\bullet} (T\textsl{C}^{\infty}\Gamma (\eta))$, even a monomorphism in the case of an axiom \textbf{(A)} fiber bundle.
\end{theorem}

\textbf{Proof}. Let us verify first that (see(\ref{154}));
\begin{equation}\label{183}
(T\varTheta\circ(\varOmega^{s})^{-1}\circ J^{s}(\varXi)\circ f_{sr})(j^{r+s}_{\zeta}\sigma)\in T_{\varTheta(j^{s}_{\zeta}\sigma)}(T\eta_{\zeta}) 
\end{equation}
and that
\begin{equation}\label{184}
\omega_{\eta_{\zeta}}((T\varTheta\circ(\varOmega^{s})^{-1}\circ J^{s}(\varXi)\circ f_{sr})(j^{r+s}_{\zeta}\sigma))\in 
T_{\varXi(j^{r}_{\zeta}\sigma)}(T\eta_{\zeta}).
\end{equation}
For the first we see that $J^{s}(\varXi)(j^{s}_{\zeta}(j^{r}\sigma)) = j^{s}_{\zeta}(\varXi(j^{r}\sigma))$ and that\\ 
$[(\varOmega^{s}_{\zeta})_{j^{s}_{\zeta}\sigma}]^{-1}j^{s}_{\zeta}(\varXi(j^{r}\sigma))\in T_{j^{s}_\zeta\sigma} J^{s}_{\zeta}\eta$. Finally
$T\varTheta (T_{j^{s}_\zeta\sigma} J^{s}_{\zeta}\eta)\subseteq T_{\varTheta(j^{s}_{\zeta}\sigma)}(T\eta_{\zeta})$, which proves (\ref{183}).
To prove (\ref{184}) we apply on the left hand side $\tau_{T\eta_{\zeta}}$ aiming to obtain $\varXi(j^{r}_{\zeta}\sigma)$. We use then (\ref{167})
and $\omega_{M}\circ\omega_{M} = \mathrm{id}_{T(TM)}$, so that $\tau_{T\eta_{\zeta}}\circ\omega_{\eta_{\zeta}} = T\tau_{\eta_{\zeta}}$. We note 
the relation
\begin{equation}\label{185}
 \tau_{\eta_{\zeta}}\circ\varTheta = \pi^{s}_{0}
\end{equation}
(where $\pi^{s}_{0} = \pi^{s}_{0}(\eta) : J^{s}\eta\longrightarrow\eta$) and also that 
\begin{equation}\label{186}
 T\pi^{s}_{0}(\eta) = \pi^{s}_{0}(TF\eta)\circ\varOmega^{s},
\end{equation}
if $\pi^{s}_{0}(TF\eta) : J^{s}(TF\eta)\longrightarrow TF\eta$. In this way we obtain that 
$T\tau_{\eta_{\zeta}}\circ T\varTheta\circ(\varOmega^{s})^{-1} = \pi^{s}_{0}(TF\eta)$, wherefrom the desired result. Analogously, simply by 
interchanging $r$ with $s$ and $\varXi$ with $\varTheta$, we get
\begin{equation}\label{187}
(T\varXi\circ(\varOmega^{r})^{-1}\circ J^{r}(\varTheta)\circ f_{rs})(j^{r+s}_{\zeta}\sigma)\in T_{\varXi(j^{r}_{\zeta}\sigma)}(T\eta_{\zeta})
\end{equation}
and
\begin{equation}\label{188}
 \omega_{\eta_{\zeta}}((T\varXi\circ(\varOmega^{r})^{-1}\circ J^{r}(\varTheta)\circ f_{rs})(j^{r+s}_{\zeta}\sigma))\in 
 T_{\varTheta(j^{s}_{\zeta}\sigma)}(T\eta_{\zeta}),
\end{equation}
which, together with Proposition 6, gives the meaning to the definition (\ref{178}) (and shows moreover that $\big[\varXi,\varTheta\big]
= - \big[\varTheta,\varXi\big]$ in this relation). Thus in order to verify (\ref{179}) it is enough to show that
\begin{equation}\label{189}
 \big[\varXi\circ\pi^{r+k}_{r},\varTheta\big] = \big[\varXi,\varTheta\big]\circ\pi^{r+s+k}_{r+s}.
\end{equation}
But this is the easy consequence of the following commutation relations
\begin{equation}\label{190}
 f_{sr}\circ\pi^{s+r+k}_{s+r} = J^{s}(\pi_{r}^{r+k})\circ f_{s,r+k},
\end{equation}
\begin{equation}\label{191}
 \varOmega^{r}\circ T\pi^{r+k}_{r}(\eta) = \pi^{r+k}_{r}(TF\eta)\circ\varOmega^{r+k}
\end{equation}
(which extends (\ref{186})) and
\begin{equation}\label{192}
 \pi^{r+k}_{r}(TF\eta)\circ J^{r+k}(\varTheta)\circ f_{r+k,s} = J^{r}(\varTheta)\circ f_{rs}\circ\pi^{s+r+k}_{s+r}(\eta).
\end{equation}
In proving (\ref{182}) we have to consider a vector bundle neighbourhood $E$ of $\sigma$. For the fiber bundles over $\eta$:
\begin{equation}\label{193}
 \pi^{r+s}_{0}(\eta) : J^{r+s}\eta\longrightarrow\eta\,,
\end{equation}
\begin{equation}\label{194}
 \pi^{r}_{0}(\eta)\circ\pi^{s}_{0}(J^{r}\eta) : J^{s}(J^{r}\eta)\longrightarrow\eta\,,
\end{equation}
\begin{equation}\label{195}
 \tau_{\eta}\circ\pi^{s}_{0}(TF\eta) : J^{s}(TF\eta)\longrightarrow\eta\,,
\end{equation}
\begin{equation}\label{196}
 \pi^{s}_{0}(\eta)\circ\tau_{J^{s}\eta} : TF(J^{s}\eta)\longrightarrow\eta\,,
\end{equation}
\begin{equation}\label{197}
 \tau_{\eta}\circ T\tau_{\eta} : TF(TF\eta)\longrightarrow\eta\,,
\end{equation}
where
\begin{equation}\label{198}
 [TF(TF\eta)]_{\zeta} = T((TF\eta)_{\zeta}) = T(T(\eta_{\zeta}))\,,
\end{equation}
the mappings
\begin{equation}\label{199}
 f_{sr} : J^{r+s}\eta\longrightarrow J^{s}(J^{r}\eta)\,,
\end{equation}
\begin{equation}\label{200}
 J^{s}(\varXi) : J^{s}(J^{r}\eta)\longrightarrow J^{s}(TF\eta)\,,
\end{equation}
\begin{equation}\label{201}
 (\varOmega^{s})^{-1} : J^{s}(TF\eta)\longrightarrow TF(J^{s}\eta)\,,
\end{equation}
\begin{equation}\label{202}
 T\varTheta : TF(J^{s}\eta)\longrightarrow TF(TF\eta)\,,
\end{equation}
\begin{equation}\label{203}
 \omega_{\eta} : TF(TF\eta)\longrightarrow TF(TF\eta)
\end{equation}
preserve fibers over $\eta$. We use here the fact that $\varXi$ and $\varTheta$ preserve fibers over $\eta$ and that
\begin{equation}\label{204}
 \tau_{\eta}\circ T\tau_{\eta} = \tau_{\eta}\circ\tau_{T\eta}.
\end{equation}
Analogously, we consider also $\varOmega^{r},\,J^{r}(\varTheta),\,T\varXi,\,f_{rs}$ as fiber preserving maps over $\eta$. If 
we use the generic notation $\pi_{F,\eta} : F(\eta)\longrightarrow\eta$ for any of the fiber bundles (\ref{193})-(\ref{197}),
we see that, for $\xi$ open in $\eta$, $F(\xi) = \pi_{F,\eta}^{-1}(\xi),\;\pi_{F,\xi} = \pi_{F,\eta}\rvert_{F(\xi)}$. Taking into account 
that the mappings (\ref{199})-(\ref{203}) preserve fibers over $\eta$, their restrictions to the fibers over $\xi$ give the respective to $\xi$ 
mappings. This proves (\ref{180}). On the other hand, this shows that it is enough to prove (\ref{182}) for $\eta = E$ vector bundle over $D$.\\
For $X,\,Y$ given by (\ref{181}) and $\tilde{X},\,\tilde{Y}$ defined according to (\ref{87}) we consider 
$\tilde{\varXi},\,\tilde{\varTheta}$ such that
\begin{equation}\label{205}
 \tilde{\varXi}(j^{r}_{\zeta}\sigma) = \tilde{X}_{\sigma}(\zeta),\;\tilde{\varTheta} (j^{s}_{\zeta}\sigma) = \tilde{Y}_{\sigma}(\zeta).
\end{equation}
Thus $\tilde{\varXi} : J^{r}E\longrightarrow E$ is defined by $\varXi$:
\begin{equation}\label{206}
 \varXi(j^{r}_{\zeta}\sigma) = \varPsi^{E_{\zeta}}(\sigma_{\zeta},\tilde{\varXi}(j^{r}_{\zeta}\sigma))\,,
\end{equation}
being fiber preserving over $D$:
\begin{equation}\label{207}
 \tilde{\varXi}_{\zeta} =: \tilde{\varXi}\rvert_{J^{r}_{\zeta}E}\,,\; \tilde{\varXi}_{\zeta} : J^{r}_{\zeta}E\longrightarrow E_{\zeta},\;
 \zeta\in D.
\end{equation}
Then $\tilde{X} :\textsl{C}^{k+r}\Gamma (E)\longrightarrow\textsl{C}^{k}\Gamma (E)),\;
\tilde{Y} :\textsl{C}^{k+s}\Gamma (E)\longrightarrow\textsl{C}^{k}\Gamma (E)),\;\;\forall k\geqslant 0$ and\\
$<\tilde{Y}^{\prime}(\sigma);\tilde{X}(\sigma)>_{\zeta} = 
\dfrac{\mathrm{d}}{\mathrm{d}t}\tilde{Y}(\sigma + t\tilde{X} (\sigma))_{\zeta}\rvert_{t=0} = 
\dfrac{\mathrm{d}}{\mathrm{d}t}\tilde{\varTheta}_{\zeta}(j^{s}_{\zeta}\sigma + t j^{s}_{\zeta}(\tilde{\varXi}(j^{r}\sigma)))\rvert_{t=0} =\\
= <\tilde{\varTheta}^{\prime}_{\zeta}(j^{s}_{\zeta}\sigma);j^{s}_{\zeta}(\tilde{\varXi}(j^{r}\sigma))>$.\\
Here $\tilde{\varTheta}_{\zeta} : J^{s}_{\zeta}E\longrightarrow E_{\zeta}$ acts between vector spaces and then 
$\tilde{\varTheta}^{\prime}_{\zeta}(j^{s}_{\zeta}\sigma)\in \mathrm{Hom}(J^{s}_{\zeta}E, E_{\zeta})$. Therefore, in virtue of the definition
(\ref{88}), (\ref{90}) we have $\widetilde{\big[X,Y\big]} = \big[\tilde{X},\tilde{Y}\big]$ and
\begin{equation}\label{208}
 \big[\tilde{X},\tilde{Y}\big](\sigma)_{\zeta} = <\tilde{\varTheta}^{\prime}_{\zeta}(j^{s}_{\zeta}\sigma);j^{s}_{\zeta}(\tilde{\varXi}(j^{r}\sigma))> -\\
 <\tilde{\varXi}^{\prime}_{\zeta}(j^{r}_{\zeta}\sigma);j^{r}_{\zeta}(\tilde{\varTheta}(j^{s}\sigma))>.
\end{equation}
In the right hand side of (\ref{178}) we have $\varXi(j^{r}_{\zeta}\sigma)\in T_{\sigma_{\zeta}}E_{\zeta}\,,\,\varXi(j^{r}_{\zeta}\sigma) =\\
= \varPsi^{E_{\zeta}}(\sigma_{\zeta},\tilde{\varXi}(j^{r}_{\zeta}\sigma)) = 
\dfrac{\mathrm{d}}{\mathrm{d}t}(\sigma + t \tilde{\varXi}(j^{r}\sigma))_{\zeta}\rvert_{t=0}$. Therefore 
$[(\varOmega^{s}_{\zeta})_{j^{s}_{\zeta}}]^{-1}(j^{s}_{\zeta}(\varXi(j^{r}\sigma))) =$ \\

$= \dfrac{\mathrm{d}}{\mathrm{d}t}(j^{s}_{\zeta}\sigma + t j^{s}_{\zeta}(\tilde{\varXi}(j^{r}\sigma)))\rvert_{t=0}$ and 
$T_{j^{s}_{\zeta}\sigma}\varTheta\cdot [(\varOmega^{s}_{\zeta})_{j^{s}_{\zeta}}]^{-1}(j^{s}_{\zeta}(\varXi(j^{r}\sigma))) =$\\

$= \dfrac{\mathrm{d}}{\mathrm{d}t}\varTheta (j^{s}_{\zeta}\sigma + t j^{s}_{\zeta}(\tilde{\varXi}(j^{r}\sigma)))\rvert_{t=0} = 
\dfrac{\mathrm{d}}{\mathrm{d}t}\varPsi^{E_{\zeta}}(\sigma_{\zeta} + 
t \tilde{\varXi}(j^{r}_{\zeta}\sigma),\tilde{\varTheta}(j^{s}_{\zeta}\sigma + t j^{s}_{\zeta}(\tilde{\varXi}(j^{r}\sigma))))\rvert_{t=0} =$\\

$= \dfrac{\mathrm{d}}{\mathrm{d}t}\big (\varPsi^{E_{\zeta}}(\sigma_{\zeta},\tilde{\varTheta}(j^{s}_{\zeta}\sigma)) + 
t \varPsi^{E_{\zeta}}(\tilde{\varXi}(j^{r}_{\zeta}\sigma),
<\tilde{\varTheta}^{\prime}_{\zeta}(j^{s}_{\zeta}\sigma);j^{s}_{\zeta}(\tilde{\varXi}(j^{r}\sigma))>)\big )\rvert_{t=0} =$\\

$=\varPsi^{TE_{\zeta}}(\varPsi^{E_{\zeta}}(\sigma_{\zeta},\tilde{\varTheta}(j^{s}_{\zeta}\sigma)), 
\varPsi^{E_{\zeta}}(\tilde{\varXi}(j^{r}_{\zeta}\sigma),
<\tilde{\varTheta}^{\prime}_{\zeta}(j^{s}_{\zeta}\sigma);j^{s}_{\zeta}(\tilde{\varXi}(j^{r}\sigma))>))$\\

in virtue of the linearity of $\varPsi^{E_{\zeta}} : E_{\zeta}\times E_{\zeta}\longrightarrow TE_{\zeta}$. Therefore, with the definitions 
(\ref{154}) and (\ref{177}) we obtain\\

$\omega_{\eta_{\zeta}}((T\varTheta\circ(\varOmega^{s})^{-1}\circ J^{s}(\varXi)\circ f_{sr})(j^{r+s}_{\zeta}\sigma))
- (T\varXi\circ(\varOmega^{r})^{-1}\circ J^{r}(\varTheta)\circ f_{rs})(j^{r+s}_{\zeta}\sigma) =$\\

$= \varPsi^{TE_{\zeta}}(\varPsi^{E_{\zeta}}(\sigma_{\zeta},\tilde{\varXi}(j^{r}_{\zeta}\sigma)),$\\

$,\varPsi^{E_{\zeta}}(0,
<\tilde{\varTheta}^{\prime}_{\zeta}(j^{s}_{\zeta}\sigma);j^{s}_{\zeta}(\tilde{\varXi}(j^{r}\sigma))> - 
<\tilde{\varXi}^{\prime}_{\zeta}(j^{r}_{\zeta}\sigma);j^{r}_{\zeta}(\tilde{\varTheta}(j^{s}\sigma))>)) =$\\

$= \varPsi^{T_{\sigma_{\zeta}}E_{\zeta}}(\varPsi^{E_{\zeta}}(\sigma_{\zeta},\tilde{\varXi}(j^{r}_{\zeta}\sigma)),$\\

$,\varPsi^{E_{\zeta}}(\sigma_{\zeta},
<\tilde{\varTheta}^{\prime}_{\zeta}(j^{s}_{\zeta}\sigma);j^{s}_{\zeta}(\tilde{\varXi}(j^{r}\sigma))> - 
<\tilde{\varXi}^{\prime}_{\zeta}(j^{r}_{\zeta}\sigma);j^{r}_{\zeta}(\tilde{\varTheta}(j^{s}\sigma))>)) =$\\

$= \varPsi^{T_{\sigma_{\zeta}}E_{\zeta}}(\varXi(j^{r}_{\zeta}\sigma),
\varPsi^{E_{\zeta}}(\sigma_{\zeta},[\tilde{X},\tilde{Y}](\sigma)_{\zeta})) = \varPsi^{T_{\sigma_{\zeta}}E_{\zeta}}(\varXi(j^{r}_{\zeta}\sigma),
\big[X,Y\big]_{\sigma}(\zeta))$\\

according to (\ref{208}). We have used also (\ref{171}) and (\ref{173}). Finally, in virtue of the axiom \textbf{(A)} the map 
$\varXi\longmapsto X$ is injective. The injectivity, in turn, ensures that the bracket (\ref{178}) satisfies the Jacobi identity.
This ends the proof of the Theorem $\blacksquare$

\textbf{Remark}. The formula (\ref{178}) of the Lie bracket in $\textsl{C}^{\infty}FB_{\eta}(J^{\bullet}\eta,TF\eta)$ may be not of great value 
as such; it only shows that the vector fields of finite differential order form a Lie subalgebra in 
$\textsl{C}^{\infty}\Gamma_{\bullet} (T\textsl{C}^{\infty}\Gamma (\eta))$ and that the Lie bracket of their total symbols enjoys of the 
localization property (\ref{180}). It results that the set of germs of sections of the vector bundle $\textsl{V}(\eta)$ (see (\ref{153})) bears a 
structure of sheaf of graduated Lie algebras. (The stalks are Lie algebras and not the fibers of $\textsl{V}(\eta)$, like in the case of the 
vector bundle $T\eta$ over $\eta$). The vector bundle $\textsl{V}(\eta)$ depends only on the completely integrable vector subbundle $TF\eta$ of 
$T\eta$ since $J^{r}\eta$ is determined by the local foliation given by $TF\eta$. It does not depend on the projection $\pi$ $\blacksquare$\\

\section{An additive formula for the curvature tensors of two supplementary subbundles}

The previous digression is used, in particular, to argue the next result which refers to the case $\eta = D\times M$ as a fiber bundle over $D$, 
when $\textsl{C}^{k}\Gamma ( D\times M)\widetilde{\longrightarrow}\textsl{C}^{k}(D,M)$, through the correspondence $\sigma\mapsto \beta,\,\;
\sigma_{\zeta} = (\zeta, \beta (\zeta)),\;\;\zeta\in D$ (see(\ref{91})). In this case we use the notations
\begin{equation}\label{209}
 J^{r}(D,M) := J^{r}\eta\,,\;(\zeta,j^{r}_{\zeta}\beta) = j^{r}_{\zeta}\sigma\,,\;
 J^{r}(D,M)_{(\zeta,p)} = \{(\zeta,j^{r}_{\zeta}\beta)\,\rvert\,\beta(\zeta)=p\},
\end{equation}
\begin{equation}\label{210}
\textsl{C}^{\infty}FB_{M}(J^{r}(D,M),TM) := \textsl{C}^{\infty}FB_{\eta}(J^{r}\eta,TF\eta)
\end{equation}
since in this case $J^{r}(D,M)$ becomes a fiber bundle over $M$ through $(\zeta,j^{r}_{\zeta}\beta) \mapsto \beta (\zeta)$ and 
$TF_{(\zeta,p)}(D\times M) \tilde{=} \{0_{T_{\zeta}D}\}\times T_{p} M$ makes $TF(D\times M)$ a vector bundle over $M$. Thus for $\varXi\in
\textsl{C}^{\infty}FB_{M}(J^{r}(D,M),TM) $ we have $\varXi(\zeta, j^{r}_{\zeta}\beta)\in T_{\beta(\zeta)}M\,,\;\forall\zeta\in D\,,\;
\forall\beta\in\textsl{C}^{r}(D,M)$. Also, in this case, we denote
\begin{equation}\label{211}
 \textsl{C}^{\infty}FB_{M}(J^{\bullet}(D,M),TM) := \textsl{C}^{\infty}FB_{\eta}(J^{\bullet}\eta,TF\eta)
\end{equation}
(see (\ref{151})) and
\begin{equation}\label{212}
 \textsl{V}^{r}(D,M)_{(\zeta,p)} = \textsl{C}^{\infty}(J^{r}(D,M)_{(\zeta,p)},T_{p}M),
\end{equation}
\begin{equation}\label{213}
  \textsl{V}(D,M)_{(\zeta,p)} =  \varinjlim_{r}\textsl{V}^{r}(D,M)_{(\zeta,p)}.
\end{equation}
For $A\in\textsl{C}^{\infty}\Gamma(T^{\ast}M\otimes TM)$ we denote
\begin{equation}\label{214}
 (l(A))_{\beta}(\zeta) := A_{\beta(\zeta)}\in T^{\ast}_{\beta(\zeta)}\otimes T_{\beta(\zeta)}M
\end{equation}
and if $X\in\textsl{C}^{\infty}\Gamma (TD)$ we put 
\begin{equation}\label{215}
(l(A)\cdot r(X))_{\beta}(\zeta) = A_{\beta(\zeta)}\cdot T_{\zeta}\beta\cdot X_{\zeta}\,,\;\zeta\in D\,,\;\beta\in\textsl{C}^{\infty}(D,M)
\end{equation}
(compare with (\ref{106})). In this way, $l(A)\cdot r(X)$ is a vector field of differential order $1$ on $\textsl{C}^{\infty}(D,M)$.

 \begin{theorem}. For $X,\,Y\in\textsl{C}^{\infty}\Gamma (TD)$  and $A\in\textsl{C}^{\infty}\Gamma(T^{\ast}M\otimes TM)$ the following equality
holds
\begin{equation}\label{216}
 \big[l(A)\cdot r(X),r(Y)\big] = - l(A)\cdot r(\big[X,Y\big])
\end{equation}
as for vector fields on $\textsl{C}^{\infty}(D,M)$.
 \end{theorem}

\textbf{Proof}. Remark that for $A_{p} = I_{T_{p}M}\,,\;\forall p\in M$, (\ref{216}) becomes (\ref{109}). The localization property of the Lie 
bracket, expressed by the relations (\ref{182}) and 
(\ref{180}), allow to reduce proving the identity on a product $G\times U$, where $G$ and $U$ are domains of chart in $D$ and $M$ respectively,
that is, to show that $[l(A)\cdot r(X),r(Y)]_{\beta}(\zeta) = - A_{\beta(\zeta)}\cdot T_{\zeta}\beta\cdot [X,Y]_{\zeta}$, for $\zeta\in G$ and 
$\beta\in\textsl{C}^{\infty}(G,U)$. Taking into account the 
intrinsic geometric meaning of the operations, we may consider only the case when $G=\mathring{G}\subseteq L, \,U = \mathring{U}\subseteq V$, 
for $L$ and $V$ vector spaces. Then $A : U\longrightarrow V^{\ast}\otimes V\,,\; X\,,\,Y : G\longrightarrow L$
and $(l(A)\cdot r(X))_{\beta}(\zeta) = A(\beta(\zeta))\cdot \beta^{\prime}(\zeta)\cdot X(\zeta)$. In this case\\

$<(l(A)\cdot r(X))_{\beta}^{\prime}(\beta);r(Y)_{\beta}>(\zeta) = <A^{\prime}(\beta(\zeta));r(Y)_{\beta}(\zeta)>\cdot \beta^{\prime}(\zeta)
\cdot X(\zeta) + $\\

$+ A(\beta(\zeta))\cdot <(r(Y)_{\beta})^{\prime}_{\zeta}(\zeta);X(\zeta)> = <A^{\prime}(\beta(\zeta));\beta^{\prime}(\zeta)Y(\zeta)>\cdot
\beta^{\prime}(\zeta)\cdot X(\zeta) + $\\

$+ A(\beta(\zeta))\cdot\beta^{\prime\prime}(\zeta)(X(\zeta),Y(\zeta)) + 
A(\beta(\zeta))\cdot <\beta^{\prime}(\zeta);Y^{\prime}(\zeta)X(\zeta)>$\\

and

$<r(Y)^{\prime}_{\beta}(\beta);(l(A)\cdot r(X))_{\beta}>(\zeta) = <((l(A)\cdot r(X))_{\beta})^{\prime}_{\zeta}(\zeta);Y(\zeta)> = $\\

$= <A^{\prime}(\beta(\zeta));\beta^{\prime}(\zeta)Y(\zeta)>\cdot\beta^{\prime}(\zeta)\cdot X(\zeta) + 
A(\beta(\zeta))\cdot \beta^{\prime\prime}(\zeta)(Y(\zeta),X(\zeta)) + $\\

$+ A(\beta(\zeta))\cdot <\beta^{\prime}(\zeta);X^{\prime}(\zeta)Y(\zeta)>$.

Therefore \\

$\big[l(A)\cdot r(X),r(Y)\big]_{\beta}(\zeta) = <r(Y)^{\prime}_{\beta}(\beta);(l(A)\cdot r(X))_{\beta}>(\zeta) -$\\

$- <(l(A)\cdot r(X))_{\beta}^{\prime}(\beta);r(Y)_{\beta}>(\zeta) = $\\

$= A(\beta(\zeta))\cdot<\beta^{\prime}(\zeta);X^{\prime}(\zeta)Y(\zeta) - Y^{\prime}(\zeta)X(\zeta)> = - 
(l(A)\cdot r(\big[X,Y\big]))_{\beta}(\zeta)\;\;
\blacksquare$

\begin{corollary}. In the same conditions on $X, Y$ and $A$ we have
\begin{multline}\label{217}
 \big[l(A)\cdot r(X),l(A)\cdot r(Y)\big] + l(A)\cdot r(\big[X,Y\big]) = \big[l(I-A)\cdot r(X),l(I-A)\cdot r(Y)\big] +\\
 + l(I-A)\cdot r(\big[X,Y\big]).
\end{multline}
\end{corollary}

\textbf{Proof}. From (\ref{216}) we get also
\begin{equation}\label{218}
  \big[r(X),l(A)\cdot r(Y)\big] = - l(A)\cdot r(\big[X,Y\big]).
\end{equation}
As $ l(I-A)\cdot r(X) = r(X) - l(A)\cdot r(X)$, (\ref{217}) reduces to $l(A)\cdot r(\big[X,Y\big]) = \\
= - \big[l(A)\cdot r(X),r(Y)\big] - \big[r(X),l(A)\cdot r(Y)\big] - l(A)\cdot r(\big[X,Y\big])\;\blacksquare$\\

The interior direct sum of two subspaces was denoted by
\begin{equation}\label{219}
 H \dot{+} K = V\,.
\end{equation}
In this case
\begin{equation}\label{220}
 P^{K}_{H} : V\longrightarrow H\subseteq V
\end{equation}
will stand for the canonical projection. Of course, in $V^{\ast}\otimes V$ we have
\begin{equation}\label{221}
  P^{K}_{H} + P^{H}_{K} = I_{V}.
\end{equation}
Also 
\begin{equation}\label{222}
 Q^{K}_{H} : V/K\widetilde{\longrightarrow} H
\end{equation}
will denote the canonical isomorphism. In the case of two supplementary subbundles $H$ and $K$ of $TM$
\begin{equation}\label{223}
 H_{m} \dot{+} K_{m} = T_{m}M,\;\;\forall m\in M\,,
\end{equation}
$P^{K}_{H}$ will denote the
section of $T^{\ast}M\otimes TM$
\begin{equation}\label{224}
 (P^{K}_{H})_{m} = P^{K_{m}}_{H_{m}}\,,\; m\in M,
\end{equation}
and $Q^{K}_{H}$ the section of $\mathrm{Hom}(TM/K,H)$
\begin{equation}\label{225}
 (Q^{K}_{H})_{m} = Q^{K_{m}}_{H_{m}}\,,\; m\in M.
\end{equation}
Of course (see (\ref{21}))
\begin{equation}\label{226}
 (Q^{K}_{H})_{m}\cdot P^{K}_{m} = (P^{K}_{H})_{m}\,,\;\forall m\in M.
\end{equation}
From Theorem 2 and Corollary 1 we obtain the important

\begin{theorem}.
  If $H$ and $K$ are two supplementary vector subbundles of $TM$
\begin{multline}\label{227}
 Q^{H_{\beta(\zeta)}}_{K_{\beta(\zeta)}}C^{H}_{\beta(\zeta)}(P^{K_{\beta(\zeta)}}_{H_{\beta(\zeta)}}\; T_{\zeta}\beta\; X_{\zeta}\wedge
 P^{K_{\beta(\zeta)}}_{H_{\beta(\zeta)}}\; T_{\zeta}\beta\; Y_{\zeta}) + \\
 + Q^{K_{\beta(\zeta)}}_{H_{\beta(\zeta)}}C^{K}_{\beta(\zeta)}(P^{H_{\beta(\zeta)}}_{K_{\beta(\zeta)}}\; T_{\zeta}\beta\; X_{\zeta}\wedge
 P^{H_{\beta(\zeta)}}_{K_{\beta(\zeta)}}\; T_{\zeta}\beta\; Y_{\zeta}) =\\
 = \big[l(P^{H}_{K})\cdot r(X),l(P^{H}_{K})\cdot r(Y)\big]_{\beta}(\zeta) + (l(P^{H}_{K})\cdot r(\big[X,Y\big]))_{\beta}(\zeta)\,,
\end{multline}
for all $X, Y\in\textsl{C}^{\infty}\Gamma (TD)\,,\;\beta\in\textsl{C}^{\infty}(D,M)\,,\;\zeta \in D$. 
\end{theorem}

\textbf{Proof}. In the Corollary 1 above we take $A = P^{K}_{H}$ and apply on both sides $P^{K_{\beta(\zeta)}}_{H_{\beta(\zeta)}}$ thus 
obtaining\\
$l(P^{K}_{H})\cdot \big[l(P^{K}_{H})\cdot r(X),l(P^{K}_{H})\cdot r(Y)\big] + l(P^{K}_{H})\cdot r(\big[X,Y\big]) =
l(P^{K}_{H})\cdot \big[l(P^{H}_{K})\cdot r(X),l(P^{H}_{K})\cdot r(Y)\big]$.\\
In this way, using Theorem 2, (\ref{226}), (\ref{221}) and Corollary 1 we get

$Q^{K_{\beta(\zeta)}}_{H_{\beta(\zeta)}}C^{K}_{\beta(\zeta)}(P^{H_{\beta(\zeta)}}_{K_{\beta(\zeta)}}\; 
T_{\zeta}\beta\; X_{\zeta}\wedge
P^{H_{\beta(\zeta)}}_{K_{\beta(\zeta)}}\; T_{\zeta}\beta\; Y_{\zeta}) = P^{K_{\beta(\zeta)}}_{H_{\beta(\zeta)}}
(\big[l(P^{K}_{H})\cdot r(X),l(P^{K}_{H})\cdot r(Y)\big]_{\beta}(\zeta) +$\\

$+ (l(P^{K}_{H})\cdot r(\big[X,Y\big]))_{\beta}(\zeta)) = 
P^{K_{\beta(\zeta)}}_{H_{\beta(\zeta)}}(\big[l(P^{H}_{K})\cdot r(X),l(P^{H}_{K})\cdot r(Y)\big]_{\beta}(\zeta) +$\\

$+ (l(P^{H}_{K})\cdot r(\big[X,Y\big]))_{\beta}(\zeta))$.

Transposing here $H$ and $K$, adding these relations and using again (\ref{221}) we obtain (\ref{227}). The symmetry in the pair $(H,K)$, plain 
in the left hand side of (\ref{227}), comes in the right hand side from Corollary 1 $\blacksquare$

\textbf{Remark}. The formula (\ref{227}) may be red as an expression for the Lie bracket\\ 
$\big[l(P^{H}_{K})\cdot r(X),l(P^{H}_{K})\cdot r(Y)\big]$; like (\ref{216}) for $\big[l(A)\cdot r(X),r(Y)\big]$,
it shows that this vector field is of differential order at most $1$, when the
general formula (\ref{178}), of no use here, would give a differential order at most $2$
$\blacksquare$

\section{The parallel transport in a supplementary\\
 vector subbundle along a tangent path to\\
 the vector subbundle under study}

Here $P$ stands for the projection
\begin{equation}\label{228}
 P = P^{H}_{K}
\end{equation}
corresponding to a direct sum decomposition (\ref{223}). Conversely, starting from a smooth section $P$ of $T^{\ast}M\otimes TM$ with 
$P_{m}^{2} = P_{m}\,,\;\forall\,m$, we get the smooth subbundles $H$ and $K$ and the decomposition. Remark that the continuity in $m$ ensures
that $P_{m}$ has a constant rank, since the rank of both $P_{m}$ and $I-P_{m}$ may only increase in a neighbourhood of $m$.
It is easy to show that \textit{for a smooth section} $P$ \textit{in projections of} $T_{m} M\,,\;\forall\,m\in M\,,$ \textit{its lift to  
$T(TM)$ as
\begin{equation}\label{229}
 \omega_{M}\circ TP\circ\omega_{M}
\end{equation}
is again a linear projection in $T_{X} TM$ and moreover
\begin{equation}\label{230}
 \omega_{M}\circ TP\circ\omega_{M} + \omega_{M}\circ T(I-P)\circ\omega_{M} = I_{T_{X} TM}\,, 
\end{equation}
$\forall X\in TM$}. Here $P$ is considered as a smooth map 
\begin{equation}\label{231}
 P : TM\longrightarrow TM\,,\;P(X) = P_{\tau_{M}(X)}\cdot X\,,
\end{equation}
so that $TP : T(TM)\longrightarrow T(TM)$. In what follows $H$ and $K$ (see (\ref{228})) will denote the total space of the respective vector 
bundles, hence the respective submanifolds of $TM$. 

\begin{theorem}. Let $H, K, P$ be as in (\ref{228}) and smooth on $M$. Then the restriction of $T\tau_{M} + \tau_{TM}$ from $T(TM)$
\begin{equation}\label{232}
 T\tau_{M} + \tau_{TM} : TK\cap\omega_{M}(TH)\longrightarrow TM
\end{equation}
is a bijection of inverse $R_{P}$ that verifies more precisely
\begin{equation}\label{233}
 \tau_{TM}\circ R_{P} = P,\;\;T\tau_{M}\circ R_{P} = I-P.
\end{equation}
Then, for all $P$
\begin{equation}\label{234}
 \omega_{M}\circ R_{P} = R_{(I-P)}.
\end{equation}
It results that $T_{k}K\cap\omega_{M}(TH)$ is a vector subspace in $T_{k}K$ such that
\begin{equation}\label{235}
 R_{P}(k+ (\cdot)) : H_{\tau_{M}(k)}\widetilde{\longrightarrow} T_{k}K\cap\omega_{M}(TH)
\end{equation}
is a linear isomorphism and section for
\begin{equation}\label{236}
 T_{k}\tau_{M} : N_{k}\longrightarrow H_{\tau_{M}(k)}\,,
\end{equation}
where 
\begin{equation}\label{237}
 N_{k} =: T_{k}K\cap (T_{k}\tau_{M})^{-1}(H_{\tau_{M}(k)}).
\end{equation}
Finally, $N_{k}$ is invariant for $\omega_{M}\circ T(I-P)\circ\omega_{M}$, and
\begin{equation}\label{238}
\omega_{M}\circ T(I-P)\circ\omega_{M}\rvert_{N_{k}} = R_{P}(k+ (\cdot))\cdot T_{k}\tau_{M}\rvert_{N_{k}}, 
\end{equation}
$\forall k\in K$, that is $\omega_{M}\circ T(I-P)\circ\omega_{M}$ projects $N_{k}$ on $T_{k}K\cap\omega_{M}(TH)$. It results that 
$\omega_{M}\circ TP\circ\omega_{M}$ projects $N_{k}$ on $T_{k}K_{\tau_{M}(k)}$.
\end{theorem}

\textbf{Proof}. We start with $W, Z\in T(TM)$ such that
\begin{equation}\label{239}
 TP\cdot W = \omega_{M}(T(I-P)\cdot Z)
\end{equation}
and show that, if $Y\in T(TM)$ satisfies
\begin{equation}\label{240}
 \omega_{M}(Y) = Y,
\end{equation}
\begin{equation}\label{241}
 P(\tau_{TM}(Y)) =  P(\tau_{TM}(W))\,,\;\;(I- P)(\tau_{TM}(Y)) =  (I-P)(\tau_{TM}(Z)) 
\end{equation}
then (see (\ref{239}))
\begin{multline}\label{242}
  TP\cdot W = \omega_{M}(T(I-P)\cdot Z) =\\
  = (TP\circ\omega_{M}\circ T(I-P))(Y) = (\omega_{M}\circ T(I-P) \circ\omega_{M}\circ TP)(Y).
\end{multline}
Remark that the condition (\ref{240}) is equivalent to
\begin{equation}\label{243}
 \tau_{TM}(Y) = T\tau_{M}(Y)
\end{equation}
and that for each $X\in TM$ there are (a diffeomorphic to the model vector space of $M$ set of)
$Y$-s such that
\begin{equation}\label{244}
\tau_{TM}(Y) = T\tau_{M}(Y) = X. 
\end{equation}
Therefore there are as many $Y$-s satisfying (\ref{240}) and (\ref{241}). As the hypotheses (\ref{239})-(\ref{241}) are local in character 
and intrinsic we may consider only the case when $M = U$ is open in $V$ vector space and $P(x)\in V^{\ast}\otimes V\,,\, P(x)^{2} = P(x)\,,
\;\forall x\in U$. Let then
\begin{equation}\label{245}
 W = (x,y;u,v)\,,\; Z = (\xi,\eta;\varphi,\psi)
\end{equation}
so that
\begin{equation}\label{246}
 TP\cdot W = (x,P(x)y;u,<P^{\prime}(x);u>y + P(x)v),
\end{equation}
$T(I-P)\cdot Z = (\xi,(I-P(\xi))\eta;\varphi,-<P^{\prime}(\xi);\varphi>\eta + (I-P(\xi))\psi)$ and therefore
\begin{equation}\label{247}
\omega_{M}(T(I-P)\cdot Z) = (\xi,\varphi;(I-P(\xi))\eta,-<P^{\prime}(\xi);\varphi>\eta + (I-P(\xi))\psi).
\end{equation}
Then the equality (\ref{239}) gives
\begin{equation}\label{248}
 \xi = x,\;\varphi = P(x)y,\;u = (I-P(x))\eta,
\end{equation}
$<P^{\prime}(x);(I-P(x))\eta>y + P(x)v =-<P^{\prime}(x);P(x)y>\eta + (I-P(x))\psi$.\\
Applying here $P(x)$ on both sides we get\\
$P(x)v = -P(x)<P^{\prime}(x);P(x)y>\eta - P(x)<P^{\prime}(x);(I-P(x))\eta>y$,\\
so that finally
\begin{multline}\label{249}
 TP\cdot W = (x,P(x)y;(I-P(x))\eta,(I-P(x))<P^{\prime}(x);(I-P(x))\eta>y\, -\\
 - P(x)<P^{\prime}(x);P(x)y>\eta).
\end{multline}
From (\ref{248}) we see also that
\begin{equation}\label{250}
 W = (x,y;(I-P(x))\eta,v)\,,\; Z = (x,\eta;P(x)y,\psi).
\end{equation}
Therefore, with the definition (\ref{240})-(\ref{241}) we find that
\begin{equation}\label{251}
 Y = (x,P(x)y+(I-P(x))\eta;P(x)y+(I-P(x))\eta,r)
\end{equation}
with $r\in V$ arbitrary chosen. Then\\
$T(I-P)\cdot Y = (x,(I-P(x))\eta;P(x)y+(I-P(x))\eta,\\
- <P^{\prime}(x);P(x)y+(I-P(x))\eta>(P(x)y+(I-P(x))\eta) + (I-P(x))r)$,\\
$\omega_{M}(T(I-P)\cdot Y) = (x,P(x)y+(I-P(x))\eta;(I-P(x))\eta,\\
- <P^{\prime}(x);P(x)y+(I-P(x))\eta>(P(x)y+(I-P(x))\eta) + (I-P(x))r)$\\
and finally
\begin{multline}\label{252}
 TP(\omega_{M}(T(I-P)\cdot Y)) = (x,P(x)y;(I-P(x))\eta,\\
 <P^{\prime}(x);(I-P(x))\eta>(P(x)y+(I-P(x))\eta) - \\
 - P(x)<P^{\prime}(x);P(x)y+(I-P(x))\eta>(P(x)y+(I-P(x))\eta).
\end{multline}
Let us compare now (\ref{252}) with (\ref{249}). Taking the derivative in the identity $P(x)^{2} = P(x)$ we get first
\begin{equation}\label{253}
 P(x)<P^{\prime}(x);a>P(x)b = 0\,,\;\forall a, b\in V,
\end{equation}
and using it for $I-P$ instead of $P$
\begin{equation}\label{254}
 (I-P(x))<P^{\prime}(x);a>(I-P(x))b = 0.
\end{equation}
These are the keys to check the equality of the fourth terms in the expressions of $TP\cdot W$ and $ TP(\omega_{M}(T(I-P)\cdot Y))$. 
If we transpose now $W$ with $Z$ and also $P$ with $I-P$, we see that the hypothesis (\ref{239}) is still satisfied and according to it 
definition (\ref{240})-(\ref{241}) works with the same $Y$. Therefore $Y$ will verify also $(T(I-P)\circ\omega_{M}\circ TP)(Y) = 
T(I-P)\cdot Z =\\
= \omega_{M}(TP\cdot W)$, wherefrom the third equality in (\ref{242}).\\
We prove now that, as soon as $Y\in T(TM)$ satisfies (\ref{240}) it also verifies
\begin{equation}\label{255}
(TP\circ\omega_{M}\circ T(I-P))(Y) = (\omega_{M}\circ T(I-P) \circ\omega_{M}\circ TP)(Y). 
\end{equation}
Likewise, if locally,
\begin{equation}\label{256}
 Y = (x,y;y,z)\,,\; x\in U\,,\;y, z\in V,
\end{equation}
comparing with (\ref{252}) we get
\begin{multline}\label{257}
 (TP\circ\omega_{M}\circ T(I-P))(Y) = ((x,P(x)y;(I-P(x))y,<P^{\prime}(x);(I-P(x))y>y +\\
 + P(x)<(I-P)^{\prime}(x);y>y
\end{multline}
and then interchanging $P$ and $I-P$ and applying $\omega_{M}$ on both sides we come to
\begin{multline}\label{258}
(\omega_{M}\circ T(I-P) \circ\omega_{M}\circ TP)(Y) = (x,P(x)y;(I-P(x))y,<(I-P)^{\prime}(x);P(x)y>y +\\
 + (I-P(x))<P^{\prime}(x);y>y.
\end{multline}
It is easy to check the identity of the fourth terms in the right hand sides of these last two relations.\\
Let us define then, for $X\in TM$:
\begin{multline}\label{259}
 R_{P}(X) = (TP\circ\omega_{M}\circ T(I-P))(Y) = (\omega_{M}\circ T(I-P) \circ\omega_{M}\circ TP)(Y)
\end{multline}
if $Y\in T(TM)$ satisfies
\begin{equation}\label{260}
 T\tau_{M}(Y) = \tau_{TM}(Y) = X.
\end{equation}
From the computations above it results that $R_{P}(X)$ does not depend on $Y$ chosen with (\ref{260}) (see (\ref{256}), (\ref{257})), that
$R_{P} : TM\longrightarrow TK\cap\omega_{M}(TH)$, is surjective and satisfies (\ref{233}) and (\ref{234}). We see also that $R_{P}$ 
is given locally by
\begin{multline}\label{261}
 R_{P}(x;y) = ((x,P(x)y;(I-P(x))y,<P^{\prime}(x);(I-P(x))y>y +\\
  + P(x)<(I-P)^{\prime}(x);y>y
\end{multline}
(see (\ref{257}). Since $\tau_{TM} + T\tau_{M} : T_{k}K\cap\omega_{M}(TH)\longrightarrow \{k\}+ H_{\tau_{M}(k)}$ is a bijection and 
$T\tau_{M}$ is linear, it results that $T_{k}K\cap\omega_{M}(TH)$ is a vector subspace of $T_{k}K$ and that
\begin{equation}\label{262}
 T_{k}\tau_{M} : T_{k}K\cap\omega_{M}(TH)\widetilde{\longrightarrow}  H_{\tau_{M}(k)}
\end{equation}
is an isomorphism of inverse $R_{P}(k+(\cdot))$ (see (\ref{235})). Let us verify now (\ref{238}). If locally $k = (x;y)$ with $P(x)y = y$, 
we find (see (\ref{237})):
\begin{equation}\label{263}
 N_{(x;y)} = \{(x,y;z,w)\,\rvert\,P(x)z = 0\,,\;(I-P(x))w = <P^{\prime}(x);z>y\}
\end{equation}
and, for $P(x)y = y\,,\; P(x)z = 0$\,:
\begin{equation}\label{264}
 (\omega_{M}\circ TP\circ\omega_{M})(x,y;z,w) = (x,y;0,w\, + <P^{\prime}(x);y>z\,- <P^{\prime}(x);z>y),
\end{equation}
\begin{equation}\label{265}
 R_{P}(x;y+z) = (x,y;z,<P^{\prime}(x);z>y \,- <P^{\prime}(x);y>z)
\end{equation}
(see (\ref{261})). Note that $(I-P(x))<P^{\prime}(x);y>z = 0$ and $ P(x)<P^{\prime}(x);z>y = 0$, in virtue of (\ref{253}) and (\ref{254}), and
therefore in (\ref{264})\\
$(I-P(x))(w\, + <P^{\prime}(x);y>z\,- <P^{\prime}(x);z>y) = 0$.\\
The proof is complete $\blacksquare$\\
Our aim is to show that there exists a canonical linear connection on the restriction of the vector bundle $K$ to a path tangent to $H$, 
therefore defining a parallel transport in $K$ along it.\\
Recall that an Ehresmann connection on a nonlinear fiber bundle is given by the \textit{horizontal vector bundle}, in fact a smooth supplementary 
vector subbundle to the tangent at the nonlinear fiber in the tangent to the total space of the fiber bundle. The vector subbundle of the 
tangents to the nonlinear fibers is also called the vertical vector bundle. If\\
$\pi : \eta\longrightarrow D$ is the nonlinear fiber 
bundle, the connection is then defined by the linear projections
\begin{equation}\label{266}
 V_{e} : T_{e}\eta\longrightarrow T_{e}\eta\,,\;\mathrm{Im}\, V_{e} = \mathrm{Ker}\, T_{e}\pi = T_{e}\eta_{\pi (e)}
\end{equation}
on the vertical tangent spaces, smoothly depending on $e\in\eta$.\\
In the case of a vector bundle $\pi : E\longrightarrow D$, the mapping 
$T\pi : TE\longrightarrow TD$ delimits itself a \textit{tangent vector bundle structure} defined as follows: let 
$s : E\times_{D}E\longrightarrow E$ denote the sum in the fibers of $E$ and $Ts : T( E\times_{D}E)\longrightarrow TE$. Remark that
$T_{(e,f)}( E\times_{D}E) = T_{e}E\times_{T_{p}D} T_{f} E$, if $\pi(e) = \pi(f) = p$. Then $T_{(e,f)}s\cdot(X,Y)$ is defined in the case that
$T_{e}\pi\cdot X = T_{f}\pi\cdot Y\in T_{p}D$, and is precisely the sum between $X$ and $Y$ in the fiber of $T\pi$.\\
The connection on $\pi : E\longrightarrow D$ is linear in the case that the vertical projection $V$ enjoys of a second linearity property 
in the following sense. We remark that, in virtue of (\ref{266}) 
$T_{e}\pi\cdot V_{e} Z = 0_{T_{\pi(e)}D}\,,\;\forall Z\in T_{e}E$, and if $T_{e}\pi\cdot Z = X$ then 
$V_{e} Z\in (T\pi)^{-1}(\{0_{T_{\tau_{D}(X)}D}\})$ since $\tau_{D}(X) = \pi(e)$. So, to be linear the connection, $V$ should be also linear 
\begin{equation}\label{267}
 V : (T\pi)^{-1}(\{X\})\longrightarrow (T\pi)^{-1}(\{0_{T_{\tau_{D}(X)}D}\})
\end{equation}
$\forall X\in TD$, with respect to the tangent vector bundle structure. The linear connection defines the covariant derivative of sections
$\sigma$ of $\pi$ with respect to vectors $X\in T_{p}D$ by 
\begin{equation}\label{268}
 (\nabla_{X}\sigma)_{p} = [\varPsi^{E_{p}}(\sigma_{p},\cdot)]^{-1}\cdot V_{\sigma_{p}}\cdot T_{p}\sigma\cdot X.
\end{equation}
In this way, for $X\in T_{p}D$ we have $(\nabla_{X}\sigma)_{p}\in E_{p}$. Of course, for a vector field $X$ we get a new section 
$\nabla_{X}\sigma$ of $\pi$. It satisfies 
\begin{equation}\label{269}
\nabla_{fX}\sigma = f\cdot\nabla_{X}\sigma 
\end{equation}
and
\begin{equation}\label{270}
 \nabla_{X}f\sigma = f\cdot\nabla_{X}\sigma + (\textsl{L}_{X}f)\cdot\sigma,
\end{equation}
for $f$ scalar function, and is biadditive with respect to $X$ and $\sigma$. Conversely, such a bilinear operator with (\ref{269}) and (\ref{270})
defines the vertical projection $V$ from (\ref{268}) and thus the linear connection.\\
In our special case, the covariant derivative arises as the linearization of a vector field of finite differential order on 
$\textsl{C}^{\infty}(D,M)$ in a critical point. From the local definition (\ref{88}) we see that, if $\sigma_{0}$ is such that 
\begin{equation}\label{271}
 X_{\sigma_{0}} = 0_{T_{\sigma_{0}}\textsl{C}^{\infty}\Gamma(\eta)}, 
\end{equation}
it is defined its \textit{linearization at the critical point} $\sigma_{0}$ 
\begin{equation}\label{272}
 X^{\prime}(\sigma_{0}) : \textsl{C}^{\infty}\Gamma(\sigma^{\ast}_{0}(TF\eta))\longrightarrow\textsl{C}^{\infty}\Gamma(\sigma^{\ast}_{0}(TF\eta))
\end{equation}
by
\begin{equation}\label{273}
  X^{\prime}(\sigma_{0})\cdot W = \big[ Y,X\big]_{\sigma_{0}}\,,\;Y_{\sigma_{0}} = W;
\end{equation}
(see (\ref{65}) where we identify $T_{\sigma}\textsl{C}^{\infty}\Gamma(\eta)$ with $\textsl{C}^{\infty}\Gamma(\sigma^{\ast}(TF\eta))$). This 
means that the definition does not depend on the vector field $Y$ of finite order chosen as an extension of $W$ in a neighbourhood of 
$\sigma_{0}$.

\begin{theorem}. Let $H, K, P$ be as before smooth on $M$ and $D$ be a smooth compact curve (i.e. $\mathrm{dim}\,D = 1$). If 
$\beta_{0} : D\longrightarrow M$ is a smooth path tangent to $H$, i.e.
\begin{equation}\label{274}
 T_{\zeta}\beta_{0}\cdot\xi\in H_{\beta_{0}(\zeta)}\,,\;\forall\zeta\in D\,,\;\forall\xi\in T_{\zeta} D,
\end{equation}
then on the vector bundle $p : \beta_{0}^{\ast} K\longrightarrow D$ there is a canonical linear connection given by the horizontal vector 
subbundle
\begin{equation}\label{275}
 HT_{(\zeta,k)}\beta_{0}^{\ast} K =: (T_{(\zeta,k)} q)^{-1}(T_{k}K\cap\omega_{M}(TH)),
\end{equation}
where $q : \beta_{0}^{\ast} K\longrightarrow K$ is canonical (see(\ref{235})). Let $\xi$ be any smooth vector field on $D$.For the vector field
\begin{equation}\label{276}
 (l(P)\cdot r(\xi))_{\beta}(\zeta) = P_{\beta(\zeta)}\cdot T_{\zeta}\beta\cdot \xi_{\zeta}\,,\;\beta\in\textsl{C}^{\infty}(D,M)\,,\;\zeta\in D,
\end{equation}
(see (\ref{215})) $\beta_{0}$ is a critical point, 
$\textsl{C}^{\infty}\Gamma(\beta_{0}^{\ast} K)\subseteq\textsl{C}^{\infty}\Gamma(\beta_{0}^{\ast}(TM))$ is an invariant subspace for the 
linearization $(l(P)\cdot r(\xi))^{\prime}(\beta_{0})$ and the covariant derivative corresponding to the linear connection (\ref{275}) is
\begin{equation}\label{277}
 \nabla_{\xi}\sigma = (l(P)\cdot r(\xi))^{\prime}(\beta_{0})\cdot\sigma\,,\;\forall\sigma\in\textsl{C}^{\infty}\Gamma(\beta_{0}^{\ast} K).
\end{equation}
\end{theorem}

\textbf{Proof}. We start by computing the linearization $(l(P)\cdot r(\xi))^{\prime}(\beta_{0})$ using the definition (\ref{273}). As 
$l(P)\cdot r(\xi)$ is a vector field of differential order $1$ we will apply Theorem 3, having in mind to consider the extension 
$Y$ of $W$ also of finite differential order. Suppose then that in (\ref{178}), written for $\big[\varTheta,\varXi\big]$, interchanging also
$r$ with $s$, we have $X_{\sigma_{0}} = 0$, i.e. $\varXi(j^{r}_{\zeta}\sigma_{0}) = 0\,,\;\forall\zeta\in D$ (see also (\ref{181})). Then we
find that\\
$(T\varTheta\circ(\varOmega^{s})^{-1})(j^{s}_{\zeta}(\varXi(j^{r}\sigma_{0})))= 0_{T_{Y_{\sigma_{0}}(\zeta)}(T\eta_{\zeta})}$,
hence zero in the space $T_{Y_{\sigma_{0}}(\zeta)}(T\eta_{\zeta})$ where the difference is taken. Therefore, in the case that 
$X_{\sigma_{0}} = 0$, we have
\begin{equation}\label{278}
 \big[\varTheta,\varXi\big](j^{r+s}_{\zeta}\sigma_{0}) = 
 [\varPsi^{T_{\sigma_{0\zeta}}\eta_{\zeta}}(\varTheta(j^{s}_{\zeta}\sigma_{0}),\cdot)]^{-1}\cdot
 \omega_{\eta_{\zeta}}((T\varXi\circ(\varOmega^{r})^{-1}\circ J^{r}(\varTheta)\circ f_{rs})(j^{r+s}_{\zeta}\sigma_{0})).
\end{equation}
Consider now the case when $\eta = D\times M,\;\mathrm{dim}\,D = 1,\;r =1$ and let $s = 0$, that is
\begin{equation}\label{279}
 Y_{\beta}(\zeta) = \varTheta(\zeta,\beta(\zeta))\,,\;\zeta\in D\,,\; \beta(\zeta))\in M.
\end{equation}
Remember that
\begin{equation}\label{280}
 X_{\beta}(\zeta) = P_{\beta(\zeta)}\cdot T_{\zeta}\beta\cdot\xi_{\zeta}.
\end{equation}
The vector field $\xi$ on $D$ defines the fiber preserving mappings
\begin{multline}\label{281}
 \varepsilon^{0}_{\xi} : J^{1}(D\times M)\longrightarrow D\times TM,\\
 (\varepsilon^{0}_{\xi})_{\zeta} : J^{1}_{\zeta}(D\times M)\longrightarrow TM\,,\;
  (\varepsilon^{0}_{\xi})_{\zeta}(j^{1}_{\zeta}\beta) = T_{\zeta}\beta\cdot\xi_{\zeta},
\end{multline}
\begin{multline}\label{282}
 \varepsilon^{1}_{\xi} : J^{1}(D\times TM)\longrightarrow D\times T(TM),\\
(\varepsilon^{1}_{\xi})_{\zeta} : J^{1}_{\zeta}(D\times TM)\longrightarrow  T(TM)\,,\;
  (\varepsilon^{1}_{\xi})_{\zeta}(j^{1}_{\zeta} X) = T_{\zeta} X\cdot\xi_{\zeta}.
\end{multline}
Therefore 
\begin{equation}\label{283}
 X_{\beta}(\zeta) = \varXi(j^{1}_{\zeta}\beta)\,,\;\varXi = P\circ\varepsilon^{0}_{\xi}.
\end{equation}
As $TF(D\times M) = D\times TM$ we have
\begin{equation}\label{284}
 \varOmega^{1}_{\zeta} : T(J^{1}_{\zeta}(D\times M))\longrightarrow J^{1}_{\zeta}(D\times TM)
\end{equation}
(see (\ref{158})). Then the following relation holds
\begin{equation}\label{285}
 T(\varepsilon^{0}_{\xi})_{\zeta}\circ (\varOmega^{1}_{\zeta})^{-1} = \omega_{M}\circ(\varepsilon^{1}_{\xi})_{\zeta}\,,\;\forall\zeta\in D. 
\end{equation}
This can be easily verified first when $M = U$ open in a vector space $V$, wherefrom in the general case. Then in our case the formula (\ref{278})
reads\\
$\big[\varTheta,\varXi\big](j^{1}_{\zeta}\beta_{0}) = 
 [\varPsi^{T_{\beta_{0}(\zeta)}M}(\varTheta(\zeta,\beta_{0}(\zeta)),\cdot)]^{-1}\cdot
 \omega_{M}((TP\circ T(\varepsilon^{0}_{\xi})_{\zeta}\circ(\varOmega^{1}_{\zeta})^{-1}) (j^{1}(\varTheta))(\zeta,\beta_{0}(\zeta))) =$\\
$=  [\varPsi^{T_{\beta_{0}(\zeta)}M}(\varTheta(\zeta,\beta_{0}(\zeta)),\cdot)]^{-1}\cdot 
 (\omega_{M}\circ TP\circ\omega_{M})((\varepsilon^{1}_{\xi})_{\zeta} (j^{1}(\varTheta(\zeta,\beta_{0}(\zeta))))) =$\\
$=  [\varPsi^{T_{\beta_{0}(\zeta)}M}(\sigma_{\zeta},\cdot)]^{-1}\cdot 
 (\omega_{M}\circ TP\circ\omega_{M})(T_{\zeta}\sigma\cdot\xi_{\zeta})$ \\
where  
\begin{equation}\label{286}
 \sigma_{\zeta} =  \varTheta(\zeta,\beta_{0}(\zeta)) =  Y_{\beta_{0}}(\zeta)\in T_{\beta_{0}}M.
\end{equation}
We thus found, for $X = l(P)\cdot r(\xi)$ and $X_{\beta_{0}} = 0$
\begin{equation}\label{287}
 (X^{\prime}(\beta_{0})\cdot\sigma)(\zeta) =  [\varPsi^{T_{\beta_{0}(\zeta)}M}(\sigma_{\zeta},\cdot)]^{-1}\cdot 
 (\omega_{M}\circ TP\circ\omega_{M})\cdot T_{\zeta}\sigma\cdot\xi_{\zeta}
\end{equation}
for $\sigma\in\textsl{C}^{\infty}\Gamma(\beta_{0}^{\ast}(TM)) = T_{\beta_{0}}\textsl{C}^{\infty}(D,M)$.\\
 In the case that $X_{\beta_{0}}(\zeta) = 0\,,\;\forall \zeta$, and $Y_{\beta_{0}}(\zeta)\in K_{\beta_{0}(\zeta)}\,,\;
\forall\zeta$, from Theorem 2 we infer that $\big[X,Y\big]_{\beta_{0}}(\zeta)\in K_{\beta_{0}(\zeta)}\,,\;\forall\zeta\in D$. Hence 
$\textsl{C}^{\infty}\Gamma(\beta_{0}^{\ast} K)$ is invariant for $X^{\prime}(\beta_{0})$.\\
 In the formula above, however, $\sigma$ is seen as a function 
$\sigma : D\longrightarrow K$, hence when $\sigma$ is thought as a section of $p : \beta_{0}^{\ast}K\longrightarrow D$ we have to consider 
instead $q\circ\sigma$. Next we verify that
$T_{(\zeta,k)}q :T_{(\zeta,k)}\beta_{0}^{\ast}K\longrightarrow N_{k}\,,\;\forall (\zeta,k)\in\beta_{0}^{\ast}K$.
But $N_{k}$ is invariant for $\omega_{M}\circ TP\circ\omega_{M}$, which projects $N_{k}$ on 
$T_{k}K_{\beta_{0}(\zeta)}$, according to Theorem 6. And we have to take into account also the isomorphism 
$(TF_{(\zeta,k)}q)^{-1} : T_{k}K_{\beta_{0}(\zeta)}\widetilde{\longrightarrow}TF_{(\zeta,k)}\beta_{0}^{\ast}K$, in order to come back in 
$\beta_{0}^{\ast}K$. Finally, for $\sigma\in\textsl{C}^{\infty}\Gamma(\beta_{0}^{\ast} K)$ we 
may replace in (\ref{287}) $[\varPsi^{T_{\beta_{0}(\zeta)}M}(\sigma_{\zeta},\cdot)]^{-1}$ by  
$[\varPsi^{K_{\beta_{0}(\zeta)}}(\sigma_{\zeta},\cdot)]^{-1}$ to get
\begin{equation}\label{288}
  (X^{\prime}(\beta_{0})\cdot\sigma)(\zeta) =  [\varPsi^{(\beta_{0}^{\ast}K)_{\zeta}}(\sigma_{\zeta},\cdot)]^{-1}\cdot 
(TF_{\sigma_{\zeta}}q)^{-1}\cdot(\omega_{M}\circ TP\circ\omega_{M})\cdot T_{\sigma_{\zeta}}q\cdot T_{\zeta}\sigma\cdot\xi_{\zeta},
\end{equation}
and then define
\begin{equation}\label{289}
V_{(\zeta,k)} =  (TF_{(\zeta,k)}q)^{-1}\cdot(\omega_{M}\circ TP\circ\omega_{M})\cdot T_{(\zeta,k)}q\,,\;(\zeta,k)\in\beta_{0}^{\ast}K.
\end{equation}
The formula (\ref{288}), of the form (\ref{268}), gives indeed a covariant derivative on $\beta_{0}^{\ast}K\longrightarrow D$, for the vertical 
projection $V_{(\zeta,k)}$ above, if this $V$  is linear with respect to the tangent vector bundle structure (see (\ref{267})). We note first 
that, for any vector bundle morphism $\varPhi :E\longrightarrow F$, the tangent mapping $T\varPhi : TE\longrightarrow TF$ is a vector bundle 
morphism with respect to the tangent vector bundle structures; next, the fact that\\
$\omega_{M} : T(TM)\longrightarrow T(TM)$ is an isomorphism between the two vector bundle structures delimited by $\tau_{TM}$ and $T\tau_{M}$ on 
$T(TM)$. In (\ref{289})
$q : \beta_{0}^{\ast}K\longrightarrow TM$ is a vector bundle morphism and then $Tq$ is a vector bundle morphism from $Tp$ to $T\tau_{M}$.
Next, $TP$ is a vector bundle morphism from $\tau_{TM}$ to $\tau_{TM}$ and then $\omega_{M}\circ TP\circ\omega_{M}$ is a vector bundle morphism 
from $T\tau_{M}$ to $T\tau_{M}$. Therefore $\omega_{M}\circ TP\circ\omega_{M}\circ Tq$ is a vector bundle morphism from $Tp$ to $T\tau_{M}$. More 
precisely, it gives a morphism from $(Tp)^{-1}(\{\xi\})$ to $(T\tau_{M})^{-1}(0_{T_{\beta_{0}(\zeta)}M})$, if $\xi\in T_{\zeta}D$. As
$\iota : K\longrightarrow TM$ is a vector bundle embedding, $(T\tau_{M}\rvert_{K})^{-1}(0_{T_{\beta_{0}(\zeta)}M})$ is a subspace in  
$(T\tau_{M})^{-1}(0_{T_{\beta_{0}(\zeta)}M})$ with the respective vector space structure and contains the image of that morphism. We remark that
 $(T\tau_{M}\rvert_{K})^{-1}(0_{T_{\beta_{0}(\zeta)}M}) = T(K_{\beta_{0}(\zeta)})$. (Note that, in general, the vector space structure of $TV$, 
given by $\varPsi^{V}$, for $V$ vector space (see p. 24) coincides with the tangent vector space structure if $V = E_{\zeta}$ and 
$\pi : E\longrightarrow D$ is a vector bundle so that $TV = (T\pi)^{-1}(\{0_{T_{\zeta}D}\})$). We thus found that 
$\omega_{M}\circ TP\circ\omega_{M}\circ Tq$ gives a vector space morphism from $(Tp)^{-1}(\{\xi\})$ to $T(K_{\beta_{0}(\zeta)})$, if 
$\xi\in T_{\zeta}D$. As $K_{\beta_{0}(\zeta)} = (\beta_{0}^{\ast}K)_{\zeta}$, $V_{(\zeta,k)}$ is a morphism from $(Tp)^{-1}(\{\xi\})$ to
$(Tp)^{-1}(\{0_{T_{\zeta}D}\})$.\\
The covariant derivative
\begin{equation}\label{290}
 (\nabla_{\xi}\sigma)_{\zeta} = [\varPsi^{(\beta_{0}^{\ast}K)_{\zeta}}(\sigma_{\zeta},\cdot)]^{-1}\cdot 
(TF_{\sigma_{\zeta}}q)^{-1}\cdot(\omega_{M}\circ TP\circ\omega_{M})\cdot T_{\sigma_{\zeta}}q\cdot T_{\zeta}\sigma\cdot\xi_{\zeta}
\end{equation}
corresponds to the horizontal subbundle
\begin{equation}\label{291}
 HT_{(\zeta,k)}\beta_{0}^{\ast}K = \mathrm {Ker}\, V_{(\zeta,k)}
\end{equation}
(see (\ref{289})). It remains only to identify this space from (\ref{235}) and (\ref{238}) of Theorem 6, getting so (\ref{275}) (recall 
also (\ref{230})) $\blacksquare$

\textbf{Remark}. In spite of its form (\ref{287}), $X^{\prime}(\beta_{0})\cdot\sigma$ does not give a covariant derivative on 
$\beta_{0}^{\ast}(TM)\longrightarrow D$ with respect to $\xi$, since, for $\sigma\in\textsl{C}^{\infty}\Gamma(\beta_{0}^{\ast}(TM))$ and $f$ 
smooth scalar function, we have\\
$(X^{\prime}(\beta_{0})\cdot(f\cdot\sigma))(\zeta) = f(\zeta)\cdot (X^{\prime}(\beta_{0})\cdot\sigma)(\zeta) +  
(\textsl{L}_{\xi}f)(\zeta)\cdot P(\sigma_{\zeta})$,\\
where $P$ acts on $TM$ as in (\ref{231}). Compare with the required property (\ref{270}). Of course $P(\sigma_{\zeta}) = \sigma_{\zeta}$,
for $\sigma_{\zeta}\in K_{\beta_{0}(\zeta)}$ $\blacksquare$
 
 In the case of an Ehresmann connection on the nonlinear fiber bundle $\pi : \eta\longrightarrow B$, we consider the transversal vector subbundle
\begin{equation}\label{292}
 K_{e} := \mathrm{ker}\,T_{e}\pi = TF_{e}\eta = T_{e}\eta_{\pi(e)}\,,\;e\in\eta,
\end{equation}
to the horizontal vector subbundle $H$ of $T\eta$ (in $M = \eta$). If $\beta_{0} : I\longrightarrow\eta$ is a horizontal path (on a compact 
interval $I$ of $\textbf{R})$, i.e. $\dot{\beta_{0}}(t)\in H_{\beta_{0}(t)}\,,\;\forall t\in I$, we consider the path in the base $B$
\begin{equation}\label{293}
 \gamma := \pi\circ\beta_{0}
\end{equation}
and the parallel transport along it
\begin{equation}\label{294}
 \gamma_{s}^{t} : G_{\beta_{0}(s)}\longrightarrow G_{\beta_{0}(t)}
\end{equation}
from a neighbourhood $G_{\beta_{0}(s)}$ of $\beta_{0}(s)$ in $\eta_{\gamma(s)}$ to a neighbourhood $G_{\beta_{0}(t)}$ of $\beta_{0}(t)$ in
$\eta_{\gamma(t)}$, for $s, t\in I$. On the other hand, according to the previous Theorem 7, we have a linear parallel transport
\begin{equation}\label{295}
 \tau_{s}^{t}\in\mathrm{Hom}(K_{\beta_{0}(s)},K_{\beta_{0}(t)})\,,\;s, t\in I,
\end{equation}
determined by the linear connection in $\beta_{0}^{\ast}K\longrightarrow I$. Then we have the following result
\begin{theorem}. In the hypothesis (\ref{292}), with the notations (\ref{293}), (\ref{294}) and (\ref{295}), we have
\begin{equation}\label{296}
 \tau_{s}^{t} = T_{\beta_{0}(s)}\gamma_{s}^{t},
\end{equation}
as $K_{\beta_{0}(t)} = T_{\beta_{0}(t)}\eta_{\gamma(t)}\,,\;\forall t\in I$. It results that, for $\eta = E$ vector bundle with a linear 
connection and $\nabla$ the corresponding covariant derivative, if we denote $\nabla^{\beta_{0}^{\ast}K}$ the covariant derivative (\ref{277}) on
$\textsl{C}^{\infty}\Gamma(\beta_{0}^{\ast}K)$, the following equality holds
\begin{equation}\label{297}
 \varPsi^{E_{\gamma(t)}}(\beta_{0}(t),\nabla_{\dot{\gamma}(t)}\alpha_{\gamma(t)}) = 
 \nabla^{\beta_{0}^{\ast}K}_{\frac{\partial}{\partial t}}\varPsi^{E_{\gamma(t)}}(\beta_{0}(t),\alpha_{\gamma(t)})
\end{equation}
for $\alpha\in\textsl{C}^{\infty}\Gamma(E)$. Note that 
$\varPsi^{E_{\gamma(t)}}(\beta_{0}(t),\alpha_{\gamma(t)})\,,\; \varPsi^{E_{\gamma(t)}}(\beta_{0}(t),\nabla_{\dot{\gamma}(t)}\alpha_{\gamma(t)})
\in T_{\beta_{0}(t)} E_{\gamma(t)} =\\
= K_{\beta_{0}(t)}$.
\end{theorem}

\textbf{Proof}. The multiplicative properties
\begin{equation}\label{298}
 \gamma_{s}^{t} = \gamma_{u}^{t}\circ\gamma_{s}^{u}\,,\;\tau_{s}^{t} = \tau_{u}^{t}\cdot\tau_{s}^{u}\,,\; s, u, t\in I,
\end{equation}
allow to reduce proving (\ref{296}) for $s$ and $t$ as close to each other as we want. In this case, we can take $\eta = D\times U$, $D$ open in 
$V$, $U$ open in $W$, $V, W$ vector spaces, $B = D$, $\pi: D\times U\longrightarrow D$ canonical. Also $H_{(x,y)} = 
\mathrm{graph}\,\varGamma(x,y)\,,\;\varGamma(x,y)\in\mathrm{Hom}(V,W)$, as in (\ref{6}). Then, if
\begin{equation}\label{299}
 \beta_{0}(t) = (x(t),y(t))\,,\;\gamma (t) = x(t),
\end{equation}
$\gamma_{s}^{t}(y)$ is the solution of the problem in $W$
\begin{equation}\label{300}
 \dfrac{\partial}{\partial t}\gamma_{s}^{t}(y) = \varGamma(x(t),\gamma_{s}^{t}(y))\cdot \dot{x}(t)\,,\;\gamma_{s}^{s}(y) = y\,,\; y\in U.
\end{equation}
The hypothesis on $\beta_{0}$ being horizontal gives
\begin{equation}\label{301}
 \gamma_{s}^{t}(y(s)) = y(t)\,,\;\forall t.
\end{equation}
Taking the derivative with respect to $y$ in (\ref{300}), in $y=y(s)$, and taking into account (\ref{301}), we find that 
$w(t) := <\dfrac{\partial \gamma_{s}^{t}}{\partial y}(y(s));w_{0}>$ is the unique solution of the problem
\begin{equation}\label{302}
 \dot{w}(t) = <\dfrac{\partial\varGamma}{\partial y}(x(t),y(t));w(t)>\cdot\dot{x}(t)\,,\;w(s) = w_{0}.
\end{equation}
In this case\\
$K = \{(x,y;0_{V},w)\rvert\, x\in D, y\in U, w\in W\}$\,,\;
$H = \{(x,y;v,\varGamma(x,y)v)\rvert\,x\in D, y\in U, v\in V\}$\\
and we find that
\begin{equation}\label{303}
 T_{(x,y;0_{V},w)}K\cap\omega_{D\times U}(TH) = 
 \{(x,y,0_{V},w;v,\varGamma(x,y)v,0_{V},<\dfrac{\partial\varGamma}{\partial y}(x,y);w>v)\,\rvert\,v\in V\}.
\end{equation}
Since in trivialization
\begin{equation}\label{304}
 \beta_{0}^{\ast}K \tilde{=} \{(t,w)\,\rvert\,t\in I, w\in W\}
\end{equation}
we have
\begin{equation}\label{305}
 q(t,w) = (x(t),y(t);0_{V},w)
\end{equation}
and
\begin{equation}\label{306}
 Tq(t,w;\theta,\kappa) = (x(t),y(t);0_{V},w;\theta\dot{x}(t),\theta\dot{y}(t),0_{V},\kappa).
\end{equation}
Therefore, according to (\ref{275}) 
\begin{equation}\label{307}
 HT_{(t,w)}\beta_{0}^{\ast} K = 
 \{(t,w;\theta,\theta<\dfrac{\partial\varGamma}{\partial y}(x(t),y(t));w>\cdot\dot{x}(t)\,\rvert\,\theta\in\textbf{R}\}.
\end{equation}
Then $(t,w(t))$ is horizontal in $\beta_{0}^{\ast} K$ if and only if $(1,\dot{w}(t))\in HT_{(t,w(t))}\beta_{0}^{\ast} K$, that is, when
$\dot{w}(t) = <\dfrac{\partial\varGamma}{\partial y}(x(t),y(t));w(t)>\cdot\dot{x}(t)$.
Comparing with (\ref{302}), we get (\ref{296}).\\
Recall that for a linear connection we have
\begin{equation}\label{308}
 \nabla_{\dot{\gamma}(t)}\alpha_{\gamma(t)} = \dfrac{\mathrm{d}}{\mathrm{d}h}\gamma_{t+h}^{t}\alpha_{\gamma(t+h)}\rvert_{h=0}.
\end{equation}
On the other hand, for $A\in\mathrm{Hom}\,(V,W)$
\begin{equation}\label{309}
 T_{v}A\cdot\varPsi^{V}(v,u) = \varPsi^{W}(Av,Au)\,,\;\forall v, u\in V.
\end{equation}
In this way\\

$\varPsi^{E_{\gamma(t)}}(\beta_{0}(t), \nabla_{\dot{\gamma}(t)}\alpha_{\gamma(t)}) = \varPsi^{E_{\gamma(t)}}(\beta_{0}(t),
\dfrac{\mathrm{d}}{\mathrm{d}h}\gamma_{t+h}^{t}\alpha_{\gamma(t+h)}\rvert_{h=0}) =$\\

$= \dfrac{\mathrm{d}}{\mathrm{d}h}\varPsi^{E_{\gamma(t)}}(\beta_{0}(t),\gamma_{t+h}^{t}\alpha_{\gamma(t+h)})\rvert_{h=0} =
 \dfrac{\mathrm{d}}{\mathrm{d}h}\varPsi^{E_{\gamma(t)}}(\gamma_{t+h}^{t}\beta_{0}(t+h),\gamma_{t+h}^{t}\alpha_{\gamma(t+h)})\rvert_{h=0} =$\\

$= \dfrac{\mathrm{d}}{\mathrm{d}h}T_{\beta_{0}(t+h)}\gamma_{t+h}^{t}\cdot\varPsi^{E_{\gamma(t+h)}}(\beta_{0}(t+h),\alpha_{\gamma(t+h)})
\rvert_{h=0} =$\\

$= \dfrac{\mathrm{d}}{\mathrm{d}h}\tau_{t+h}^{t}\varPsi^{E_{\gamma(t+h)}}(\beta_{0}(t+h),\alpha_{\gamma(t+h)})\rvert_{h=0} =
 \nabla^{\beta_{0}}_{\frac{\partial}{\partial t}}\varPsi^{E_{\gamma(t)}}(\beta_{0}(t),\alpha_{\gamma(t)})$,\\
 
 the outside derivatives being taken in $ T_{\beta_{0}(t)} E_{\gamma(t)} = K_{\beta_{0}(t)}$. This ends the proof $\blacksquare$\\

 \textbf{Remark}. We can take in (\ref{297}) 
\begin{equation}\label{310}
 \beta_{0}(t) = 0_{E_{\gamma(t)}}\;,
\end{equation}
for $\gamma : I\longrightarrow B$ arbitrary, since the connection is linear. Then the equality shows that $\nabla_{\dot{\gamma}(t)}$ is 
determined by $\nabla^{0_{E_{\gamma(t)}}}_{\frac{\partial}{\partial t}}$, the connection given in Theorem 7 being thus the right 
generalization of a linear connection on a vector bundle.\\
On the other hand, the local parallel transport maps (\ref{294}) are still defined and the equality (\ref{296}) holds under the only hypothesis 
that $C^{K} = 0$ in a neighbourhood in $M$ of the tangent to $H$ path $\beta_{0}(I)$  $\blacksquare$\\

Let us consider again $\pi : \eta\longrightarrow D$ a nonlinear fiber bundle, $B\subseteq D$ submanifold, both $B$ and $D$ compact with 
boundary, $\eta\rvert_{B}$ and 
\begin{equation}\label{311}
 \rho : \textsl{C}^{\infty}\Gamma(\eta)\longrightarrow\textsl{C}^{\infty}\Gamma(\eta\rvert_{B})\,,\;\rho(\sigma) := \sigma\rvert_{B},
\end{equation}
the natural restriction. It is easy to verify that for $\sigma\in \textsl{C}^{\infty}\Gamma(\eta)$
\begin{equation}\label{312}
 \sigma^{\ast}(TF\eta)\rvert_{B} = (\sigma\rvert_{B})^{\ast}(TF(\eta\rvert_{B}))
\end{equation}
and that, for $k\geqslant 0$
\begin{equation}\label{313}
 T_{\sigma}\rho :  T_{\sigma}\textsl{C}^{k}\Gamma(\eta)\longrightarrow T_{\sigma\rvert_{B}}\textsl{C}^{k}\Gamma(\eta\rvert_{B})\,,\;
 T_{\sigma}\rho\cdot X = X\rvert_{B}
\end{equation}
with the same meaning. We will need in the sequel the following
\begin{proposition}. Let $X^{D}\in\textsl{C}^{\infty}\Gamma_{\bullet} (T\textsl{C}^{\infty}\Gamma (\eta)),
X^{B}\in\textsl{C}^{\infty}\Gamma_{\bullet} (T\textsl{C}^{\infty}\Gamma (\eta\rvert_{B}))$ be smooth vector fields of finite order such that
\begin{equation}\label{314}
 X^{B}_{\sigma\rvert_{B}} = X^{D}_{\sigma}\rvert_{B}\,,\;\forall\sigma\in\textsl{C}^{\infty}\Gamma (\eta).
\end{equation}
If $\sigma_{0}$ is a critical point for $X^{D}$ then $\sigma_{0}\rvert_{B}$ is a critical point for $X^{B}$ and
\begin{equation}\label{315}
 <(X^{B})^{\prime}(\sigma_{0}\rvert_{B});Z\rvert_{B}> = <(X^{D})^{\prime}(\sigma_{0});Z>\rvert_{B}\,,\;\forall 
 Z\in\textsl{C}^{\infty}\Gamma(\sigma_{0}^{\ast}(TF\eta)).
\end{equation}
\end{proposition}
\textbf{Proof}. We can use a vector bundle neighbourhood of $\sigma_{0}(D)$ and the corresponding neighbourhood and chart on 
$\textsl{C}^{\infty}\Gamma (\eta)$ $\blacksquare$\\

Note that this fact is an infinite dimensional analogue of the following:\\

\textit{Let $M, N$ be smooth manifolds, $\rho : M\longrightarrow N$ smooth map, $X^{M}, X^{N}$ smooth vector fields, on $M, N$ respectively,
such that
\begin{equation}\label{316}
 X^{N}_{\rho(x)} = T_{x}\rho\cdot X^{M}_{x}\,,\;\forall x\in M.
\end{equation}
If $x_{0}$ is a critical point for $X^{M}$, then $\rho(x_{0})$ is a critical point for $X^{N}$ and
\begin{equation}\label{317}
 (X^{N})^{\prime}(\rho(x_{0}))\cdot T_{x_{0}}\rho = T_{x_{0}}\rho\cdot (X^{M})^{\prime}(x_{0})
\end{equation}
on $T_{x_{0}} M$.}\\

Remember that we denote the covariant derivative (\ref{277}) given by Theorem 7
\begin{equation}\label{318}
\nabla^{\beta_{0}^{\ast}K}_{\xi}\sigma = (l(P)\cdot r(\xi))^{\prime}(\beta_{0})\cdot\sigma\,,
\;\forall\sigma\in\textsl{C}^{\infty}\Gamma(\beta_{0}^{\ast} K). 
\end{equation}
Also, if $D = I,\;I$ interval in $\textbf{R}$ 
\begin{equation}\label{319}
 (\beta_{0}^{\ast}K)_{s}^{t}\in\mathrm{Hom}\,(K_{\beta_{0}(s)},K_{\beta_{0}(t)})\,,\;s, t\in I,
\end{equation}
will stand for the parallel transport determined by this linear connection.\\

On \textit{the infinitesimal variation of the tangent paths to $H$} we have first
\begin{theorem}. Let $H$ and $K$ be two supplementary subbundles of $TM$, $I$ compact interval of $\textbf{R}$, $s_{0}\in\textbf{R}$, 
$\varepsilon> 0$ and 
\begin{equation}\label{320}
 \beta : (s_{0}-\varepsilon,s_{0}+\varepsilon)\times I\longrightarrow M
\end{equation}
be a smooth map such that
\begin{equation}\label{321}
 \dfrac{\partial\beta}{\partial t}(s,t)\in H_{\beta(s,t)}\,,\;\forall s\in (s_{0}-\varepsilon,s_{0}+\varepsilon)\,,\;\forall t\in I.
\end{equation}
Then
\begin{equation}\label{322}
 \nabla_{\frac{\partial}{\partial t}}^{\beta(s_{0},\cdot)^{\ast}K} P^{H_{\beta(s_{0},t)}}_{K_{\beta(s_{0},t)}} 
 \dfrac{\partial\beta}{\partial s}(s_{0},t) = Q^{H_{\beta(s_{0},t)}}_{K_{\beta(s_{0},t)}}
 C^{H}_{\beta(s_{0},t)}( P^{K_{\beta(s_{0},t)}}_{H_{\beta(s_{0},t)}} \dfrac{\partial\beta}{\partial s}(s_{0},t)\wedge 
 \dfrac{\partial\beta}{\partial t}(s_{0},t)),
\end{equation}
or, equivalently
\begin{multline}\label{323}
 \dfrac{\partial}{\partial t}\{(\beta(s_{0},\cdot)^{\ast}K)^{t_{0}}_{t}\cdot P^{H_{\beta(s_{0},t)}}_{K_{\beta(s_{0},t)}} 
 \dfrac{\partial\beta}{\partial s}(s_{0},t)\} =\\
 = (\beta(s_{0},\cdot)^{\ast}K)^{t_{0}}_{t}Q^{H_{\beta(s_{0},t)}}_{K_{\beta(s_{0},t)}}
 C^{H}_{\beta(s_{0},t)}( P^{K_{\beta(s_{0},t)}}_{H_{\beta(s_{0},t)}} \dfrac{\partial\beta}{\partial s}(s_{0},t)\wedge 
 \dfrac{\partial\beta}{\partial t}(s_{0},t)).
\end{multline}
\end{theorem}

\textbf{Proof}. We take in Theorem 5 $X = \dfrac{\partial}{\partial s}\,,\; Y = \dfrac{\partial}{\partial t}$ on suitable domain 
$\varOmega$ with smooth boundary such that, for an 
$\eta > 0$, $(s_{0}-\eta,s_{0}+\eta)\times I\subset\bar{\varOmega}\subset(s_{0}-\varepsilon,s_{0}+\varepsilon)\times I$ and consider 
$D = \bar{\varOmega}$. In virtue of (\ref{321})
\begin{equation}\label{324}
 (l(P^{H}_{K})\cdot r(\dfrac{\partial}{\partial t}))_{\beta}(s,t) = 0\,,\;\forall (s,t)\in D,
\end{equation}
and therefore (\ref{227}) gives, for $\zeta = (s,t)$
\begin{equation}\label{325}
 \big[l(P^{H}_{K})\cdot r(\dfrac{\partial}{\partial s}),l(P^{H}_{K})\cdot r(\dfrac{\partial}{\partial t})\big]_{\beta}(s,t) = 
 Q^{H_{\beta(s,t)}}_{K_{\beta(s,t)}}
 C^{H}_{\beta(s,t)}( P^{K_{\beta(s,t)}}_{H_{\beta(s,t)}} \dfrac{\partial\beta}{\partial s}(s,t)\wedge 
 \dfrac{\partial\beta}{\partial t}(s,t)).
\end{equation}
Now we apply the preceding Proposition 7 for $\eta = D\times M\,,\;B = \{(s_{0},t)\rvert\,t\in I\}\subset D$
\begin{equation}\label{326}
 X^{D}_{\beta}(s,t) = P^{H_{\beta(s,t)}}_{K_{\beta(s,t)}} \dfrac{\partial\beta}{\partial t}(s,t)\,,\; X^{B}_{\gamma}(t) = 
 P^{H_{\gamma(t)}}_{K_{\gamma(t)}}\dot{\gamma}(t). 
\end{equation}
Then (\ref{315}) gives
\begin{equation}\label{327}
 <(l(P^{H}_{K})\cdot r(\dfrac{\partial}{\partial t}))^{\prime}(\beta);Z>(s_{0},t) =  
 <(l(P^{H}_{K})\cdot r(\dfrac{\partial}{\partial t}))^{\prime}(\beta(s_{0},\cdot));Z\rvert_{s=s_{0}}>(t).
\end{equation}
According to Theorem 7 
\begin{equation}\label{328}
<(l(P^{H}_{K})\cdot r(\dfrac{\partial}{\partial t}))^{\prime}(\beta(s_{0},\cdot));Z\rvert_{s=s_{0}}>(t) = 
\nabla_{\frac{\partial}{\partial t}}^{\beta(s_{0},\cdot)^{\ast}K}Z_{(s_{0},t)}. 
\end{equation}
We take here $Z_{(s,t)} =  P^{H_{\beta(s,t)}}_{K_{\beta(s,t)}}\dfrac{\partial\beta}{\partial s}(s,t)$ and then (\ref{322}) comes from (\ref{325}),
(\ref{327}) and the definition (\ref{273}). Finally, (\ref{322}) and (\ref{323}) are equivalent in virtue of the equality
\begin{equation}\label{329}
\dfrac{\partial}{\partial t}(\gamma^{t_{0}}_{t}\alpha_{\gamma(t)})= \gamma^{t_{0}}_{t} \nabla_{\dot{\gamma}(t)}\alpha_{\gamma(t)} 
\end{equation}
which holds for every linear connection and section $\alpha$ (see also (\ref{308})) $\blacksquare$\\

Coming back to the notations (\ref{219}) - (\ref{222}) and (\ref{21}), i.e.
\begin{equation}\label{330}
 P^{H} : V\longrightarrow V/H
\end{equation}
for the canonical projection, we recall (see, for instance, Narasimhan [4], for charts on Grassmann manifolds) the natural bijection between the 
supplementary subspaces $K$ of $H$ and the linear sections 
$S\in\mathrm{Hom}\,(V/H,V)$ of $P^{H}$ given by
\begin{equation}\label{331}
 K\widetilde{\longmapsto} S = Q^{H}_{K}\,;
\end{equation}
and also the affine structure of the space of these sections, modelled on the vector space $\mathrm{Hom}\,(V/H,H)$, where
\begin{equation}\label{332}
 \overrightarrow{K_{1}K_{2}} :=  Q^{H}_{K_{2}} - Q^{H}_{K_{1}}.
\end{equation}
Note the equality
\begin{equation}\label{333}
  Q^{H}_{K_{2}} - Q^{H}_{K_{1}} = P^{K_{1}}_{H}Q^{H}_{K_{2}}.
\end{equation}
In the framework of (\ref{319}) we will denote
\begin{equation}\label{334}
 [\beta_{0}^{\ast}K]^{t}_{s} := \big(Q^{H_{\beta_{0}(t)}}_{K_{\beta_{0}(t)}}\big)^{-1}\cdot 
 (\beta_{0}^{\ast}K)_{s}^{t}\cdot Q^{H_{\beta_{0}(s)}}_{K_{\beta_{0}(s)}} 
\end{equation}
the operators induced in quotients
\begin{equation}\label{335}
 [\beta_{0}^{\ast}K]^{t}_{s}\in\mathrm{Hom}\, (T_{\beta_{0}(s)}M/H_{\beta_{0}(s)},T_{\beta_{0}(t)}M/H_{\beta_{0}(t)})\,,\;s, t\in I.
\end{equation}
They keep the multiplicative properties
\begin{equation}\label{336}
  [\beta_{0}^{\ast}K]^{t}_{u}\cdot [\beta_{0}^{\ast}K]^{u}_{s} =  [\beta_{0}^{\ast}K]^{t}_{s}\,,\;\;\;
  [\beta_{0}^{\ast}K]^{s}_{t} = \big([\beta_{0}^{\ast}K]^{t}_{s}\big)^{-1}\,,\; s, u, t\in I.
\end{equation}
The operators $ [\beta_{0}^{\ast}K]^{t}_{s}$ represent the parallel transport of the fibers of $\beta_{0}^{\ast}(TM/H)\longrightarrow I$
determined by the linear connection obtained through the vector bundle isomorphism $Q^{H_{\beta_{0}(t)}}_{K_{\beta_{0}(t)}}\,,\;t\in I$, from the 
connection on $\beta_{0}^{\ast}K\longrightarrow I$.\\
The operators $[\beta_{0}^{\ast}K]^{t}_{s}$ still depend on $K$, except for the case when the curvature of $H$ is zero, as we can see from

\begin{theorem}. Let $K^{1}\,,\;K^{2}$ be two smooth supplementary subbundles to the same vector subbundle $H$ of $TM$ and 
$\beta_{0} : I\longrightarrow M$ smooth map such that
\begin{equation}\label{337}
 \dot{\beta_{0}}(t)\in H_{\beta_{0}(t)}\,,\;\forall t\in I.
\end{equation}
Then for any $t_{0}\in I$ fixed and $v\in T_{\beta_{0}(t_{0})}M/H_{\beta_{0}(t_{0})}$ arbitrary
\begin{equation}\label{338}
 \dfrac{\mathrm{d}}{\mathrm{d} t}\{[\beta_{0}^{\ast}K^{1}]^{t_{0}}_{t}\cdot[\beta_{0}^{\ast}K^{2}]^{t}_{t_{0}}\cdot v\} = 
 [\beta_{0}^{\ast}K^{1}]^{t_{0}}_{t}\cdot C^{H}_{\beta_{0}(t)}\big(P^{K^{1}_{\beta_{0}(t)}}_{H_{\beta_{0}(t)}}
 Q^{H_{\beta_{0}(t)}}_{K^{2}_{\beta_{0}(t)}}[\beta_{0}^{\ast}K^{2}]^{t}_{t_{0}}v\wedge \dot{\beta_{0}}(t)\big).
\end{equation}
\end{theorem}

\textbf{Proof}. Remark first that, if (\ref{338})
holds for some fixed $t_{0}$ and $t$ and all $v$, then (\ref{338}) holds for any other $t_{1}$ instead of $t_{0}$, the same $t$ and all $v\in
 T_{\beta_{0}(t_{1})}M/H_{\beta_{0}(t_{1})}$. To see this, it is enough to replace $v$ by $[\beta_{0}^{\ast}K^{2}]^{t_{0}}_{t_{1}}v$, to apply on 
both sides $[\beta_{0}^{\ast}K^{1}]^{t_{1}}_{t_{0}}$ and use the multiplicative property (\ref{336}). In this way, we may suppose in (\ref{338})
$t$ and $t_{0}$ as close to each other as we want - even equal.\\
The second step consists in proving (\ref{338}) for $K^{2}$ defined by a 
submersion $\pi : U\longrightarrow B$, $U$ neighbourhood of $\beta_{0}(t_{0})$,
\begin{equation}\label{339}
 K^{2}_{p} = \mathrm{ker}\,T_{p}\pi\,,\;\forall p\in U,
\end{equation}
(as in (\ref{292})). We consider $v\in T_{\beta_{0}(t_{0})}M/H_{\beta_{0}(t_{0})}$ given and let $\delta(s)$
be a short path, defined for $s$ around $s_{0}$, such that
\begin{equation}\label{340}
 \pi(\delta(s)) = \pi(\beta_{0}(t_{0}))\,,\;\forall s\,,\;\delta(s_{0}) = \beta_{0}(t_{0})\,,\;\dot{\delta}(s_{0}) =  
 Q^{H_{\beta_{0}(t_{0})}}_{K^{2}_{\beta_{0}(t_{0})}}v.
\end{equation}
Next, let $\gamma := \pi\circ\beta_{0}$ and 
\begin{equation}\label{341}
 \beta(s,t) = \gamma^{t}_{t_{0}}(\delta(s)),
\end{equation}
where $\gamma^{t}_{t_{0}}$ is the parallel transport of the nonlinear fibers $\pi^{-1}(\{b\})\,,\;b\in B$, along $\gamma$, given by the 
horizontal vector bundle $H_{p}\,,\;p\in U$. Then $\beta(s,t)\in\pi^{-1}(\{\gamma(t)\})$, for all $s$ and all $t$, as 
$\beta(s,t_{0})\in\pi^{-1}(\{\gamma(t_{0})\})$ for all $s$ (see (\ref{340})). Moreover
\begin{equation}\label{342}
 \beta(s_{0},t) = \beta_{0}(t)\,,\;\forall t\in I,
\end{equation}
since $\beta(s_{0},t) = \gamma^{t}_{t_{0}}(\delta(s_{0})) =  \gamma^{t}_{t_{0}}(\beta_{0}(t_{0})) = \beta_{0}(t)$, by the hypothesis (\ref{337}). 
Therefore (\ref{341}) gives a variation of $\beta_{0}$ through tangent to $H$ paths. And
\begin{equation}\label{343}
 \dfrac{\partial\beta}{\partial s}(s_{0},t) =  Q^{H_{\beta_{0}(t)}}_{K^{2}_{\beta_{0}(t)}}[\beta_{0}^{\ast}K^{2}]^{t}_{t_{0}}v,
\end{equation}
in virtue of Theorem 8, (\ref{296}), (\ref{340}) and definition (\ref{334}). We apply Theorem 9 to $H$, $K = K^{1}$ and this 
$\beta$; mutiplying both sides of (\ref{323}) by $\big( Q^{H_{\beta_{0}(t_{0})}}_{K^{1}_{\beta_{0}(t_{0})}}\big)^{-1}$ and using the equality
$P^{H}_{K^{1}}\cdot Q^{H}_{K^{2}} =  Q^{H}_{K^{1}}$, we obtain (\ref{338}) for any $K^{1}$ and $K^{2}$ of the form (\ref{339}).\\
Third, we show that if (\ref{338}) holds for $K^{1}$ and $K^{2}$, it is verified also for the transposed pair, $K^{2}$ and $K^{1}$. Indeed\\

$\dfrac{\mathrm{d}}{\mathrm{d}t}\{[\beta_{0}^{\ast}K^{2}]^{t_{0}}_{t}\cdot[\beta_{0}^{\ast}K^{1}]^{t}_{t_{0}}\cdot v\} =
\dfrac{\mathrm{d}}{\mathrm{d}t}\{[\beta_{0}^{\ast}K^{1}]^{t_{0}}_{t}\cdot[\beta_{0}^{\ast}K^{2}]^{t}_{t_{0}}\}^{-1}v =$\\

$= - [\beta_{0}^{\ast}K^{2}]^{t_{0}}_{t}\cdot[\beta_{0}^{\ast}K^{1}]^{t}_{t_{0}}\cdot
\dfrac{\mathrm{d}}{\mathrm{d}t}\{[\beta_{0}^{\ast}K^{1}]^{t_{0}}_{t}\cdot[\beta_{0}^{\ast}K^{2}]^{t}_{t_{0}}\}\cdot
[\beta_{0}^{\ast}K^{2}]^{t_{0}}_{t}\cdot[\beta_{0}^{\ast}K^{1}]^{t}_{t_{0}}v =$\\

$= - [\beta_{0}^{\ast}K^{2}]^{t_{0}}_{t}\cdot C^{H}_{\beta_{0}(t)}\big(P^{K^{1}_{\beta_{0}(t)}}_{H_{\beta_{0}(t)}}
 Q^{H_{\beta_{0}(t)}}_{K^{2}_{\beta_{0}(t)}}[\beta_{0}^{\ast}K^{1}]^{t}_{t_{0}}v\wedge \dot{\beta_{0}}(t)\big) = $\\

$= [\beta_{0}^{\ast}K^{2}]^{t_{0}}_{t}\cdot C^{H}_{\beta_{0}(t)}\big(P^{K^{2}_{\beta_{0}(t)}}_{H_{\beta_{0}(t)}}
 Q^{H_{\beta_{0}(t)}}_{K^{1}_{\beta_{0}(t)}}[\beta_{0}^{\ast}K^{1}]^{t}_{t_{0}}v\wedge \dot{\beta_{0}}(t)\big)$,\\
since from (\ref{333}) we get $P^{K_{2}}_{H}Q^{H}_{K_{1}} = - P^{K_{1}}_{H}Q^{H}_{K_{2}}$. Then the equality (\ref{338}) holds also for $K^{1}$
of the same form (\ref{339}) and $K^{2}$ arbitrary.\\
Finaly, we show that, if (\ref{338}) holds for the pairs $(K^{1},K^{3})$ and $(K^{3},K^{2})$, then it is verfied by the pair $(K^{1},K^{2})$.
Indeed,\\

$\dfrac{\mathrm{d}}{\mathrm{d}t}\{[\beta_{0}^{\ast}K^{1}]^{t_{0}}_{t}\cdot [\beta_{0}^{\ast}K^{2}]^{t}_{t_{0}}\cdot v\} = 
\dfrac{\mathrm{d}}{\mathrm{d}t}\{[\beta_{0}^{\ast}K^{1}]^{t_{0}}_{t}\cdot [\beta_{0}^{\ast}K^{3}]^{t}_{t_{0}}\cdot 
[\beta_{0}^{\ast}K^{3}]^{t_{0}}_{t}\cdot [\beta_{0}^{\ast}K^{2}]^{t}_{t_{0}}\cdot v\} =$\\

$=\dfrac{\mathrm{d}}{\mathrm{d}t}\{[\beta_{0}^{\ast}K^{1}]^{t_{0}}_{t}\cdot [\beta_{0}^{\ast}K^{3}]^{t}_{t_{0}}\}\cdot 
[\beta_{0}^{\ast}K^{3}]^{t_{0}}_{t}\cdot [\beta_{0}^{\ast}K^{2}]^{t}_{t_{0}}\cdot v +$\\

$+ [\beta_{0}^{\ast}K^{1}]^{t_{0}}_{t}\cdot [\beta_{0}^{\ast}K^{3}]^{t}_{t_{0}}\cdot
\dfrac{\mathrm{d}}{\mathrm{d}t}\{[\beta_{0}^{\ast}K^{3}]^{t_{0}}_{t}\cdot [\beta_{0}^{\ast}K^{2}]^{t}_{t_{0}}\cdot v\} =$\\

$= [\beta_{0}^{\ast}K^{1}]^{t_{0}}_{t}\cdot C^{H}_{\beta_{0}(t)}\big(P^{K^{1}_{\beta_{0}(t)}}_{H_{\beta_{0}(t)}}
 Q^{H_{\beta_{0}(t)}}_{K^{3}_{\beta_{0}(t)}}[\beta_{0}^{\ast}K^{2}]^{t}_{t_{0}}v\wedge \dot{\beta_{0}}(t)\big) +$

$+ [\beta_{0}^{\ast}K^{1}]^{t_{0}}_{t}\cdot C^{H}_{\beta_{0}(t)}\big(P^{K^{3}_{\beta_{0}(t)}}_{H_{\beta_{0}(t)}}
 Q^{H_{\beta_{0}(t)}}_{K^{2}_{\beta_{0}(t)}}[\beta_{0}^{\ast}K^{2}]^{t}_{t_{0}}v\wedge \dot{\beta_{0}}(t)\big) =$\\
 
$= [\beta_{0}^{\ast}K^{1}]^{t_{0}}_{t}\cdot C^{H}_{\beta_{0}(t)}\big(P^{K^{1}_{\beta_{0}(t)}}_{H_{\beta_{0}(t)}}
 Q^{H_{\beta_{0}(t)}}_{K^{2}_{\beta_{0}(t)}}[\beta_{0}^{\ast}K^{2}]^{t}_{t_{0}}v\wedge \dot{\beta_{0}}(t)\big),$ \\

as (\ref{333}) gives $ P^{K_{1}}_{H}Q^{H}_{K_{3}} +  P^{K_{3}}_{H}Q^{H}_{K_{2}} =  P^{K_{1}}_{H}Q^{H}_{K_{2}}$. The proof ends by combining these 
facts $\blacksquare$\\
The equation (\ref{338}) allows to determine the parallel transport of the connection induced on $\beta_{0}^{\ast}K^{2}\longrightarrow I$,
when the operators $(\beta_{0}^{\ast}K^{1})_{s}^{t}\,$, corresponding to another supplementary vector subbundle $K^{1}$, are known.

\section{The equation of infinitesimal variation through\\
tangent paths to the given vector subbundle}

First, the following consequence of the previous relation (\ref{338}).
\begin{theorem}. In the same hypothesis (\ref{337}), the operator
\begin{equation}\label{344}
 (\textsl{D} X)_{t}  =: [\beta_{0}^{\ast}K]_{t_{0}}^{t}\dfrac{\mathrm{d}}{\mathrm{d}t}\big([\beta_{0}^{\ast}K]_{t}^{t_{0}}
 P^{H_{\beta_{0}(t)}} X_{t}\big) - C^{H}_{\beta_{0}(t)}\big(P^{K_{\beta_{0}(t)}}_{H_{\beta_{0}(t)}} X_{t}\wedge\dot{\beta_{0}}(t)\big)\,,
\end{equation}
acting as
\begin{equation}\label{345}
 \textsl{D} : \textsl{C}^{\infty}\Gamma(\beta_{0}^{\ast}(TM))\longrightarrow  \textsl{C}^{\infty}\Gamma(\beta_{0}^{\ast}(TM/H))\,,
\end{equation}
depends neither on $t_{0}$, nor on the smooth, supplementary to $H$, vector subbundle $K$.
\end{theorem}

\textbf{Proof}. Let $K^{1}, K^{2}$ be two such subbundles of $TM$; then the difference of the expressions (\ref{344}) corresponding to them
gives

$[\beta_{0}^{\ast}K^{1}]_{t_{0}}^{t}\dfrac{\mathrm{d}}{\mathrm{d}t}\big([\beta_{0}^{\ast}K^{1}]_{t}^{t_{0}}
 P^{H_{\beta_{0}(t)}} X_{t}\big) - C^{H}_{\beta_{0}(t)}\big(P^{K^{1}_{\beta_{0}(t)}}_{H_{\beta_{0}(t)}} X_{t}\wedge\dot{\beta_{0}}(t)\big)
 - [\beta_{0}^{\ast}K^{2}]_{t_{0}}^{t}\dfrac{\mathrm{d}}{\mathrm{d}t}\big([\beta_{0}^{\ast}K^{2}]_{t}^{t_{0}}
 P^{H_{\beta_{0}(t)}} X_{t}\big) +$\\

$+ C^{H}_{\beta_{0}(t)}\big(P^{K^{2}_{\beta_{0}(t)}}_{H_{\beta_{0}(t)}} X_{t}\wedge\dot{\beta_{0}}(t)\big) = 
 [\beta_{0}^{\ast}K^{1}]^{t}_{t_{0}}\cdot
\dfrac{\mathrm{d}}{\mathrm{d}t}\{[\beta_{0}^{\ast}K^{1}]^{t_{0}}_{t}\cdot[\beta_{0}^{\ast}K^{2}]^{t}_{t_{0}}\}\cdot
[\beta_{0}^{\ast}K^{2}]^{t_{0}}_{t} P^{H_{\beta_{0}(t)}} X_{t} +$\\

$+ C^{H}_{\beta_{0}(t)}\big((P^{K^{2}_{\beta_{0}(t)}}_{H_{\beta_{0}(t)}} - 
P^{K^{1}_{\beta_{0}(t)}}_{H_{\beta_{0}(t)}}) X_{t}\wedge\dot{\beta_{0}}(t)\big) = 
C^{H}_{\beta_{0}(t)}\big(P^{K^{1}_{\beta_{0}(t)}}_{H_{\beta_{0}(t)}} Q^{H_{\beta_{0}(t)}}_{K^{2}_{\beta_{0}(t)}} P^{H_{\beta_{0}(t)}} X_{t}
\wedge\dot{\beta_{0}}(t)\big) +$

$+ C^{H}_{\beta_{0}(t)}\big((P_{K^{1}_{\beta_{0}(t)}}^{H_{\beta_{0}(t)}} - 
P_{K^{2}_{\beta_{0}(t)}}^{H_{\beta_{0}(t)}}) X_{t}\wedge\dot{\beta_{0}}(t)\big) =
C^{H}_{\beta_{0}(t)}\big(( Q^{H_{\beta_{0}(t)}}_{K^{2}_{\beta_{0}(t)}} - 
Q^{H_{\beta_{0}(t)}}_{K^{1}_{\beta_{0}(t)}}) P^{H_{\beta_{0}(t)}} X_{t}
\wedge\dot{\beta_{0}}(t)\big) +$

$+  C^{H}_{\beta_{0}(t)}\big((P_{K^{1}_{\beta_{0}(t)}}^{H_{\beta_{0}(t)}} - 
P_{K^{2}_{\beta_{0}(t)}}^{H_{\beta_{0}(t)}}) X_{t}\wedge\dot{\beta_{0}}(t)\big)$,\\

where we have used (\ref{338}) and also (\ref{333}). But this makes zero, as
\begin{equation}\label{346}
 P^{H}_{K} = Q^{H}_{K}\cdot P^{H}
\end{equation}
(see also (\ref{226}))$\blacksquare$\\
Next we prove the reciprocal of Theorem 9.

\begin{theorem}. If $\beta_{0} : I\longrightarrow M$ is tangent to $H$, then $X\in\textsl{C}^{\infty}\Gamma(\beta_{0}^{\ast}(TM))$ satisfies
\begin{equation}\label{347}
 (\textsl{D} X)_{t} = 0\,,\,\forall t\in I,
\end{equation}
if and only if there exists a smooth variation of $\beta_{0}$ on $I$
\begin{equation}\label{348}
 \beta : (-\varepsilon,\varepsilon)\times I\longrightarrow M\,,\;\beta(0,\cdot) = \beta_{0}\,,
\end{equation}
in tangent to $H$ paths, i.e.
\begin{equation}\label{349}
 \dfrac{\partial\beta}{\partial t}(s,t)\in H_{\beta(s,t)}\,,\,\forall t\in I\,,\;\forall s\in(-\varepsilon,\varepsilon)\,,
\end{equation}
such that
\begin{equation}\label{350}
 X_{t} =  \dfrac{\partial\beta}{\partial s}(0,t)\,,\;\forall t\in I.
\end{equation}
\end{theorem}

\textbf{Proof}. it is easy to verify that (\ref{323}) can be written as $\textsl{D}\,\dfrac{\partial\beta}{\partial s}(0,\cdot) = 0$ on $I$.\\
Conversely, $X$ being given with (\ref{347}), we have to construct $\beta$ with (\ref{348}), (\ref{349}) and (\ref{350}).
First, we consider the case when
\begin{equation}\label{351}
 \beta_{0}(I)\subset U\,,\;\chi : U\longrightarrow V\times W\,,\;\chi (U) = D\times G,
\end{equation}
$\chi$ being a chart of the form (\ref{5}) where $H$ is of the form (\ref{6}) ($D$ open in $V$, $G$ open in $W$). Then, disregarding the 
chart $\chi$, we consider $K$ and $H$ as in (\ref{303}), i.e.
\begin{equation}\label{352}
K = \{(x,y;0_{V},w)\rvert\, x\in D, y\in G, w\in W\},\;
H = \{(x,y;v,\varGamma(x,y)v)\rvert\,x\in D, y\in G, v\in V\}.
\end{equation}
Let us denote
\begin{equation}\label{353}
 \beta_{0}(t) = (x_{0}(t),y_{0}(t))\in V\times W\,,\;X(t) = (A(t),B(t))\in V\times W;
\end{equation}
then the equation (\ref{347}) becomes
\begin{multline}\label{354}
 \dfrac{\mathrm{d}B}{\mathrm{d}t}(t) - \varGamma(x_{0}(t),y_{0}(t)) \dfrac{\mathrm{d}A}{\mathrm{d}t}(t) - 
 <\dfrac{\partial\varGamma}{\partial x}(x_{0}(t),y_{0}(t));A(t)>\dfrac{\mathrm{d}x_{0}}{\mathrm{d}t}(t) -\\
 - <\dfrac{\partial\varGamma}{\partial y}(x_{0}(t),y_{0}(t));B(t)>\dfrac{\mathrm{d}x_{0}}{\mathrm{d}t}(t) = 0.  
\end{multline}
If, taking into account (\ref{348}), we denote
\begin{equation}\label{355}
 \beta(s,t) = (x(s,t),y(s,t))\,,\;x(0,t) = x_{0}(t)\,,\;y(0,t) = y_{0}(t),
\end{equation}
the condition (\ref{349}) reads
\begin{equation}\label{356}
 \dfrac{\partial y}{\partial t}(s,t) = \varGamma (x(s,t),y(s,t))\dfrac{\partial x}{\partial t}(s,t).
\end{equation}
It is clear that (\ref{354}) is the consequence of (\ref{356}) when (\ref{350}) is satisfied, i.e.
\begin{equation}\label{357}
 \dfrac{\partial x}{\partial s}(0,t) = A(t)\,,\;\dfrac{\partial y}{\partial s}(0,t) = B(t).
\end{equation}
So, taking (\ref{354}) as hypothesis, together with 
\begin{equation}\label{358}
  \dfrac{\mathrm{d}y_{0}}{\mathrm{d}t}(t) = \varGamma(x_{0}(t),y_{0}(t)) \dfrac{\mathrm{d}x_{0}}{\mathrm{d}t}(t),
\end{equation}
we look for $x(s,t), y(s,t)$ with (\ref{356}), (\ref{355}) and (\ref{357}). $x(s,t)$ can be chosen arbitrary with (\ref{355}) and (\ref{357}).
Next, for $t_{0}\in I$ fixed and $y(s,t_{0})$ arbitrary such that
\begin{equation}\label{359}
 y(0,t_{0}) = y_{0}(t_{0})\,,\;\dfrac{\partial y}{\partial s}(0,t_{0}) = B(t_{0}),
\end{equation}
we consider $y(s,t)$ defined by (\ref{356}) and $y(s,t_{0})$ so prescribed. Then $y(s,t)$ will satisfy $y(0,t) = y_{0}(t)$ and 
$\dfrac{\partial y}{\partial s}(0,t) = B(t)\,,\;\forall t\in I$, since both $\dfrac{\partial y}{\partial s}(0,t)$ and $ B(t)$ satisfy as $Z(t)$
the differential equation\\

$\dfrac{\mathrm{d}Z}{\mathrm{d}t}(t) - <\dfrac{\partial\varGamma}{\partial y}(x_{0}(t),y_{0}(t));Z(t)>\dfrac{\mathrm{d}x_{0}}{\mathrm{d}t}(t) - 
\varGamma(x_{0}(t),y_{0}(t)) \dfrac{\mathrm{d}A}{\mathrm{d}t}(t) -\\
- <\dfrac{\partial\varGamma}{\partial x}(x_{0}(t),y_{0}(t));A(t)>\dfrac{\mathrm{d}x_{0}}{\mathrm{d}t}(t) = 0$\\

with the same initial condition $Z(t_{0}) = B(t_{0})$.\\
It remains however to show that $\exists\, \varepsilon > 0$ such that the solution of the equation (\ref{356}) be defined on all of $I$ for 
$\lvert s\rvert <\varepsilon$. We have in mind to use the following fact: \textit{for a smooth vector field $X$, on a manifold $M$, of local flow 
$\varphi^{t}(x),\;x\in M$ and $t$ in a neighbourhood $J(x)$ of $0$ in $\textbf{R}$, $\varphi^{0}(x) = x$, if $\varphi^{t}(x_{0})$ is defined 
for $t$ in a compact interval $I$ containing $0$, there exists a neighbourhood $U$ of $x_{0}$ in $M$ such that $\varphi^{t}(x)$ is defined on all of $I$ 
for $x\in U$.} (For $M$ compact, $\varphi^{t}(x)$ is defined $\forall\, t\in\textbf{R}$ and $\forall\, x\in M$. If $M$ is not compact we consider
$\psi\in \textsl{C}^{\infty}_{0}(M,\textbf{R})\,,\;\psi = 1$ in a neighbourhood of $\{\varphi^{t}(x_{0})\,\rvert\, t\in I\}$. Then the
flow of $\psi X$ is global and coincides with the flow of $X$ on $I$ for $x$ in a neighbourhood of $x_{0}$.) We write the equation 
(\ref{356}), with $y(s,t_{0})$ given, in the form

$\dfrac{\mathrm{d}z}{\mathrm{d}t}(t) = \varGamma(x(\sigma(t),\tau(t)),z(t))\dfrac{\partial x}{\partial\tau}(\sigma(t),\tau(t))\,,\;
\dfrac{\mathrm{d}\sigma}{\mathrm{d}t}(t) = 0\,,\;\dfrac{\mathrm{d}\tau}{\mathrm{d}t}(t) = 1;$

$z(0) = y(s,t_{0})\,,\;\sigma(0) = s\,,\,\tau(0) = t_{0}$.

Then, of course, $\sigma(t) = s\,,\,\tau(t) = t+t_{0}\,,\;z(t) = y(s,t+t_{0})\,,\;\forall t\in I-t_{0}$, and the solution will be defined on all 
of $I-t_{0}$ for $s$ in a neighbourhood of $0$.\\
In the general case, in order to find $\varepsilon > 0$ such that $\beta$ is defined as in (\ref{348}), with (\ref{349}) and (\ref{350}),
satisfied on $I$, we proceed as follows. For all $t_{0}$ we find $J\ni t_{0}$ interval, open as a subset of the compact interval $I$, such that
$\overline{J}$ has the properties (\ref{351}) of the interval $I$. We then extract a finite and minimal subcovering for $I$ with such intervals
$J$ (that can be of one of the forms: $J = [a,t_{2})\,,\;J = (t_{1},t_{2})\,,\;J = (t_{1},b]\,,\; J = [a,b]$, if $I = [a,b]$). When the covering 
is minimal, the maps $J\longmapsto\mathrm{inf}\,J$ and $J\longmapsto\mathrm{sup}\,J$ are injective. Let us enumerate the intervals of the 
covering such that
\begin{equation}\label{360}
 \mathrm{sup}\,J_{k}<\mathrm{sup}\,J_{k+1}\,,\;1\leqslant k< n,
 \end{equation}
if $n> 1$ is their number ($n=1$ corresponds to the previous situation). It is easy to see that
\begin{equation}\label{361}
 \mathrm{inf}\,J_{1} = a\,,\;\mathrm{sup}\,J_{n} = b.
\end{equation}
Also, to check the inequality 
\begin{equation}\label{362}
 \mathrm{inf}\,J_{k+1} < \mathrm{sup}\,J_{k}\,,\;1\leqslant k< n,
\end{equation}
in the contrary case the points from $[\mathrm{sup}\,J_{k},\mathrm{inf}\,J_{k+1}]$ being not covered.\\
Then the solution $\beta(s,t)$ can be constructed recurrently: first on $\overline{J_{1}}$ for $\lvert s \rvert < \varepsilon_{1}$; if the 
solution is already defined on $\overline{\bigcup_{k=1}^{q}J_{k}} = [a,\mathrm{sup}\,J_{q}]$, for $\lvert s \rvert < \varepsilon_{q}$, it can be 
extended to $\overline{\bigcup_{k=1}^{q+1}J_{k}}$, since according to (\ref{362}) $\mathrm{inf}\,J_{q+1} < \mathrm{sup}\,J_{q}$ and on the 
interval $[\mathrm{inf}\,J_{q+1},\mathrm{sup}\,J_{q}]$ the solution is already constructed using the chart $\chi_{q}$. With $x(s,t)$ suitably
extended, $y(s,\frac{1}{2}(\mathrm{inf}\,J_{q+1} + \mathrm{sup}\,J_{q}))$ inherited and the chart $\chi_{q+1}$, we get the extension of the 
solution to $\overline{\bigcup_{k=1}^{q+1}J_{k}} = [a,\mathrm{sup}\,J_{q+1}]$, for a certain $\varepsilon_{q+1}\leqslant\varepsilon_{q}$ and 
$\lvert s \rvert < \varepsilon_{q+1}$; and so on, up to $q+1 = n$. The theorem is proven $\blacksquare$\\
Of course, the equation (\ref{347}) can be written in the form (see (\ref{322})) 
\begin{equation}\label{363}
 \nabla_{\frac{\partial}{\partial t}}^{\beta_{0}^{\ast}K} P^{H_{\beta_{0}(t)}}_{K_{\beta_{0}(t)}}X_{t} 
 = Q^{H_{\beta_{0}(t)}}_{K_{\beta_{0}(t)}}
 C^{H}_{\beta_{0}(t)}( P^{K_{\beta_{0}(t)}}_{H_{\beta_{0}(t)}} X_{t}\wedge 
 \dot{\beta}_{0}(t)).
\end{equation}
The equation (\ref{347}) appears to be the root for the \textit{Jacobi equation of infinitesimal variation of geodesics}. In that case we have 
a linear connection on the vector bundle $TM$ over a smooth manifold $M$. In this analysis we will make use of
\begin{proposition}. For every smooth path
\begin{equation}\label{364}
 Z : I\longrightarrow TM\,,
\end{equation}
on an interval $I\subseteq\textbf{R}$, if we denote
\begin{equation}\label{365}
 \gamma := \tau_{M}\circ Z,
\end{equation}
for $X\in TM$, $V_{X} : T_{X}(TM)\longrightarrow T_{X}(T_{\tau_{M}(X)}M)$ the vertical projection defining the connection (see (\ref{266})) and
$T$ the respective torsion tensor, the following equality holds $\forall t\in I$
\begin{equation}\label{366}
 \varPsi^{T_{\gamma(t)}M}(\dot{\gamma}(t),\cdot)^{-1}\; V_{\dot{\gamma}(t)}\;\omega_{M}(\dot{Z}(t)) - \varPsi^{T_{\gamma(t)}M}(Z(t),\cdot)^{-1}
 \;V_{Z(t)}\;\dot{Z}(t) = T_{\gamma(t)}(Z(t),\dot{\gamma}(t))\,.
\end{equation}
\end{proposition}
\textbf{Proof}. As this equality (\ref{366}) is local and of intrinsic meaning, we may suppose $M = U$ open subset in the vector space $V$, so 
that $TM = U\times V$. In that case \\
$H_{(x,y)}\subset T_{(x,y)}(U\times V) = \{(x,y)\}\times V\times V$,
\begin{equation}\label{367}
 H_{(x,y)} = \{(x,y;v,\varGamma(x,y)\,v)\rvert\,v\in V\}
\end{equation}
where the liniarity of the connection corresponds to the linearity of $\varGamma(x,y)$ in $y$, or to the existence of linear $\varGamma(x)$ on 
$V\otimes V$ such that
\begin{equation}\label{368}
 \varGamma(x,y)\,v = \varGamma(x)(y\otimes v)\,,\; x\in U\,,\;y, v\in V.
\end{equation}
Then
\begin{equation}\label{369}
 V_{(x,y)}(x,y;v,w) =  (x,y;0_{V},w - \varGamma(x)(y\otimes v)).
\end{equation}
In these local coordinates the torsion tensor becomes (see Kobayashi + Nomizu, vol. I [2])
\begin{equation}\label{370}
 T_{x}(v,w) = \varGamma(x)(v\otimes w - w\otimes v).
\end{equation}
For (see (\ref{364}) and (\ref{365}))
\begin{equation}\label{371}
 Z(t) = (x(t),y(t))\,,\;\gamma(t) = x(t),
\end{equation}
we have
\begin{multline}\label{372}
 \dot{\gamma}(t) = (x(t);x^{\prime}(t))\,,\;\dot{Z}(t) = (x(t),y(t);x^{\prime}(t),y^{\prime}(t))\,,\\
 \omega_{M}(\dot{Z}(t)) =  (x(t),x^{\prime}(t);y(t),y^{\prime}(t)).
\end{multline}
Using (\ref{369}) we get\\

$ \varPsi^{T_{\gamma(t)}M}(\dot{\gamma}(t),\cdot)^{-1}\; V_{\dot{\gamma}(t)}\;\omega_{M}(\dot{Z}(t)) = 
(x(t);y^{\prime}(t) - \varGamma(x(t))(x^{\prime}(t)\otimes y(t)))$,

$\varPsi^{T_{\gamma(t)}M}(Z(t),\cdot)^{-1}\;V_{Z(t)}\;\dot{Z}(t) = 
(x(t);y^{\prime}(t) - \varGamma(x(t))(y(t)\otimes x^{\prime}(t)))$,

wherefrom the result $\blacksquare$\\

It is easy to establish \textit{the relation between the Riemann tensor $R$ and the curvature tensor $C$ of the horizontal subbundle $H$} of 
$T(TM)$ defining the linear connection, namely
\begin{equation}\label{373}
 R_{\tau_{M}(X)}(T_{X}\tau_{M}\;A,T_{X}\tau_{M}\;B)\;X = -\, \varPsi^{T_{\tau_{M}(X)}M}(X,\cdot)^{-1}\; Q_{K_{X}}^{H_{X}}\; 
 C_{X}(A\wedge B),
\end{equation}
if $X\in TM\,,\;A, B\in H_{X}$ and 
\begin{equation}\label{374}
 K_{X} = T_{X}(T_{\tau_{M}(X)}M).
\end{equation}
Recall that $T_{X}\tau_{M} : H_{X}\widetilde{\longrightarrow} T_{\tau_{M}(X)}M\,,\;
\varPsi^{T_{\tau_{M}(X)}M}(X,\cdot) :  T_{\tau_{M}(X)}M\widetilde{\longrightarrow} K_{X}$ and \\
$ Q_{K_{X}}^{H_{X}} : T_{X}(TM)/H_{X}\widetilde{\longrightarrow} K_{X}$ are isomorphisms.\\

And the link with the Jacobi equation (see Kobayashi + Nomizu vol. II [2]) is given in
\begin{theorem}. Let $M$ be smooth manifold with a linear connection on $TM$ given by $\nabla$. Let $\gamma_{0} : I\longrightarrow M$ be a 
geodesic for $\nabla$. Then $Z\in\textsl{C}^{\infty}\Gamma(\gamma_{0}^{\ast}(TM)) = T_{\gamma_{0}}\textsl{C}^{\infty}(I,M)$ is a solution for the 
Jacobi equation 
\begin{equation}\label{375}
 \nabla_{\dot{\gamma}_{0}(t)}^{2}Z(t) +  \nabla_{\dot{\gamma}_{0}(t)}(T(Z(t),\dot{\gamma}_{0}(t))) + 
 R(Z(t),\dot{\gamma}_{0}(t))\;\dot{\gamma}_{0}(t) = 0
\end{equation}
if and only if
\begin{equation}\label{376}
 X_{t} := \omega_{M}(\dot{Z}(t))
\end{equation}
satisfies the equation (\ref{347}) for $H$ the horizontal subbundle of $T(TM)$ corresponding to $\nabla$ and
\begin{equation}\label{377}
 \beta_{0} = \dot{\gamma}_{0}.
\end{equation}
\end{theorem}

\textbf{Proof}. From (\ref{297}), for $ P^{H_{\beta_{0}(t)}}_{K_{\beta_{0}(t)}}X_{t} = \varPsi^{E_{\gamma(t)}}(\beta_{0}(t),\alpha_{\gamma(t)})$
and $E = TM$, we obtain \\
$\varPsi^{T_{\gamma_{0}(t)}M}(\beta_{0}(t),\cdot)^{-1}\;\nabla_{\frac{\partial}{\partial t}}^{\beta_{0}^{\ast}K} 
P^{H_{\beta_{0}(t)}}_{K_{\beta_{0}(t)}}X_{t} = \nabla_{\dot{\gamma}_{0}(t)}\varPsi^{T_{\gamma_{0}(t)}M}(\beta_{0}(t),\cdot)^{-1}
P^{H_{\beta_{0}(t)}}_{K_{\beta_{0}(t)}}X_{t}$\\
Next, using (\ref{363}) and (\ref{373}) we infer
\begin{equation}\label{378}
 \nabla_{\dot{\gamma}_{0}(t)}\varPsi^{T_{\gamma_{0}(t)}M}(\beta_{0}(t),\cdot)^{-1}
P^{H_{\beta_{0}(t)}}_{K_{\beta_{0}(t)}}X_{t} = -  R(T_{\beta_{0}(t)}\tau_{M}\;P^{H_{\beta_{0}(t)}}_{K_{\beta_{0}(t)}}X_{t}\,,\dot{\gamma}_{0}(t))\;
\dot{\gamma}_{0}(t).
\end{equation}
But, for $X_{t}$ from (\ref{376}), $T_{\beta_{0}(t)}\tau_{M}\;P^{H_{\beta_{0}(t)}}_{K_{\beta_{0}(t)}}X_{t} = Z(t)$, since 
$T_{\beta_{0}(t)}\tau_{M}\;P^{H_{\beta_{0}(t)}}_{K_{\beta_{0}(t)}}X_{t} = T_{\beta_{0}(t)}\tau_{M}\;X_{t} =\\
= (T\tau_{M}\circ\omega_{M})(\dot{Z}(t)) = \tau_{TM}(\dot{Z}(t)) = Z(t)$. As in our case $V_{X} = P^{H_{X}}_{K_{X}}$, in the left hand side of 
(\ref{378}) we can use (\ref{366}) from Proposition 8 and the usual notation
\begin{equation}\label{379}
 \nabla_{\dot{\gamma}_{0}(t)}Z(t) = \varPsi^{T_{\gamma_{0}(t)}M}(Z(t),\cdot)^{-1}\;V_{Z(t)}\;\dot{Z}(t).
\end{equation}
In this way we get\\
$\nabla_{\dot{\gamma}_{0}(t)}(\nabla_{\dot{\gamma}_{0}(t)}Z(t) + T_{\gamma_{0}(t)}(Z(t),\dot{\gamma}_{0}(t))) = 
- R(Z(t),\dot{\gamma}_{0}(t))\;\dot{\gamma}_{0}(t)$,\\
or the equation (\ref{375}) $\blacksquare$

\end{document}